\documentstyle[12pt,leqno]{article}
\textheight9in \def\DATE{July 8, 1997}
\newtheorem{theorem}{Theorem}[section]
\newtheorem{definition}[theorem]{Definition}

\newtheorem{observation}[theorem]{Observation}

\newtheorem{example}[theorem]{Example}
\newtheorem{lemma}[theorem]{Lemma}

\newtheorem{proposition}[theorem]{Proposition}
\newtheorem{remark}[theorem]{Remark}

\catcode`\@=11

\def\ps@myheadings{\let\@mkboth\@gobbletwo
\def\@oddhead{\ifnum\count0=1 \hfill\else
\rightmark \hfil \rm\thepage\fi}%
\def\@oddfoot{\ifnum\count0=1 \hfill \rm 1 \hfill \else
\hfill\fi}
\def\@evenhead%
{\rm\leftmark\hfil\rm\thepage}%
\def\@evenfoot{}\def\sectionmark##1{}
\def\subsectionmark##1{}}

\def\@begintheorem#1#2{\it \trivlist \item[\hskip
 \labelsep{\bf #1\ #2.}]}
\def\@opargbegintheorem#1#2#3{\it \trivlist\item[\hskip%
 \labelsep{\bf #1\ #2.\ (#3)}]}
\def\@endtheorem{\endtrivlist}

\def\@listI{\leftmargin\leftmargini \parsep 1pt plus 2.5pt
 minus 1pt\topsep 5pt plus 2pt minus 3pt%
 \itemsep 0pt plus 2.5pt minus 1pt}
\let\@listi\@listI
\@listi

\def\@sect#1#2#3#4#5#6[#7]#8{\ifnum #2>\c@secnumdepth%
 \def \@svsec {}\else \refstepcounter {#1}\edef \@svsec%
 {\csname the#1\endcsname. \hskip .1em }\fi \@tempskipa%
 #5\relax \ifdim \@tempskipa >\z@ \begingroup #6\relax%
 \@hangfrom {\hskip #3\relax \@svsec }{\interlinepenalty%
 \@M #8.\par }\endgroup \csname #1mark\endcsname {#7}%
 \addcontentsline {toc}{#1}{\ifnum #2>\c@secnumdepth%
 \else \protect \numberline {\csname the#1\endcsname. }%
 \fi #7}\else \def \@svsechd {#6\hskip #3\@svsec #8.%
 \csname #1mark\endcsname {#7}\addcontentsline {toc}{#1}%
 {\ifnum #2>\c@secnumdepth \else \protect \numberline%
 {\csname the#1\endcsname. }\fi #7}}\fi \@xsect {#5}}

\def\section{\@startsection {section}{1}{\z@ }%
 {-3.5ex plus -1ex minus -.2ex}{2.3ex plus .2ex}{\bf }}

\def\thebibliography#1{%
 \section *{References.\@mkboth {REFERENCES}{REFERENCES}}%
 \list {[\arabic {enumi}]}{\settowidth \labelwidth {[#1]}%
 \leftmargin \labelwidth \advance \leftmargin \labelsep %
 \usecounter {enumi}} \def \newblock %
 {\hskip .11em plus .33em minus -.07em} \sloppy \clubpenalty
 4000%
 \widowpenalty 4000 \sfcode`\.=1000\relax}

\def\@maketitle{%
 \newpage \null \vskip 2em
 \begin{center}{\Large\bf \@title \par }
 \vskip 1.5em
 {\large \lineskip .5em
 \begin {tabular}[t]{c}\@author
 \end{tabular}\par } \vskip .8em {June 20, 1994}
 \end{center}\par \vskip 1.5em}

\def\footnote{\@ifnextchar [{\@xfootnote }{\stepcounter
{\@mpfn }%
\begingroup \let \protect
\noexpand \xdef \@thefnmark {\hskip-3mm}%
\endgroup \@footnotetext
}}
\def\@makefnmark{}

\def\abstract{%
\if@twocolumn \section *{Abstract}
 \else \small\quotation\noindent{\bf Abstract.}\fi}

\def\fnum@figure{{\bf \figurename {} \thefigure }}

\catcode`\@=13

\def\qed{\hspace*{\fill}
\mbox{\hphantom{mm}\rule{0.25cm}{0.25cm}}\\}

\def\B{{\cal B}}
\def\I{{\cal I}}
\def\bfR{{\bf R}}
\def\BC{{{\cal B}C}}
\def\rada#1#2{#1,\ldots,#2}
\def\IC{{{\cal I}C}}
\def\C{{\cal C}}
\def\P{{\cal P}}
\def\coll#1{\{{#1}(n)\}_{n\geq 1}}
\def\colla#1{\{{#1}_n\}_{n\geq 1}}
\def\E{{\cal E}}
\def\Hom#1#2{{\rm Hom}(#1,#2)}
\def\id{\mbox{$1\!\!1$}}

\def\barB{{\overline {\cal B}}}
\def\barK{{\overline K}}
\def\barA{{\overline {\cal A}}}
\def\Ainfty{{A($\infty$)}}
\def\barBC{{\overline {{\cal B}C}}}
\def\barW{{\overline W}}
\def\M{{\cal M}}
\def\ext{\mbox{\large$\land$}}

\def\squeeze{\times\hskip-1.5mm \cdot \hskip-1mm%
 \cdot\hskip-1mm\cdot\hskip-1.5mm\times}
\def\Z{{\bf Z}}
\def\fr{{\cal F}}
\def\tr{{\cal T}}
\def\tT{{\tt T}}
\def\cyclsum{
{\raisebox{-.4mm}{\Large $\circ$}%
 \hskip-4.3mm\sum}
 }
\def\znamenko#1{{(-1)^{#1}\cdot}}
\def\A{{\cal A}}
\def\Cob{{\Omega}}

\def\Ass{{\it Ass}}
\def\Cycl{{\it Cycl}}
\def\UAss{{\it UAss}}
\def\UPoiss{{\it UPoiss}}
\def\Poiss{{\it Poiss}}
\def\Comm{{\it Comm}}
\def\UComm{{\it UComm}}
\def\UP{{\it U}{\cal P}}

\def\bk{{\bf k}}
\def\prez#1#2{\langle #1;#2\rangle}
\def\bs{{\bf s}}
\def\Span{{\rm Span}}
\def\prezmod#1#2#3{{\langle #1;#2;#3\rangle}}
\def\M{{\cal M}}
\def\susp{\uparrow\!}
\def\ss{{\bf s\hskip0mm}}
\def\Q{{\cal Q}}
\def\dcob{{\partial_{\Cob}}}
\def\sgn{{\rm sgn}}

\def\sfF{{\sf F}}
\def\varomega{{\varphi}}

\def\osfF{{\buildrel \circ \over {{\sf F}}}}
\def\ot{{\otimes}}
\def\SS{{\bf E}}
\def\skel#1#2{{#1}^{(\leq #2)}}

\def\barD{{\overline {\cal D}}}
\def\barDelta{{\overline \Delta}}
\def\dual#1{{{#1}^*}}
\def\dualI#1{{(#1)^*}}
\def\set#1{{\{#1\}}}
\def\ZZ{{\bf Z}}
\def\tttr#1#2{{\tt T}_{{#1},\ldots,{#2}}}

\def\odstintro{{\vskip2mm\noindent%
{\bf I.\the\count88.} \global\advance\count88 by 1}}

\def\Ker{{\rm Ker}}
\def\oDelta{\stackrel{\mbox{\scriptsize o}}%
 {\Delta}\hskip-1mm\rule{0mm}{2mm}}
\def\desusp{\downarrow\!}

\long\def\comment#1\endcomment{{}}

\begin{document}
\pagestyle{myheadings}
\bibliographystyle{plain}
\baselineskip20pt plus 2pt minus 1pt
\parskip3pt plus 1pt minus .5pt

\begin{center}
{\Large \bf
Simplex, associahedron, and cyclohedron}
\end{center}

\begin{center}
{\large Martin Markl}
\end{center}

\footnote{Supported by a Fulbright grant, by the grant
AV \v CR \#1019507 and by the grant GA \v CR \#201/96/0310}

\section*{Introduction.}

\vskip2mm
\odstintro
The paper deals with three types of convex polyhedra. The most
classical is the
$n$-dimensional {\em simplex\/}
$\Delta^n$~\cite[\S10.1]{switzer:75}, the basic ingredient of
simplicial topology and perhaps one of the most important
mathematical objects at all~\cite{may:67}.

Another polyhedron is the Stasheff polytope $K_n$, also
called the
{\em associahedron\/}, the basic tool for the study of
homotopy associative Hopf
spaces~\cite[page~277]{stasheff:TAMS63}.

The last type is the polyhedron $W_n$, defined as the
Axelrod-Singer
compactification of the configuration space of $n$ points on the
circle,
and introduced by R.~Bott and
C.~Taubes~\cite[page~5249]{bott-taubes:JMP94}
in connection with the study of nonperturbative link invariants,
recently dubbed by J.~Stasheff the
{\em cyclohedron\/}~\cite{stasheff:from-ops}.

\odstintro
The crucial property of
the collection $K =\{K_n\}_{n\geq 1}$ is that it forms a
cellular {\em operad\/}~\cite[page~278]{stasheff:TAMS63}.
J.~Stasheff observed in~\cite{stasheff:from-ops}
that the collection $W = \{W_n\}_{n\geq 1}$ is a
right {\em module\/}, in the sense
of~\cite[page~1476]{markl:zebrulka}, over the
operad $K$. In Theorem~\ref{vicko} we prove
that also the collection $\Delta = \{\Delta^n\}_{n\geq 0}$ is a
natural right module over the operad $\Ass$ for associative
algebras.

\odstintro
Operads were introduced
to encode varieties of algebras. We show that, in the same
spirit,
also modules over operads describe varieties of some objects. We
call
these objects
{\em traces\/} (Definition~\ref{el}),
since they naturally generalize traces on
associative algebras.

For a so-called {\em cyclic\/} operad
$\P$~\cite[Definition~2.1]{getzler-kapranov:cyclic} we construct
a natural $\P$-module $M_{\P}$, the module {\em associated\/} to
the operad $\P$ (Definition~\ref{Turmo}).
We show that $M_\P$-traces are exactly
{\em invariant bilinear forms\/} in the sense of E.~Getzler and
M.M.~Kapranov~\cite[Definition~4.1]{getzler-kapranov:cyclic}.

\odstintro
Algebras over the cellular chain operad $CC_*(K)$ of the
associahedron are \Ainfty-algebras introduced by J.~Stasheff
in~\cite[page~294]{stasheff:TAMS63}.
They can be understood as algebras with
the usual associativity condition
\[
(ab)c = a(bc)
\]
satisfied only up to a system of coherent homotopies.
In Proposition~\ref{nuzky}
we show that
the traces over the cellular chain complex $CC_*(W)$ of the
cyclohedron are {\em homotopy traces\/} on \Ainfty-algebras,
for which the usual
condition
\[
T(ab) = T(ba)
\]
is satisfied only up to a system of coherent homotopies.
The traces over the cellular chain complex $CC_*(\Delta)$ of the
simplex are described in Theorem~\ref{vicko}.

\odstintro
The cellular chain complex $CC_*(K)$ of the associahedron has
a very
effective description -- it is the {\em operadic bar
construction\/}
on the operad $\Ass$ for associative
algebras~\cite[Example~4.1]{markl:zebrulka}.
We introduce the bar construction on a {\em module\/} over an
operad (this definition was independently made by V.~Ginzburg
and
A.A.~Voronov in~\cite{ginzburg-voronov}) and show that the
cellular
chain complex $CC_*(W)$ of the cyclohedron is the bar
construction
on the $\Ass$-module $\Cycl$, which describes traces (ordinary,
not
homotopy) on associative algebras (Theorem~\ref{ucpavka}).
A fully algebraic
description of $CC_*(\Delta)$ is given in Theorem~\ref{resiz}.

\odstintro
V.Ginzburg and
M.M.~Kapranov~\cite[Definition~4.1.3]{ginzburg-kapranov:DMJ94}
introduced so-called {\em Koszul operads\/}, with all
expected nice properties, and the related notion of the
{\em Koszul
dual\/} of a
{\em quadratic
operad\/}~\cite[\S2.1.9]{ginzburg-kapranov:DMJ94}.

We introduce analogous notions for {\em modules\/} over operads,
i.e.~we introduce {\em quadratic modules\/}, their
{\em Koszul (quadratic) duals\/} and the
property of {\em Koszulness\/} for these modules. These
definitions
were
again independently made by V.~Ginzburg and
A.A.~Voronov in~\cite{ginzburg-voronov}.

\odstintro
The operad $\Ass$ for associative algebras is
Koszul~\cite[Corollary~4.2.7]{ginzburg-kapranov:DMJ94}.
Since the operadic bar construction on $\Ass$
is the cellular chain complex of the associahedron, the
Koszulness of
$\Ass$ follows from the acyclicity of $K$, which in turn
follows from
the fact that it is a convex polyhedron.

In a similar manner, we show in Theorem~\ref{Amphora1}
that the module $\Cycl$ describing
traces on associative algebras is Koszul, as a consequence
of the
convexity of the cyclohedron $W$. A more general
argument is to observe that $\Cycl$ is the module associated
to the
cyclic operad $\Ass$ (Example~\ref{kacirek})
and then apply Theorem~\ref{Katalogizacni} saying that a
module associated to a Koszul operad is Koszul.

\odstintro
We show in Lemma~\ref{whoop} that,
for each module over an operad, there exists a
{\em spectral sequence\/}, converging to the homology of the bar
construction. We also prove in Proposition~\ref{myska}
that for modules over Koszul
operads this spectral sequence collapses.

Our spectral sequence, applied to an $\Ass$-module $\Cycl$,
carries
a strong geometrical message -- the initial term is the cellular
chain complex of the cyclohedron, while the next term is
the cellular
chain complex of the simplex. If we interpret the cyclohedron
as the compactification of the simplex constructed by a
sequence of
blow-ups~\cite[page~5249]{bott-taubes:JMP94},
then the spectral sequence
describes the inverse process -- `deblowing-up' the cyclohedron
back
to the simplex, see Section~\ref{22}.

\odstintro
{\em Some further suggestions.\/}
Consider the following `standard situation' closely related to
the topological quantum field theory. Let $C(S^m)$ be
the Axelrod-Singer
compactification~\cite[Section~5]{axelrod-singer:preprint}
of the configuration
space of distinct points in the sphere $S^m$, and ${\sf F}_m$
the compactification
of the moduli space of configurations of distinct
points in ${\bf R}^m$~\cite[\S3.2]{getzler-jones:preprint}.

It is known that ${\sf F}_n$ is a topological
operad~\cite[\S3.2]{getzler-jones:preprint} and that this
operad acts on
the
right module $C(S^m)$~\cite[Theorem~5.2]{markl:cf}
(to be precise, if $m \not= 1,3,7$,
the sphere $S^m$ is not parallelizable and we
need a suitable framed versions of the objects above).
The homology operad $H_*({\sf F}_m)$ describes a form of graded
Poisson algebras~\cite[Theorem~3.1]{getzler-jones:preprint}
(or Batalin-Vilkovisky algebras, in the framed
case)~\cite[Section~4]{G2}, and it is not difficult to see that
$H_*(C(S^n))$
is the module associated to the cyclic operad $H_*({\sf F}_m)$
in the
sense of our Definition~\ref{Turmo}.
Our paper deals with the above situation for $n=1$, while
all the
machinery cries for an application to a general situation.

Another suggestion for further research is the following.
E.~Getzler and M.M.~Kapranov introduced
in~\cite[Definition~5.2]{getzler-kapranov:cyclic} the cyclic
homology
of an algebra over a cyclic operad as the (left, nonabelian)
derived
functor of the universal invariant bilinear form functor
$\lambda(\P,-)$. We propose to study, for a (noncyclic)
operad $\P$
and a $\P$-module $M$, the derived functor of the universal
$M$-{\em trace\/} as a natural generalization of the
cyclic homology. The cyclic homology
will be then a special case for $\P$ cyclic and $M = M_\P$.

There are two ways to read the paper -- either as an exposition
of the properties of the associahedron, cyclohedron and simplex,
with some generalizations, or as a paper on general theory
of modules
over operads, with a special attention paid to the three
examples above.

\noindent
{\em Acknowledgment:\/}
I would like to express my thanks to Jim Stasheff for numerous
discussions
and hospitality during my stay at the University of North
Carolina.
Also the communication with Sasha Voronov, who was working
independently on~\cite{ginzburg-voronov}, was very useful. I
am also
very grateful to Steve Shnider and
the referee for careful reading the manuscript and
many useful remarks.

\section*{Plan of the paper:}

\noindent
\hangindent=5mm
\hangafter=1
{\em Section~\ref{1968}:
Associahedron and the cyclohedron as a truncation of the
simplex.\/}
We recall the convex realization of the associahedron as a
truncation of the simplex, due to S.~Shnider and S.~Sternberg,
and
construct a similar realization of the cyclohedron.

\noindent
\hangindent=5mm
\hangafter=1
{\em Section~\ref{bolehlav}: Cyclohedron as a module over the
associahedron.\/}
We recall (right) modules over operads and introduce traces as
algebraic objects described by these modules. We introduce
the module
$\Cycl$ for traces on associative algebras. We prove that the
cyclohedron is a module over the associahedron
and describe the corresponding traces.

\noindent
\hangindent=5mm
\hangafter=1
{\em Section~\ref{hrnicek1}: Simplex as a module over the
operad for
associative algebras.\/}
We show that the simplex is a module over the
operad $\Ass$ for associative algebras. We prove that the
associated
cellular chain complex is free and describe the corresponding
traces.

\noindent
\hangindent=5mm
\hangafter=1
{\em Section~\ref{cervena-tuzka}:
Quadratic operads and modules;
modules associated to cyclic operads.\/}
We present a class of
operads and modules having a particularly easy description. We
recall
cyclic operads and introduce the module associated to a cyclic
operad.

\noindent
\hangindent=5mm
\hangafter=1
{\em Section~\ref{penezenka}: Cyclohedron as the cobar
construction.\/} We introduce the cobar construction on a
module over
an operad. We define quadratic Koszul modules. We show that the
cellular chain complex of the cyclohedron is the cobar
construction
on the module $\Cycl$ and deduce from this fact that $\Cycl$ is
Koszul.

\noindent
\hangindent=5mm
\hangafter=1
{\em Section~\ref{22}: Cyclohedron as a compactification of the
simplex.\/}
We view the cyclohedron as a compactification of the simplex,
constructed as a sequence of blow-ups. We show that the spectral
sequence related to the cobar construction `deflates' the
cyclohedron
back to the simplex.

\noindent
\hangindent=5mm
\hangafter=1
{\em Appendix: Traces versus invariant bilinear forms.\/}
We show that traces over the module associated to a cyclic
operad are
exactly invariant bilinear forms of E.~Getzler and
M.M.~Kapranov.

\section{Associahedron and the cyclohedron as a truncation
of the
simplex}
\label{1968}

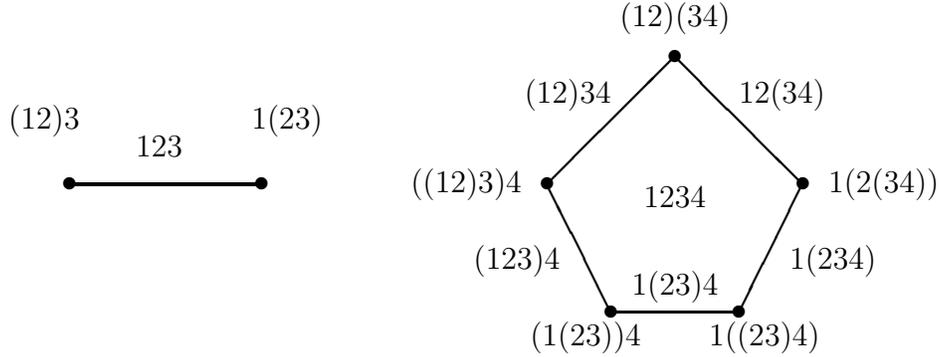
\begin{figure}[hb]
\begin{center}
%TexCad Options
%\grade{\off}
%\emlines{\off}
%\beziermacro{\off}
%\reduce{\on}
%\snapping{\off}
%\quality{2.00}
%\graddiff{0.01}
%\snapasp{1}
%\zoom{2.65}
\unitlength 1.70mm
\thicklines
\begin{picture}(54.33,25.50)
\put(32.33,12.50){\line(1,1){10.00}}
\put(42.33,22.50){\line(1,-1){10.00}}
\put(52.33,12.50){\line(-1,-2){5.00}}
\put(47.33,2.50){\line(-1,0){10.00}}
\put(37.33,2.50){\line(-1,2){5.00}}
\put(42.33,22.50){\makebox(0,0)[cc]{$\bullet$}}
\put(52.33,12.50){\makebox(0,0)[cc]{$\bullet$}}
\put(47.33,2.50){\makebox(0,0)[cc]{$\bullet$}}
\put(37.33,2.50){\makebox(0,0)[cc]{$\bullet$}}
\put(32.33,12.50){\makebox(0,0)[cc]{$\bullet$}}
\put(30.33,12.50){\makebox(0,0)[rc]{$((12)3)4$}}
\put(54.33,12.50){\makebox(0,0)[lc]{$1(2(34))$}}
\put(42.33,25.50){\makebox(0,0)[cc]{$(12)(34)$}}
\put(35.33,0.50){\makebox(0,0)[cc]{$(1(23))4$}}
\put(49.33,0.50){\makebox(0,0)[cc]{$1((23)4)$}}
\put(37.33,19.50){\makebox(0,0)[rc]{$(12)34$}}
\put(47.33,19.50){\makebox(0,0)[lc]{$12(34)$}}
\put(33.33,6.50){\makebox(0,0)[rc]{$(123)4$}}
\put(51.33,6.50){\makebox(0,0)[lc]{$1(234)$}}
\put(42.33,4.50){\makebox(0,0)[cc]{$1(23)4$}}
\put(42.33,11.50){\makebox(0,0)[cc]{$1234$}}
\put(-5.00,12.50){\line(1,0){15.00}}
\put(-5.00,12.50){\makebox(0,0)[cc]{$\bullet$}}
\put(10.00,12.50){\makebox(0,0)[cc]{$\bullet$}}
\put(2.00,15.50){\makebox(0,0)[cc]{$123$}}
\put(-7.00,17.50){\makebox(0,0)[cc]{$(12)3$}}
\put(12.00,17.50){\makebox(0,0)[cc]{$1(23)$}}
\end{picture}
\end{center}
\caption{$K_3$ (left) and $K_4$ (right).\label{k3andk4}}
\end{figure}

Let $\B(n)$ denote the set of all meaningful bracketings of $n$
independent variables $\rada1n$.
The {\em associahedron\/} $K_n$ is a convex $(n-2)$-dimensional
polyhedron whose faces are indexed by elements of $\B(n)$. To
be more precise, $\B(n)$ is a poset
(= partially-ordered set) ordered by saying that $b' \prec
b''$ if
$b''$ is
obtained from $b'$ by removing one or more pair of
brackets. Then
$K_n$ is a convex polyhedron whose poset of faces is (isomorphic
to) $\B(n)$. See Figure~\ref{k3andk4} for $K_3$ and $K_4$. A
nice
picture of $K_5$ can be found
in~\cite[page~151]{markl-stasheff:JofAlg94}.

We recall a very cute `linear convex realization' of $K_n$
as a truncation of the $(n-2)$-dimensional simplex, due to
S.~Shnider
and S.~Sternberg~\cite{shnider-sternberg:book}.
Our exposition follows the corrected version given
in~\cite[Appendix~B]{stasheff:from-ops}.

We need an alternative description of the poset $\B(n)$. Let
$P(n)$ denote the set of all proper subintervals of the interval
$[1,n-1] = \{\rada1{n-1}\}$.
Two intervals $I,J \in P(n)$ are
called {\em compatible\/}, if $I\cup J$ is not an interval
properly containing both $I$ and $J$, i.e.~if either $J\subset
I$, or
$I\subset J$, or $I\cup J$ is not an interval. Let $\I(n)$
be the set of all subsets $\iota$ of $P(n)$ such that $I$
and $J$ are
compatible for any $I,J \in \iota$.
The poset structure on $\I(n)$ is given by the set
inclusion: $\iota \preceq \kappa$ if $\kappa\subset \iota$.

\begin{lemma}
\label{prim}
(Shnider-Sternberg)
The posets $\B(n)$ and $\I(n)$ are isomorphic.
\end{lemma}

\noindent
{\bf Proof.}
For $I = [i,j]\in P(n)$, let $b(I)$ be the bracketing $1\cdots(i
\cdots j+1)\cdots n$. This correspondence is easily seen
to induce
a poset isomorphism $\I(n) \cong \B(n)$.\qed

\noindent
Define the function $c: P(n)\to \bfR_{>0}$ by $c(I):= 3^{\#I}$,
for $I
\in P(n)$. Let $K_n \subset \bfR^{n-1}$
be the convex polytope defined by
\[
K_n = \left\{(t_1,\ldots,t_{n-1})\in \bfR^{n-1};\
\sum_{k=1}^{n-1}t_k = c([1,n-1]),\
\sum_{k\in I}t_k \geq c(I),\
I\in P(n)\right\}.
\]
Denote also, for $I\in P(n)$, by $P_I$ the hyperplane
\[
P_I := \left\{(t_1,\ldots,t_{n-1})\in \bfR^{n-1};\
\sum_{k\in I}t_k =c(I)\right\}.
\]
The proof of the following proposition is given
in~\cite[Appendix~B]{stasheff:from-ops}.
\begin{proposition}(Shnider-Sternberg)
\label{Skoda}
The polytope $K_n$ has nonempty interior in the
$(n-2)$-dimensional
hyperplane $\{(t_1,\ldots,t_{n-1})\in \bfR^{n-1};\
\sum_{k=1}^{n-1}t_k = 3^{n-1}\}$. The intersection
\[
K_n \cap \bigcap\{P_{I},\ I \in \iota\}
\]
defines a nonempty $(n-\#I-2)$-dimensional face of $K_n$ for any
$\iota \in \I(n)$. All faces of $K_n$ are
obtained in this way.
\end{proposition}

\begin{figure}[hbtp]
\begin{center}
%TexCad Options
%\grade{\off}
%\emlines{\off}
%\beziermacro{\off}
%\reduce{\on}
%\snapping{\off}
%\quality{2.00}
%\graddiff{0.01}
%\snapasp{1}
%\zoom{1.00}
\unitlength 0.95mm
\linethickness{0.4pt}
\begin{picture}(85.00,100.33)
\thinlines
\put(0.00,4.33){\vector(1,0){100}}
\put(0.00,4.33){\vector(0,1){100.00}}
\put(-5.00,89.33){\line(1,-1){87.89}}
\put(-3.00,64.33){\line(1,0){65.00}}
\put(60.00,1.33){\line(0,1){60.11}}
\put(10.00,99.33){\line(0,-1){97.00}}
\put(-3.00,14.33){\line(1,0){88.00}}
\put(-2.00,84.33){\line(1,0){4.00}}
\put(-5.00,84.33){\makebox(0,0)[rc]{$24$}}
\put(-5.00,64.33){\makebox(0,0)[rc]{$18$}}
\put(-5.00,14.33){\makebox(0,0)[rc]{$3$}}
\put(60.00,-0.67){\makebox(0,0)[ct]{$18$}}
\put(80.00,-0.67){\makebox(0,0)[ct]{$24$}}
\put(10.00,-0.67){\makebox(0,0)[ct]{$3$}}
\put(2.00,98.33){\makebox(0,0)[lb]{$t_3$}}
\put(43.00,47.33){\makebox(0,0)[lc]{$P_{\{[2]\}}$}}
\put(38.00,67.33){\makebox(0,0)[cb]{$P_{\{[1,2]\}}$}}
\put(64.00,46.33){\makebox(0,0)[lc]{$P_{\{[2,3]\}}$}}
\put(13.00,95.33){\makebox(0,0)[lc]{$P_{\{[1]\}}$}}
\put(25.00,17.33){\makebox(0,0)[cb]{$P_{\{[3]\}}$}}
\put(80.00,1.33){\line(0,1){6.11}}
\thicklines
\put(10.00,14.33){\line(1,0){50.00}}
\put(60.00,14.33){\line(0,1){10.00}}
\put(60.00,24.33){\line(-1,1){40.00}}
\put(20.00,64.33){\line(-1,0){10.00}}
\put(10.00,64.33){\line(0,-1){50.00}}
\put(10.00,64.33){\makebox(0,0)[cc]{$\bullet$}}
\put(20.00,64.33){\makebox(0,0)[cc]{$\bullet$}}
\put(10.00,14.33){\makebox(0,0)[cc]{$\bullet$}}
\put(60.00,14.33){\makebox(0,0)[cc]{$\bullet$}}
\put(60.00,24.33){\makebox(0,0)[cc]{$\bullet$}}
\put(12.17,61.50){\makebox(0,0)[cc]{$\cdot$}}
\put(16.67,60.50){\makebox(0,0)[cc]{$\cdot$}}
\put(18.83,56.00){\makebox(0,0)[cc]{$\cdot$}}
\put(12.17,44.83){\makebox(0,0)[cc]{$\cdot$}}
\put(12.17,36.50){\makebox(0,0)[cc]{$\cdot$}}
\put(16.67,52.16){\makebox(0,0)[cc]{$\cdot$}}
\put(16.67,43.83){\makebox(0,0)[cc]{$\cdot$}}
\put(16.67,35.50){\makebox(0,0)[cc]{$\cdot$}}
\put(13.67,49.00){\makebox(0,0)[cc]{$\cdot$}}
\put(18.83,47.66){\makebox(0,0)[cc]{$\cdot$}}
\put(18.83,39.33){\makebox(0,0)[cc]{$\cdot$}}
\put(12.33,24.66){\makebox(0,0)[cc]{$\cdot$}}
\put(16.83,23.66){\makebox(0,0)[cc]{$\cdot$}}
\put(13.83,20.50){\makebox(0,0)[cc]{$\cdot$}}
\put(19.00,27.50){\makebox(0,0)[cc]{$\cdot$}}
\put(19.00,19.16){\makebox(0,0)[cc]{$\cdot$}}
\put(25.83,43.83){\makebox(0,0)[cc]{$\cdot$}}
\put(25.83,35.50){\makebox(0,0)[cc]{$\cdot$}}
\put(22.83,49.00){\makebox(0,0)[cc]{$\cdot$}}
\put(22.83,40.66){\makebox(0,0)[cc]{$\cdot$}}
\put(22.83,32.33){\makebox(0,0)[cc]{$\cdot$}}
\put(28.00,47.66){\makebox(0,0)[cc]{$\cdot$}}
\put(28.00,39.33){\makebox(0,0)[cc]{$\cdot$}}
\put(26.00,23.66){\makebox(0,0)[cc]{$\cdot$}}
\put(28.17,27.50){\makebox(0,0)[cc]{$\cdot$}}
\put(39.83,27.50){\makebox(0,0)[cc]{$\cdot$}}
\put(32.83,27.33){\makebox(0,0)[cc]{$\cdot$}}
\put(39.83,19.16){\makebox(0,0)[cc]{$\cdot$}}
\put(36.83,32.66){\makebox(0,0)[cc]{$\cdot$}}
\put(43.83,24.50){\makebox(0,0)[cc]{$\cdot$}}
\put(29.83,32.50){\makebox(0,0)[cc]{$\cdot$}}
\put(36.83,24.33){\makebox(0,0)[cc]{$\cdot$}}
\put(36.83,16.00){\makebox(0,0)[cc]{$\cdot$}}
\put(35.00,39.50){\makebox(0,0)[cc]{$\cdot$}}
\put(35.17,19.33){\makebox(0,0)[cc]{$\cdot$}}
\put(52.33,19.33){\makebox(0,0)[cc]{$\cdot$}}
\put(49.33,24.50){\makebox(0,0)[cc]{$\cdot$}}
\put(49.33,16.16){\makebox(0,0)[cc]{$\cdot$}}
\put(54.50,23.16){\makebox(0,0)[cc]{$\cdot$}}
\put(22.00,57.33){\makebox(0,0)[cc]{$\cdot$}}
\put(30.83,49.00){\makebox(0,0)[cc]{$\cdot$}}
\put(39.67,40.67){\makebox(0,0)[cc]{$\cdot$}}
\end{picture}
\end{center}
\caption{$(t_1,t_3)$-projection of the
convex realization of $K_4$.\label{realK4}}
\end{figure}
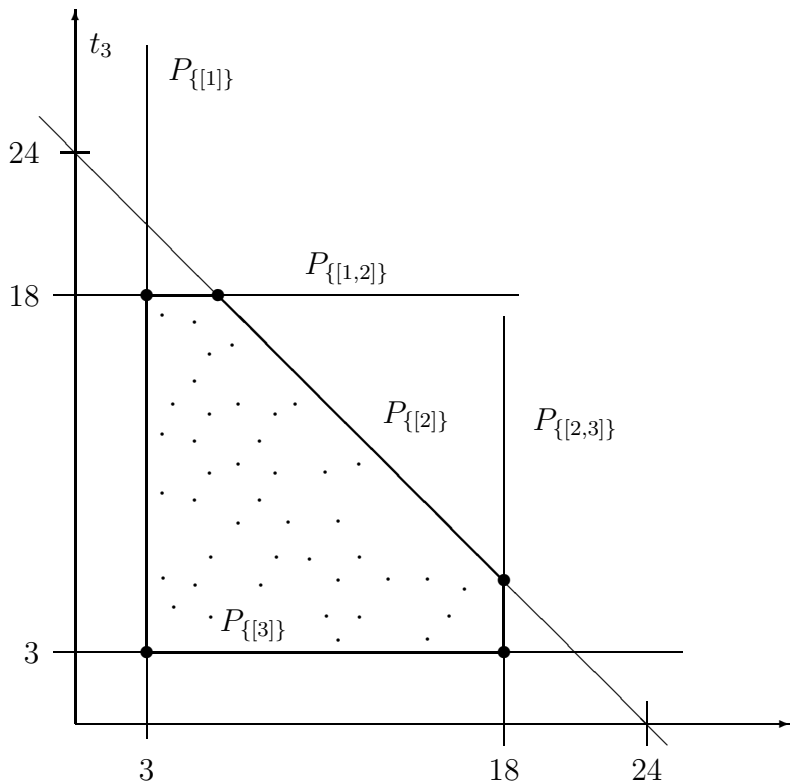

If we denote, for $\iota \in \I(n)$, by
$P_{\iota}$ the intersection $\bigcap\{P_{I},\ I \in \iota\}$,
then
the above proposition immediately implies that
the correspondence $\iota \mapsto K_n \cap P_\iota$ defines an
isomorphism of the poset $\I(n)$ and the poset of faces
of $K_n$.
This is the promised
convex realization of $K_n$. The case $n=4$ is illustrated on
Figure~\ref{realK4}.

As observed in~\cite[Appendix~B]{stasheff:from-ops},
the above construction works also
for other choices of the function $c : P(n)\to \bfR_{>0}$
provided it
is
admissible in the sense that
\[
c(I)+c(J)< c(I\cup J), \mbox{ if $I\cup J$ properly contains
both $I$
and $J$.}
\]

Let us proceed to the definition of the cyclohedron $W_n$.
As the associahedron, it will be a convex polyhedron
characterized by the poset $\BC(n)$ indexing its faces. Consider
again
$n$
independent formal variables, labeled by natural
numbers $1,\ldots,n$. The elements of $\BC(n)$ will be
equivalence
classes represented by
bracketing of the chain $\rada{\sigma(1)}{\sigma(n)}$, where
$\sigma \in \Sigma_n$ is a {\em cyclic\/} permutation. In
contrast
to the case of $\B(n)$, we allow also the bracketing which
embraces all elements. Thus, for instance, $(3(12))$
represents an
element of $\BC(n)$.

The equivalence relation is given as follows. Let $\sigma\in
\Sigma_n$
be a cyclic permutation, let $b'$ be a bracketing of
$\sigma(1)\cdots\sigma(s)$ and $b''$ a bracketing of
$\sigma(s+1)\cdots\sigma(n)$, for some $1\leq s\leq n$. Then
we identify $b'b''$ to $b''b'$. Thus, for example,
$3(12)= (12)3$ in $\BC(3)$
(but $(3(12))\not= ((12)3)$). The partial
order on $\BC(n)$ is defined, as for $\B(n)$, by deleting
pairs of
brackets.

Each element $b$ of $\BC(n)$ can be uniquely represented by
a symbol,
obtained from a representative of $b$ by forcing the
indeterminates
into the natural order. We call such symbols {\em cyclic
bracketings\/}. The formal definition will be obvious from the
following example.

\begin{example}{\rm\
The poset $\BC(2)$ contains three elements, $(12)$, $(21)$ and
$12$, where $(21)$ is represented by the cyclic bracketing
$1)(2$.
The poset structure is depicted by the interval, see
Figure~\ref{W2}.
\begin{figure}[hbtp]
\begin{center}
\unitlength=1mm
\begin{picture}(53.00,15.00)(10,10)
\thicklines
\put(20.00,10.00){\line(1,0){30.00}}
\put(20.00,10.00){\makebox(0,0)[cc]{$\bullet$}}
\put(50.00,10.00){\makebox(0,0)[cc]{$\bullet$}}
\put(35.00,15.00){\makebox(0,0)[cc]{$12$}}
\put(17.00,15.00){\makebox(0,0)[cc]{$(12)$}}
\put(53.00,15.00){\makebox(0,0)[cc]{$1)(2$}}
\end{picture}
\end{center}
\caption{$W_2$.\label{W2}}
\end{figure}

\noindent
Below are listed elements of the poset $\BC(3)$:
\[
\begin{array}[b]{|c|c||c|c|}
\hline
\mbox{elements of $\BC(3)$}&\mbox{cyclic
bracketings}&\mbox{cont.}&\mbox{\hskip11mm cont.\hskip11mm}
\\
\hline \hline
(1(23))&(1(23))&(231)&1)(23
\\
((12)3)&((12)3)&2(31)=(31)2&1)2(3
\\
((23)1)&1)((23)&(312)&12)(3
\\
(2(31))&1))(2(3&(12)3=3(12)&(12)3
\\
((31)2)&1)2)((3&(123)&(123)
\\
(3(12))&(12))(3&123=231=312&123
\\
1(23)=(23)1&1(23)&&\\
\hline
\end{array}
\]
The poset structure of $\BC(3)$ is depicted on Figure~\ref{W3}.
}\end{example}
\begin{figure}[hbtp]
\begin{center}
\unitlength2mm
\begin{picture}(62.00,35.00)(10,3)
\thicklines
\put(20.00,20.00){\line(1,1){10.00}}
\put(30.00,30.00){\line(1,0){20.00}}
\put(50.00,30.00){\line(1,-1){10.00}}
\put(60.00,20.00){\line(-1,-1){10.00}}
\put(50.00,10.00){\line(-1,0){20.00}}
\put(30.00,10.00){\line(-1,1){10.00}}
\put(20.00,20.00){\makebox(0,0)[cc]{$\bullet$}}
\put(30.00,30.00){\makebox(0,0)[cc]{$\bullet$}}
\put(50.00,30.00){\makebox(0,0)[cc]{$\bullet$}}
\put(60.00,20.00){\makebox(0,0)[cc]{$\bullet$}}
\put(50.00,10.00){\makebox(0,0)[cc]{$\bullet$}}
\put(30.00,10.00){\makebox(0,0)[cc]{$\bullet$}}
\put(28.00,35.00){\makebox(0,0)[cc]{$((12)3)$}}
\put(52.00,35.00){\makebox(0,0)[cc]{$(1(23))$}}
\put(62.00,20.00){\makebox(0,0)[lc]{$1)((23)$}}
\put(52.00,5.00){\makebox(0,0)[cc]{$1))(2(3$}}
\put(28.00,5.00){\makebox(0,0)[cc]{$1)2)((3$}}
\put(18.00,20.00){\makebox(0,0)[rc]{$(12))(3$}}
\put(40.00,35.00){\makebox(0,0)[cc]{$(123)$}}
\put(40.00,5.00){\makebox(0,0)[cc]{$1)2(3$}}
\put(27.00,25.00){\makebox(0,0)[lc]{$(12)3$}}
\put(27.00,15.00){\makebox(0,0)[lc]{$12)(3$}}
\put(53.00,25.00){\makebox(0,0)[rc]{$1(23)$}}
\put(53.00,15.00){\makebox(0,0)[rc]{$1)(23$}}
\put(40.00,20.00){\makebox(0,0)[cc]{$123$}}
\end{picture}
\end{center}
\caption{$W_3$.\label{W3}}
\end{figure}
The structure of $\BC(4)$ is indicated on Figure~\ref{W4}. The
picture
is already rather complicated, so we labeled only the vertices
(= the
minimal elements of $\BC(4)$). The label of an arbitrary face
can be easily found -- it is the least upper bound
of all vertices of the face. For example, the pentagon on the
top of $W_4$ is labeled by $123)(4$, the front hexagon is
labeled by
$(12)34$, etc.
\begin{figure}[hbtp]
\begin{center}
\unitlength=0.67mm
\begin{picture}(200.11,131.11)

%VIDITELNE HRANY:
\thicklines
\put(20.00,10.00){\line(1,0){160.05}}
\put(0.00,20.00){\line(1,1){89.87}}
\put(20.00,10.00){\line(-2,1){19.95}}
\put(89.87,109.90){\line(-3,4){9.74}}
\put(180.00,10.00){\line(2,1){19.10}}
\put(0.00,20.00){\line(0,1){23.00}}
\put(200.00,20.00){\line(0,1){23.00}}
\put(0.00,43.00){\line(1,1){80.00}}
\put(100.11,131.11){\line(5,-2){20.00}}
\put(110.00,110.00){\line(3,4){9.89}}
\put(120.00,123.00){\line(1,-1){80.11}}
\put(80.00,123.00){\line(5,2){20.11}}
\put(110.00,110.00){\line(1,-1){90.00}}
\put(110.11,110.00){\line(0,0){0.00}}
\put(90.00,110.00){\line(1,0){20.11}}

%NEVIDITELNE HRANY:
\thinlines
\put(170.00,25.00){\line(1,2){5.00}}
\put(30.00,25.00){\line(-1,2){5.00}}
\put(30.00,25.00){\line(3,1){60.00}}
\put(180.00,10.00){\line(-2,3){10.00}}
\put(85.11,55.11){\line(-3,-1){60.00}}
\put(90.00,45.00){\line(1,0){19.89}}
\put(200.00,43.00){\line(-3,-1){15.05}}
\put(175.00,35.00){\line(3,1){6.05}}
\put(109.89,45.11){\line(1,2){4.89}}
\put(100.00,65.11){\line(3,-2){15.11}}
\put(90.00,45.11){\line(-1,2){4.89}}
\put(85.00,55.00){\line(3,2){15.00}}
\put(20.00,10.00){\line(2,3){10.00}}
\put(114.78,55.11){\line(3,-1){60.22}}
\put(0.00,43.00){\line(3,-1){14.96}}
\put(25.00,35.00){\line(-3,1){5.00}}
\put(170.00,25.00){\line(-3,1){60.11}}
\put(100.00,65.00){\line(0,1){43}}
\put(100.00,131.00){\line(0,-1){17}}

%KOULE V ROZICH
\put(175.00,35.00){\makebox(0,0)[cc]{$\bullet$}}
\put(200.00,43.00){\makebox(0,0)[cc]{$\bullet$}}
\put(200.00,20.00){\makebox(0,0)[cc]{$\bullet$}}
\put(180.00,10.00){\makebox(0,0)[cc]{$\bullet$}}
\put(170.00,25.00){\makebox(0,0)[cc]{$\bullet$}}
\put(90.00,45.00){\makebox(0,0)[cc]{$\bullet$}}
\put(85.00,55.00){\makebox(0,0)[cc]{$\bullet$}}
\put(100.00,65.00){\makebox(0,0)[cc]{$\bullet$}}
\put(115.00,55.00){\makebox(0,0)[cc]{$\bullet$}}
\put(110.00,45.00){\makebox(0,0)[cc]{$\bullet$}}
\put(80.00,123.00){\makebox(0,0)[cc]{$\bullet$}}
\put(100.00,131.00){\makebox(0,0)[cc]{$\bullet$}}
\put(120.00,123.00){\makebox(0,0)[cc]{$\bullet$}}
\put(90.00,110.00){\makebox(0,0)[cc]{$\bullet$}}
\put(110.00,110.00){\makebox(0,0)[cc]{$\bullet$}}
\put(0.00,43.00){\makebox(0,0)[cc]{$\bullet$}}
\put(0.00,20.00){\makebox(0,0)[cc]{$\bullet$}}
\put(20.00,10.00){\makebox(0,0)[cc]{$\bullet$}}
\put(25.00,35.00){\makebox(0,0)[cc]{$\bullet$}}
\put(30.00,25.00){\makebox(0,0)[cc]{$\bullet$}}

%POPISKY VRCHOLU
\put(0.00,46.00){\makebox(0,0)[rb]{\scriptsize
$1)2)\!)\!(3(\!(4$}}
\put(0.00,17.00){\makebox(0,0)[rt]{\scriptsize
$(12)\!)\!)\!(3(4$}}
\put(18.00,7.00){\makebox(0,0)[ct]{\scriptsize
$(12)\!)\!(\!(34)$}}
\put(32.00,23.00){\makebox(0,0)[lc]{\scriptsize
$1)2)\!(\!(\!(34)$}}
\put(28.00,34.00){\makebox(0,0)[lc]{\scriptsize
$1)\!)2)\!(\!(3(4$}}
\put(88.00,38.00){\makebox(0,0)[cc]{\scriptsize
$1)\!)\!(2(\!(34)$}}
\put(112.00,38.00){\makebox(0,0)[cc]{\scriptsize
$1)\!(\!(\!(2(34)\!)$}}
\put(83.00,58.00){\makebox(0,0)[rc]{\scriptsize
$1)\!)\!)\!(2(3(4$}}
\put(117.00,58.00){\makebox(0,0)[lc]{\scriptsize
$1)\!(\!(\!(23)4)$}}
\put(103.00,67.00){\makebox(0,0)[lb]{\scriptsize
$1)\!)\!(\!(23)\!(4$}}
\put(78.00,126.00){\makebox(0,0)[rb]{\scriptsize
$1)2)3)\!(\!(\!(4$}}
\put(122.00,126.00){\makebox(0,0)[lb]{\scriptsize
$(1(23)\!)\!(4$}}
\put(100.00,135.00){\makebox(0,0)[cb]{\scriptsize
$1)\!(23)\!)\!(\!(4$}}
\put(112.00,110.00){\makebox(0,0)[lc]{\scriptsize
$(\!(\!12\!)3)\!)\!(\!4$}}
\put(88.00,110.00){\makebox(0,0)[rc]{\scriptsize
$(\!12\!)\!)3)\!(\!(\!4$}}
\put(180.00,5.00){\makebox(0,0)[lc]{\scriptsize
$(\!(12)\!(34)\!)$}}
\put(200.00,15.00){\makebox(0,0)[lc]{\scriptsize
$(\!(\!(12)3)4)$}}
\put(200.00,47.00){\makebox(0,0)[lc]{\scriptsize
$(\!(1(23)\!)4)$}}
\put(166.00,22.00){\makebox(0,0)[rc]{\scriptsize
$(1(2(34)\!)\!)$}}
\put(172.00,34.00){\makebox(0,0)[rc]{\scriptsize
$1(\!(23)4)\!)$}}
\end{picture}
\end{center}
\caption{$W_4$.\label{W4}}
\end{figure}

We construct, mimicking the approach of Shnider and Sternberg, a
convex realization of the poset $\BC(n)$. First some
terminology. By a {\em cyclic subinterval\/} of $[1,n]$ we mean
either a `normal' subinterval $[i,j]$, $1\leq i\leq j \leq n$,
representing the subset $\{\rada ij\}$ of
$\{\rada 1n\}$, or the symbol $i][j$,
$1\leq i <j \leq n$, representing $\{\rada 1i\}\cup \{\rada
jn\}$. We
will always suppose that the corresponding sets are proper
subsets of
$\{\rada 1n\}$, i.e.~we exclude the intervals $[1,n]$ and
$i][i+1$,
for $1\leq i <n$. Let us denote by $PC(n)$ the set of all cyclic
subintervals in the above sense. We denote by $\IC(n)$ the
set of
all subsets of $PC(n)$ consisting of {\em nested\/}
subintervals,
meaning that, for $I,J \in \iota \in \IC(n)$, either $I\subset
J$ or
$J \subset I$.
Again, $\IC(n)$ is a poset, the order being induced by the
inclusion.
We have the following analog of Lemma~\ref{prim}.
\begin{lemma}
The posets $\BC(n)$ and $\IC(n)$ are isomorphic.
\end{lemma}

\noindent
{\bf Proof.}
Define, for $I \in PC(n)$, the cyclic bracketing $b(I)\in
\BC(n)$ by
\[
b(I)=\left\{
\begin{array}{ll}
1\cdots(i\cdots j+1)\cdots n, &\mbox{ for $I = [i,j],\ j< n$,}
\\
1)\cdots(i\cdots n, &\mbox{ for $I = [i,n]$, and}
\\
1\cdots i+1)\cdots(j\cdots n, &\mbox{ for $I = i][j$.}
\end{array}
\right.
\]
This correspondence induces the desired poset
isomorphism.\qed

Our convex realization of $\BC(n)$, whose possibility
was predicted in~\cite[Appendix~B]{stasheff:from-ops}),
is defined as follows. Let
$W_n \subset \bfR^n$ be the convex polyhedron
\[
W_n = \left\{(t_1,\ldots,t_n)\in \bfR^n;\
\sum_{k=1}^nt_k = c([1,n]),\
\sum_{k\in I}t_k \geq c(I),\
I\in PC(n)\right\}.
\]
For $I\in P(n)$, let $P_I$ be the hyperplane
\[
P_I := \left\{(t_1,\ldots,t_n)\in \bfR^n;\
\sum_{k\in I}t_k =c(I)\right\}.
\]
The proof of the following proposition is a straightforward
modification of the proof of Proposition~\ref{Skoda} as given
in~\cite[Appendix~B]{stasheff:from-ops}).
\begin{proposition}
The polytope $W_n$ has nonempty interior in the
$(n-1)$-dimensional
hyperplane $\{(t_1,\ldots,t_n)\in \bfR^n;\
\sum_{k=1}^n t_k = c([1,n])\}$. The intersection
\[
W_n \cap \bigcap\{P_{I},\ I \in \iota\}
\]
defines a nonempty $(n-\#I-1)$-dimensional face of $W_n$ for any
$\iota\in \IC(n)$ and all faces of $W_n$ are
obtained in this way.
\end{proposition}

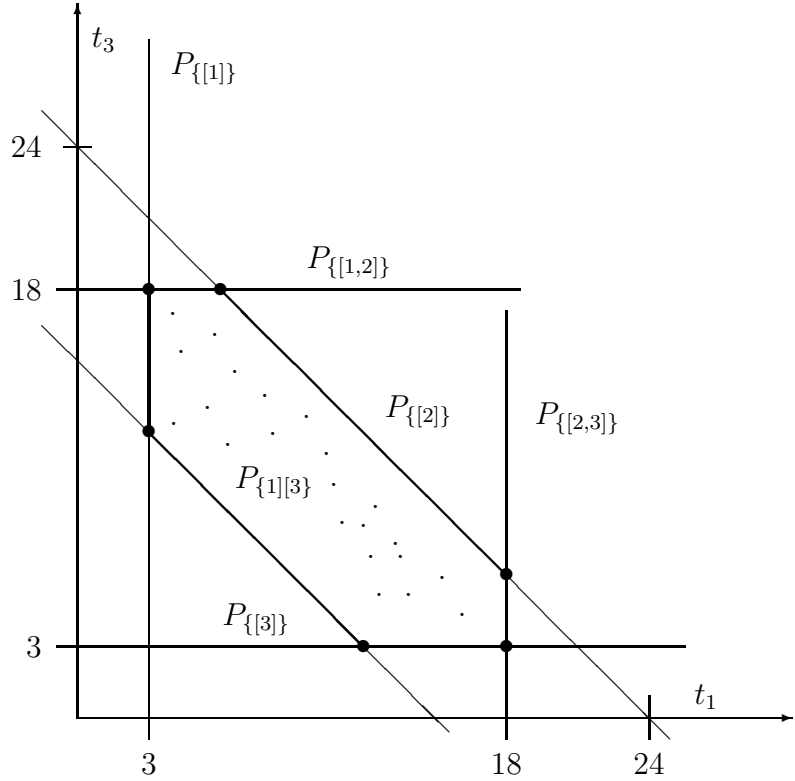
\begin{figure}[hbtp]
\begin{center}
%TexCad Options
%\grade{\off}
%\emlines{\off}
%\beziermacro{\off}
%\reduce{\on}
%\snapping{\off}
%\quality{2.00}
%\graddiff{0.01}
%\snapasp{1}
%\zoom{2.00}
\unitlength 0.95mm
\linethickness{0.4pt}
\begin{picture}(100.00,100.00)(0,-7)
\thinlines\put(-5,80.00){\makebox(0,0)[rc]{$24$}}
\put(-5.00,60.00){\makebox(0,0)[rc]{$18$}}
\put(-5.00,10.00){\makebox(0,0)[rc]{$3$}}
\put(60.00,-5.00){\makebox(0,0)[ct]{$18$}}
\put(80.00,-5.00){\makebox(0,0)[ct]{$24$}}
\put(10.00,-5.00){\makebox(0,0)[ct]{$3$}}
\put(0.00,0.00){\vector(1,0){100.00}}
\put(0.00,0.00){\vector(0,1){100.00}}
\put(-5.00,85.00){\line(1,-1){87.89}}
\put(-3.00,60.00){\line(1,0){65.00}}
\put(60.00,-3.00){\line(0,1){60.11}}
\put(10.00,95.00){\line(0,-1){98.00}}
\put(-3.00,10.00){\line(1,0){88.00}}
\put(-2.00,80.00){\line(1,0){4.00}}
\put(88.00,2.00){\makebox(0,0)[cb]{$t_1$}}
\put(2.00,94.00){\makebox(0,0)[lb]{$t_3$}}
\put(43.00,43.00){\makebox(0,0)[lc]{$P_{\{[2]\}}$}}
\put(38.00,63.00){\makebox(0,0)[cb]{$P_{\{[1,2]\}}$}}
\put(64.00,42.00){\makebox(0,0)[lc]{$P_{\{[2,3]\}}$}}
\put(13.00,91.00){\makebox(0,0)[lc]{$P_{\{[1]\}}$}}
\put(25.00,13.00){\makebox(0,0)[cb]{$P_{\{[3]\}}$}}
\put(80.00,-3.00){\line(0,1){6.11}}
\put(-5.00,55.00){\line(1,-1){57.00}}
\put(22.00,33.00){\makebox(0,0)[lc]{$P_{\{1][3\}}$}}
\thicklines
\put(60.00,10.00){\line(0,1){10.00}}
\put(60.00,20.00){\line(-1,1){40.00}}
\put(20.00,60.00){\line(-1,0){10.00}}
\put(10.00,60.00){\makebox(0,0)[cc]{$\bullet$}}
\put(20.00,60.00){\makebox(0,0)[cc]{$\bullet$}}
\put(60.00,10.00){\makebox(0,0)[cc]{$\bullet$}}
\put(60.00,20.00){\makebox(0,0)[cc]{$\bullet$}}
\put(10.00,60.00){\line(0,-1){20.00}}
\put(10.00,40.00){\line(1,-1){30.00}}
\put(40.00,10.00){\line(1,0){20.00}}
\put(10.00,40.00){\makebox(0,0)[cc]{$\bullet$}}
\put(40.00,10.00){\makebox(0,0)[cc]{$\bullet$}}
\put(13.33,56.67){\makebox(0,0)[cc]{.}}
\put(19.17,53.67){\makebox(0,0)[cc]{.}}
\put(14.50,51.33){\makebox(0,0)[cc]{.}}
\put(22.00,48.50){\makebox(0,0)[cc]{.}}
\put(26.17,45.17){\makebox(0,0)[cc]{.}}
\put(45.17,22.67){\makebox(0,0)[cc]{.}}
\put(32.00,42.17){\makebox(0,0)[cc]{.}}
\put(51.00,19.67){\makebox(0,0)[cc]{.}}
\put(27.33,39.83){\makebox(0,0)[cc]{.}}
\put(46.33,17.33){\makebox(0,0)[cc]{.}}
\put(34.83,37.00){\makebox(0,0)[cc]{.}}
\put(53.83,14.50){\makebox(0,0)[cc]{.}}
\put(35.83,32.67){\makebox(0,0)[cc]{.}}
\put(41.67,29.67){\makebox(0,0)[cc]{.}}
\put(37.00,27.33){\makebox(0,0)[cc]{.}}
\put(44.50,24.50){\makebox(0,0)[cc]{.}}
\put(40.00,27.00){\makebox(0,0)[cc]{.}}
\put(41.00,22.67){\makebox(0,0)[cc]{.}}
\put(42.17,17.33){\makebox(0,0)[cc]{.}}
\put(18.17,43.50){\makebox(0,0)[cc]{.}}
\put(13.50,41.17){\makebox(0,0)[cc]{.}}
\put(21.00,38.33){\makebox(0,0)[cc]{.}}
\end{picture}
\end{center}
\caption{$(t_1,t_3)$-projection of the convex realization of
$W_3$.\label{realW3}}
\end{figure}

For $\iota \in \IC(n)$, let
$P_{\iota}$ be the intersection $\bigcap\{P_{I},\ I \in
\iota\}$.
Then the correspondence $\iota \mapsto W_n \cap P_\iota$
defines an
isomorphism between the poset $\IC(n)$ and the poset of faces of
the polytope $W_n$. This is our
convex realization of $W_n$. The convex realization of $W_3$
is shown
on Figure~\ref{realW3}.

\begin{observation}{\rm\
The cyclohedron $W_n$ has $n(n-1)$ codimension-one faces,
represented
by the bracketings
\begin{equation}
\label{csa}
b_{k,n} :=
(\rada{\sigma(1)}{\sigma(k)})\rada{\sigma(k+1)}{\sigma(n)},\ 1 <
k\leq n,
\end{equation}
where $\sigma \in \Sigma_n$ is a cyclic permutation. The face
represented by the bracketing $b_{k,n}$ is isomorphic to the
product $W_{n-k+1}\times K_k$. For example, $W_4$ depicted on
Figure~\ref{W4}, has
\begin{itemize}
\item[-]
$4$ hexagonal faces,
corresponding to $(12)34$, $(23)41$, $(34)12$ and $(41)23$,
isomorphic to $W_3 \times K_2 = W_3 \times \mbox{point}$,
\item[-]
4 square faces, corresponding to $(123)4$, $(234)1$, $(341)2$
and
$(412)3$, isomorphic to $W_2 \times K_3$, and
\item[-]
4 pentagonal faces, corresponding to $(1234)$, $(2341)$,
$(3412)$ and
$(4123)$, isomorphic to $W_1 \times K_4 = \mbox{point}
\times K_4$.
\end{itemize}
This was observed by J.~Stasheff who realized that this is
a strong
motivation for the existence of a module structure which we will
discuss in the following Section~\ref{bolehlav}.
}\end{observation}
\begin{observation}{\rm\
It is clear from our constructions that the cyclohedron $W_n$
is a
truncation of the associahedron $K_{n+1}$, for $n\geq 1$.
This is, of course, a trivial statement -- any convex
polyhedron is
a truncation of an arbitrary other convex polyhedron of the same
dimension, so we must be more precise.

For any $n\geq 1$ there is an obvious map
$P(n+1)\hookrightarrow PC(n)$ which decomposes $PC(n)$ as
\[
PC(n)= P(n+1) \sqcup E(n),
\]
where $E(n)$ is the subset of `exotic' cyclic intervals of
the form
$i][j$, $1\leq i<j\leq n$. The polyhedron $W_n$ is then the
truncation of $K_{n+1}$ by hyperplanes $P_I$ indexed by the
`exotic'
intervals $I\in E(n)$. Compare Figures~\ref{realK4}
and~\ref{realW3}
for $n=3$. We do not
know whether this observation has any deeper meaning.
}\end{observation}

\begin{observation}{\rm
\label{sdff}
Choose, for each $t=\rada 1n$, a point $P_t$
in the interior of the codimension
one face of $W_n$ corresponding to $(\rada tn,\rada 1{t-1})$.
The convex hull of the set
$\set{\rada{P_1}{P_n}}$ is a simplex, closely related to the
`deblowing up' of $W_n$ described in Section~\ref{22}. We
will use
this simplex to introduce an orientation of $W_n$.
}\end{observation}

\section{Cyclohedron as a module over the associahedron}
\label{bolehlav}

We believe that there is no need to give a detailed definition
of an
operad. Recall only that an {\em operad\/} (in a symmetric
monoidal
category $\C= (\C,\times)$) is a
sequence $\P = \{\P(n); n\geq 1\}$ of objects of $\C$
together with morphisms
\[
\gamma=\gamma_{m_1,\ldots,m_l}:\P(l)\times
\P(m_1)\times\cdots\times\P(m_l)
\longrightarrow
\P(m_1+\cdots+m_l),
\]
given for any $l,m_1,\ldots,m_l \geq 1$,
satisfying the usual axioms~\cite[Definition~3.12]{may:1972}.
If not stated otherwise, we assume our operads to be
{\em symmetric\/}, i.e.~we assume
that each $\P(n)$ has a right action of the symmetric group
$\Sigma_n$, $n\geq 2$,
which has again to satisfy some
axioms~\cite[Definition~1.1]{may:1972}. We frequently write
$p(\rada{p_1}{p_l})$ instead of $\gamma(p,\rada{p_1}{p_l})$.

One comment concerning the action of the symmetric group is
in order
here. Our convention is determined by the conventional choice
of the
multiplication in the symmetric group. We accepted the standard
one with
$\sigma \cdot \tau$ meaning $\sigma(\tau)$, i.e.~the permutation
(= a
map) $\tau$ followed by $\sigma$. Then $\P(n)$ must
be a
{\em right\/} $\Sigma_n$-module, which is the convention used
in the
original definition of P.~May quoted above.

Recall that, for any object $V\in \C$, there
exist the so-called {\em endomorphism operad\/} $\E_V =
\coll{\E_V}$
with $\E_V(n):= \Hom {V^{\times n}}V$. If $\P$ is an operad
in $\C$,
then a {\em$\P$-algebra
structure on $V$\/} is an operad map $a: \P \to \E_V$.

\begin{example}{\rm\
\label{hrnicek}
The collection $\B = \coll {\B}$ introduced in
Section~\ref{1968}
has a structure of a (nonsymmetric) operad in the
category of posets. The composition $\gamma(b;\rada{b_1}{b_l})$
is,
for $b\in \B(l)$ and $b_i \in \B(m_i)$, defined as the
bracketing
$b(\rada{b_1}{b_l})$ obtained by inserting $b_i$ at the $i$-th
position in $b$, $1\leq i\leq l$. We believe that it is clear
what we
mean by this. For example
\begin{eqnarray*}
\gamma_{\B}(12;1(23),12)& =& 1(23)45,
\\
\gamma_{\B}((12)3;(12)3,12,(12)(24))
& =& ((12)345)(67)(89),\ \mbox{etc.}
\end{eqnarray*}
A classical result of
J.~Stasheff~\cite[page~278]{stasheff:TAMS63} says
that
the collection of the associahedra $K = \{K_n\}_{n\geq 1}$ has a
cellular
(nonsymmetric) operad structure which induces on the collection
$\B =
\coll {\B}$ of
its faces the operad structure of Example~\ref{hrnicek}.
As a consequence, cellular chains on $K$ form an
operad $\A = \coll \A$, $\A(n)= CC_*(K_n)$, in the category of
differential graded
vector spaces. Algebras over the operad $\A$ are
\Ainfty-algebras~\cite[page~294]{stasheff:TAMS63}.
}\end{example}

Probably the most effective way to describe the operad $\A$
is to say
that $\A = \Cob(\dualI{\bs \Ass})$, the operadic cobar
construction on the
dual cooperad $\dualI {\bs \Ass}$, where $\bs \Ass$ is the
suspension of the
operad for associative
algebras. This is the same, since the operad $\Ass$ is Koszul,
as to say that $\A$ is the {\em minimal model\/} of the
associative
operad $\Ass$. All this is explained
in~\cite{markl:zebrulka}. We are
going
to
make
a similar analysis for the cyclohedron.

\begin{definition}
\label{Amphora}
A (right) module over an operad $\P$ is a collection
$M = \coll M$ such that each $M(n)$ is, for $n\geq 1$, a
$\Sigma_n$-module,
together with morphisms
\begin{equation}
\label{dymka}
\nu= \nu_{\rada{m_1}{m_l}}: M(l) \times \P(m_1)\times \cdots
\times
\P(m_l)
\longrightarrow M(m_1+\cdots+ m_l),
\end{equation}
given for any $l, m_1,\ldots,m_l \geq 1$.
The structure maps $\nu_{\rada{m_1}{m_l}}$ must satisfy
the axioms
which are obtained by replacing, in the May's definition of an
operad~\cite[Definition~3.12]{may:1972}, the first occurrence
of $\P$
by $M$, see~\cite[Definition~1.3]{markl:zebrulka} for details.

\end{definition}
\begin{remark}{\rm\
Observe the resemblance of the above definition to the
definition of
an operad. This is due to the fact that right modules over
an operad
are
special cases of general modules, which are abelian groups
object in
a certain comma category of operads, see the discussion
in~\cite[page~1476]{markl:zebrulka}.
}\end{remark}

\begin{remark}{\rm\
\label{mince}
The structure map $\nu$ is, as in the
case of operads, determined by the system of `comp' maps
\begin{equation}
\label{comp}
\circ_i :M(m) \otimes \P(n) \to M(m+n-1),\ m,n \geq 1,\
1\leq i\leq m,
\end{equation}
defined by $\circ_i (x,p) := \nu(x;1,\ldots,1,p,1,\ldots,1)$
($p$ at
the $i$-th position) which have to satisfy certain
axioms~\cite[Formula~(1)]{markl:zebrulka}.
}\end{remark}

\begin{example}{\rm\
An operad $\P$ is a
module over itself. Very important nontrivial examples are
provided by the Axelrod-Singer compactification of configuration
spaces of points in a manifold. The result is a module over
the operad
of `local' configurations, see~\cite{markl:cf}.

There is the following analog of the endomorphism operad. Let
$A,W$ be objects of the category $\C$. Then the
collection $\E_{A,W}= \coll{\E_{A,W}}$ with $\E_{A,W}(n):=
\Hom{A^{\times n}}W$ is a module over the endomorphism operad
$\E_A$,
the module structure being given by the obvious composition
of maps,
as in the case of the endomorphism operad. A $\P$-algebra
structure $a: \P \to \E_A$ on $A$ induces a $\P$-module
structure on $\E_{A,W}$.
}\end{example}

We are going now to introduce
objects described, in the similar sense as algebras are
described by
operads, by {\em modules\/} over
operads. We will call them, from the reasons which will be
explained
in Example~\ref{why-traces},
{\em traces\/} over algebras.

\begin{definition}
\label{el}
Let $M$ be a $\P$-module and let $A$ be a $\P$-algebra. An
$M$-trace
over
$A$ is a map $t : M \to \E_{A,W}$ of $\P$-modules, where
$\E_{A,W}$
has the $\P$-module structure induced from the $\P$-algebra
structure
on $A$.
\end{definition}

\begin{example}{\rm\
Rather dull examples of traces are given by taking $M= \P$. For
example, an $\Ass$-trace over an associative algebra is (given
by) a
bilinear map $B : A \times A \to W$ such that $B(ab,c)
= B(a,bc)$,
$a,b,c\in A$,
i.e.~by an (not necessary symmetric) invariant bilinear form.
}\end{example}

We will need the following notation.
Let, for permutations $\sigma \in \Sigma_l$ and
$\sigma_i\in \Sigma_{m_i}$, $1\leq i\leq l$,
\begin{equation}
\label{myska1}
\sigma(\rada{\sigma_1}{\sigma_l})\in \Sigma_{m_1+\cdots+ m_l}
\end{equation}
denote the permutation $\sigma(\rada{m_1}{m_l})\cdot (\sigma_1
\oplus \cdots \oplus \sigma_l)$, where the meaning of $\sigma_1
\oplus \cdots \oplus \sigma_l$ is clear and
$\sigma(\rada{m_1}{m_l})$
permutes the blocks of $\rada{m_1}{m_l}$-elements via $\sigma$.
This defines a map $\Sigma_l
\times \Sigma_{m_1}\squeeze \Sigma_{m_l}
\to \Sigma_{m_1+\cdots +m_l}$.

\begin{example}{\rm\
\label{tuzka}
More interesting example of a trace can be constructed as
follows.
Take again the (symmetric) operad \Ass\ for associative
algebras.
Recall that $\Ass(n) = \bk[\Sigma_n]$, the group ring of the
symmetric group
over the ground field $\bk$. The operad structure map $\gamma =
\gamma_{\Ass}$ is defined by
$\gamma(\sigma;\rada{\sigma_1}{\sigma_l}) =
\sigma(\rada{\sigma_1}{\sigma_l})$, where
$\sigma(\rada{\sigma_1}{\sigma_l})$ has the same meaning as
in~(\ref{myska1}).

The group of cyclic permutations $\Z_n = \Z/n\Z$
acts from the left on
$\Sigma_n$.
The group $\Sigma_{n-1}$ is imbedded in $\Sigma_n$ as
permutations
leaving $1$ fixed. This embedding is a cross-section to
the ${\bf
Z}_n$-action, thus we can identify $\Sigma_{n-1}$ as a {\em
set\/} to
the coset space ${\bf Z}_n \backslash \Sigma_n$.
The projection $\pi_n : \Sigma_n \to
\Z_n \backslash \Sigma_n \cong \Sigma_{n-1}$ then induce on
$\bk[\Sigma_{n-1}]$ a structure of a right $\Sigma_n$-module.
Define the collection $\Cycl = \coll{\Cycl}$ by
$\Cycl(n):= \bk[\Sigma_{n-1}]$, $n\geq 1$.
The system of maps $\{\pi_n: \Sigma_n \to
\Sigma_{n-1}\}_{n\geq 1}$ is the projection $\pi :
\Ass \to \Cycl$ of collections.

\begin{lemma}
The projection $\pi : \Ass \to \Cycl$ induces on $\Cycl$
the structure
of a module over the operad $\Ass$.
\end{lemma}

\noindent
{\bf Proof.}
The structure maps $\nu = \nu_{\Cycl}$ are determined, for
$\sigma \in \Sigma_l$ and $\sigma_i \in \Sigma_{m_i}$, $1\leq i
\leq l$, by
\[
\nu(\pi(\sigma);\rada{\sigma_1}{\sigma_l}) :=
\gamma_{\Ass}(\sigma;\rada{\sigma_1}{\sigma_l}),
\]
where $\gamma_{\Ass}(\sigma;\rada{\sigma_1}{\sigma_l})=
\sigma(\rada{\sigma_1}{\sigma_l})$.
The proof is then finished by an easy verification
that, if $\sigma' \equiv \sigma''$ mod $\Z_l$, then
\[
\sigma'(\rada{\sigma_1}{\sigma_l}) \equiv
\sigma''(\rada{\sigma_1}{\sigma_l}) \mbox{ mod }
\Z_{m_1+\cdots +m_l},
\]
which we leave to the reader.\qed
}\end{example}

\begin{example}{\rm\
\label{why-traces}
A $\Cycl$-trace over an
associative algebra $A$ is (characterized by)
a map $T: A \to W$ such that
\begin{equation}
\label{nabla}
T(ab) = \znamenko{|a|\cdot |b|}T(ba),\ a,b\in A,
\end{equation}
i.e.~$T$ is a trace in the usual sense.
We postpone the verification of this statement to
Example~\ref{nabijecka}.
}\end{example}

In the rest of this section we show that the collection $W :=
\colla W$ of the cyclohedra is a natural cellular (right)
module over
the cellular operad $K = \colla K$ and describe $W$-traces on an
A($\infty$)- (= $K$)-algebra.

It is convenient to consider harmless symmetrizations. Recall
that
$\B(n)$ was the poset of all bracketings of $\rada 1n$. Take
instead
be the poset $\barB(n)$ of all bracketings of
$\rada{\sigma(1)}{\sigma(n)}$, $\sigma \in \Sigma_n$. Obviously
$\barB(n) = \Sigma_n \times \B(n)$ and $\barB(n)$ is the
poset of
faces of the $n!$-connected polyhedron $\barK_n :=
\Sigma_n\times K_n$. The collection $\barK := \colla \barK$ is a
(symmetric) cellular operad and the corresponding operad
of cellular
chains $\barA := CC_*(\barK)$ is the (symmetric) operad for
\Ainfty-algebras.

There is a similar symmetrization of the cyclohedron.
We introduced $\BC(n)$ as the poset of equivalence classes of
bracketings of
$\rada {\sigma(1)}{\sigma(n)}$ with a {\em cyclic\/} permutation
$\sigma \in \Sigma_n$. If we admit all permutations, we
obtain the
poset $\barBC(n)= \Sigma_{n-1}\times \BC(n)$ whose realization
is the $(n-1)!$-connected polyhedron $\barW_n := \Sigma_{n-1}
\times W_n$.
\begin{lemma}
\label{paska}
The collection $\barBC := \coll{\barBC}$ is a natural module
over
the operad $\barB:= \coll{\barB}$ in the symmetric monoidal
category
of posets.
\end{lemma}

\noindent
{\bf Proof.}
The easiest way to define the module structure is the following.
Let $b$ be a bracketing of $\rada{1}{l}$ representing
an element $[b]\in \BC(l)\subset \barBC(l)$. Let $b_i \in
\B(m_i)\subset \barB(m_i)$
be, for $1\leq i \leq l$ ,
a bracketing of $\rada{1}{m_i}$. Then we define
$\nu_{\barBC}(b;\rada{b_1}{b_l}) \in \barBC(m)$, $m =
m_1+\cdots m_l$, to be
the element represented by the composite (in the same sense
as in
Example~\ref{hrnicek}) bracketing $b(\rada{b_1}{b_l})$
of $m$.

The set $\barBC(l)$ is $\Sigma_l$-generated by elements of the
same form as $b$, i.e.~by elements represented by a bracketing
of
the `unpermuted' string $\rada1l$,
and the same is true also for $\barB(m_i)$, $1\leq i\leq l$.
Thus the formula for the composition of arbitrary elements is
dictated by the equivariance of the module composition map. We
leave
the verification that this definition is correct to the reader.
\qed

\begin{theorem}
\label{sirky}
The collection $\barW := \colla{\barW}$ carries a structure of
a module over the operad $\barK := \colla{\barK}$ in the
category of
cellular complexes which induces, on the level of the poset of
faces, the structure of Lemma~\ref{paska}.

The homology collection $H_*(\barW)= \{H_*(\barW_n)\}_{n\geq 1}$
coincides, as an $\Ass = H_*(\barK)$-module, to the module
$\Cycl$
introduced in Example~\ref{why-traces}.
\end{theorem}

\noindent
{\bf Proof.}
The proof is a modification of the proof of the existence of an
operad structure on the associahedron, given by J.~Stasheff
in~\cite[page~278]{stasheff:TAMS63}.
By Remark~\ref{mince}, the $\barK$-module
structure on $\barW$ is given by the `comp' maps
\[
\circ_i :\barW_m \times \barK_n \to \barW_{m+n-1},\
m,n\geq 1,\ 1\leq i\leq m.
\]
As a matter of fact, in our case it is enough to specify
\begin{equation}
\label{vrsek}
\circ_1 :W_m \times K_n \to W_{m+n-1},
\end{equation}
the remaining `comp' maps are determined by the equivariance. We
define $\circ_1$ of~(\ref{vrsek}) to be the identification
of the
product $W_m \times K_n$ to the face of $W_{m+n-1}$ indexed by
the bracketing $b_{n,n+m-1}$ of~(\ref{csa}). The second part is
immediate.\qed

\begin{observation}{\rm\
The $\barK$-module structure on the `symmetrized'
cyclohedron $\barW$ restricts to the right action of the
nonsymmetric operad $K$ on the `nonsymmetric' cyclohedron $W$.
This is the structure observed by J.~Stasheff
in~\cite[Section~4]{stasheff:from-ops}.
Another argument for the existence of the module
structure of Theorem~\ref{sirky} is the interpretation of the
cyclohedron to the compactification of a configuration space,
see
Section~\ref{22}.
}\end{observation}

Let us consider the $\barA$-module $\M:= CC_*(\barW)$ of
cellular chains on the cyclohedron.
To describe traces over $\M$,
it is convenient to accept the following
notation. For graded indeterminates $\rada{a_1}{a_n}$
and a permutation $\sigma\in
\Sigma_n$, the {\em Koszul sign\/}
$\epsilon(\sigma)=\epsilon(\sigma;\rada{a_1}{a_n})$
is defined by
\[
a_1\land\dots\land a_n = \epsilon(\sigma;a_1,\dots,a_n)
\cdot a_{\sigma(1)}\land \dots \land a_{\sigma(n)},
\]
which has to be satisfied in the free graded commutative algebra
$\ext(\rada{a_1}{a_n})$. Denote also
\[
\chi(\sigma)=\chi(\sigma;\rada{a_1}{a_n}) :=
\sgn(\sigma)\cdot \epsilon(\sigma;\rada{a_1}{a_n}).
\]
For an expression
$X(\rada {a_1}{a_n})$ in indeterminates
$\rada {a_1}{a_n}$, let the {\em cyclic sum\/}
\begin{equation}
\label{cyclsum}
\cyclsum X(\rada{a_1}{a_n}) :=
\sum_\sigma \chi(\sigma) X(\rada {a_{\sigma(1)}}{a_{\sigma(n)}})
\end{equation}
be the summation over all cyclic permutations. A convincing
example
of the use of this convention is the
(graded) Jacobi identity written as
\[
\cyclsum [a_1,[a_2,a_3]] = 0.
\]

The following proposition, whose proof we postpone after
Theorem~\ref{ucpavka},
describes $\M$-traces over \Ainfty-algebras.

\begin{proposition}
\label{nuzky}
Let $A = (A; m_1=\partial,m_2,m_3,\ldots)$ be an
\Ainfty-algebra~(\cite[page~294]{stasheff:TAMS63}, but we use
the sign
convention of~\cite[\S1.4]{markl:JPAA92}). Then
an
$\M$-trace is given by a differential graded vector space $W =
(W,\delta)$, $\deg{\delta} = -1$,
and a system $T_n :A^{\otimes n}\to W$
of degree-$(n-1)$ linear
maps, $n\geq 1$, such that, for all $\rada{a_1}{a_n} \in A$,
\begin{itemize}
\item[(i)]
$T_n(\rada{a_1}{a_n})= \chi(\sigma)\cdot
T_n(\rada{a_{\sigma(1)}}{a_{\sigma(n)}})$
for all
cyclic permutations $\sigma \in \Sigma_n$, and
\item[(ii)]
for all $n\geq 1$,
\begin{equation}
\label{Ax}
\delta T_n(\rada{a_1}{a_n}) =
\cyclsum \sum_{1\leq k\leq n}
\znamenko{k+n} T_{n-k+1}(m_k(\rada{a_1}{a_{k}}),
\rada{a_{k+1}}{a_n}).
\end{equation}
\end{itemize}
We call such objects homotopy traces over an \Ainfty-algebra
$A$.
\end{proposition}

\noindent
Let us write down the axiom~(\ref{Ax}) explicitly for small
$n$. For
$n=1$
it gives
\[
\partial T_1(a) = T_1 (\delta(a)),\ a\in A,
\]
which means that $T_1$ is a homomorphism of differential graded
spaces $(A,\partial)$ and $(W,\delta)$. For $n=2$ it becomes
\[
\delta T_2(a,b) +T_2(\partial a,b) -\znamenko{|a|\cdot |b|}
T_2(\partial b,a) =
T_1(m_2(a,b)) -\znamenko{|a|\cdot |b|}
 T_1(m_2(b,a)),\ a,b\in A,
\]
i.e.~$T_1$ is a trace for the `multiplication' $m_2$ up to
a homotopy $T_2$.
For higher $n$'s, the axiom~(\ref{Ax})
represents `coherence conditions' for
the homotopy $T_2$. An important special case is when $A$ is an
ordinary associative algebra, that is the only nontrivial
structure map is $m_2$, which is an associative multiplication
$\cdot$. The axioms for the corresponding trace are obtained by
putting, in
Proposition~\ref{nuzky},
$m_k=0$ for $k\geq 3$. We also substitute
$(-1)^\frac{n(n-1)}{2}T_n$ for
$T_n$, to get rid of the overall sign $(-1)^n$.
A homotopy trace is then a system $\{T_n:A^{\otimes n}\to W
\}_{n\geq 1}$ of degree-$(n-1)$ linear maps,
satisfying~\ref{nuzky}(i) and
\begin{eqnarray}
\delta T_1(a)\!\! &=&\!\! 0\nonumber
\\
\delta T_2 (a,b)\!\! &=& \!\!
T_1(a\cdot b)-\znamenko{|a|\cdot |b|}
T_1(b\cdot a)\nonumber
\\
\delta T_3 (a,b,c)\!\! &=&\!\!T_2(a\cdot b,c)
+\znamenko{|a|\cdot(|b|+|c|)} T_2(b\cdot c,a)
+\znamenko{|c|\cdot(|a|+|b|)} T_2(c\cdot a,b)\nonumber
\\
&\vdots&\nonumber
\\
\label{Blatter}
\hskip5mm \delta T_n(\rada {a_1}{a_n})\!\!&=& \!\!
\cyclsum T_{n-1}(a_1\cdot a_2,\rada{a_3}{a_n}),\ n\geq 4.
\end{eqnarray}
Equation~(\ref{Blatter}) describes a trace over a certain
$\Ass$-module $\barD$, closely related to the simplex. The
following
Section~\ref{hrnicek1} is devoted to the study of this module.

\section{Simplex as a module over the operad for associative
algebras}
\label{hrnicek1}

Let $\Delta_n$ be the standard $(n-1)$-dimensional simplex
(!observe
that the conventional notation for
our $\Delta_n$ is $\Delta^{n-1}$!). An
explicit description of $\Delta_n$ is the following. Denote, for
$1\leq i\leq n$, by $e_i$ the point
$(\rada 00,1,\rada 00)\in {\bf R}^n$ ($1$ at
the $i$-th position). Then $\Delta_n \subset {\bf R}^n$ is
the convex
hull of the set $\{\rada{e_1}{e_n}\}$.
Figure~\ref{simplex} shows $\Delta_n$ for $n=3$.
\begin{figure}[hbtp]
\begin{center}
%TexCad Options
%\grade{\off}
%\emlines{\off}
%\beziermacro{\off}
%\reduce{\on}
%\snapping{\off}
%\quality{2.00}
%\graddiff{0.01}
%\snapasp{1}
%\zoom{3.00}
\unitlength 1.20mm
\begin{picture}(34.44,46.16)
\unitlength 1.6mm
\thinlines
\put(13.83,13.16){\line(-5,-3){23.83}}
\put(10.61,13.28){\line(5,-3){23.83}}
\put(12.33,10.33){\line(0,1){25.83}}
\put(12.33,27.16){\makebox(0,0)[cc]{$\bullet$}}
\put(-4.00,2.33){\makebox(0,0)[cc]{$\bullet$}}
\put(28.83,2.33){\makebox(0,0)[cc]{$\bullet$}}
\put(-4.06,-0.62){\makebox(0,0)[cc]{$\{1\}$}}
\put(28.83,-0.78){\makebox(0,0)[cc]{$\{2\}$}}
\put(9.33,29.66){\makebox(0,0)[cc]{$\{3\}$}}
\put(1.66,17.66){\makebox(0,0)[cc]{$\{31\}$}}
\put(22.66,18.16){\makebox(0,0)[cc]{$\{23\}$}}
\put(10.83,-1.28){\makebox(0,0)[cc]{$\{12\}$}}
\put(10.50,6.33){\makebox(0,0)[cc]{$\{123\}$}}
\put(12.33,27.33){\line(-1,0){1.50}}
\put(10.83,27.33){\line(1,0){3.00}}
\put(16.17,13.16){\makebox(0,0)[cc]{$\cdot$}}
\put(3.83,11.00){\makebox(0,0)[cc]{$\cdot$}}
\put(13.00,8.83){\makebox(0,0)[cc]{$\cdot$}}
\put(7.33,11.50){\makebox(0,0)[cc]{$\cdot$}}
\put(11.33,16.66){\makebox(0,0)[cc]{$\cdot$}}
\put(8.17,12.33){\makebox(0,0)[cc]{$\cdot$}}
\put(6.67,15.83){\makebox(0,0)[cc]{$\cdot$}}
\put(20.50,9.33){\makebox(0,0)[cc]{$\cdot$}}
\put(21.50,5.83){\makebox(0,0)[cc]{$\cdot$}}
\put(17.33,5.00){\makebox(0,0)[cc]{$\cdot$}}
\put(10.50,18.00){\makebox(0,0)[cc]{$\cdot$}}
\put(14.00,18.50){\makebox(0,0)[cc]{$\cdot$}}
\put(0.83,6.16){\makebox(0,0)[cc]{$\cdot$}}
\put(-5.22,3.50){\line(1,-1){2.44}}
\put(27.56,0.94){\line(1,1){3.00}}
\put(2.56,3.83){\makebox(0,0)[cc]{$\cdot$}}
\put(17.67,16.17){\makebox(0,0)[cc]{$\cdot$}}
\put(11.00,22.94){\makebox(0,0)[cc]{$\cdot$}}
\thicklines
\put(12.33,27.33){\line(-2,-3){16.67}}
\put(-4.34,2.33){\line(1,0){33.33}}
\put(29.00,2.33){\line(-2,3){16.78}}
\end{picture}
\end{center}
\caption{$\Delta_3$.\label{simplex}}
\end{figure}
There is a classical
correspondence between the poset of subsets of $\{\rada1n\}$
and the
poset of faces of $\Delta_n$ given by
\[
\mbox{subset $S$ of $\{\rada1n\}$}
\longleftrightarrow
\mbox{convex hull of the set $\{e_i\}_{i\in S}\subset {\bf
R}^n$.}
\]
See~\cite[\S10.1]{switzer:75} for details.
Let $\Delta := \{\Delta_n\}_{\geq 1}$. In fact, it is more
convenient
to consider the symmetrized version $\barDelta :=
\{\barDelta_n\}_{\geq 1}$, where $\barDelta_n$ is the disjoint
union
of $(n-1)!$ copies of $\Delta_n$, indexed by cyclic orders
of its
vertices. This means that the poset of faces of $\barDelta_n$
consists of elements of the form
\begin{equation}
\label{bookshop}
\set{\rada{i_1}{i_l}}\times [\sigma],
\end{equation}
where $\set{\rada{i_1}{i_l}}$ is a subset of $\set{\rada
1n}$, $1\leq
l\leq n$, and $[\sigma]$ is an equivalence class from the
left coset
${\bf Z}_n \backslash \Sigma_n$. We define the
right action of $\Sigma_n$ by
\[
(\set{\rada{i_1}{i_l}}\times [\sigma])\cdot \rho :=
\set{\rada{\rho^{-1}(i_1)}{\rho^{-1}(i_l)}}\times [\sigma\rho].
\]

Let $e = \set{\rada{i_1}{i_l}}\times [\sigma]$ be a face (=
cell) of
$\barDelta(n)$ as in~(\ref{bookshop}). An orientation of $e$
is given
by choosing an order of elements of $\set{\rada{i_1}{i_l}}$. Two
such
orders induce the same orientation if and only if they differ
by a
permutation of signature $+1$. Thus the cellular cell complex
$\barD(n) := CC_*(\barDelta(n))$ is a vector space with
the basis
\begin{equation}
\label{Y}
\langle\rada{i_1}{i_l}\rangle\times [\sigma]
\end{equation}
where $\langle\rada{i_1}{i_l}\rangle$ denotes
the cell $\set{\rada{i_1}{i_l}}\times [\sigma]$ with the
orientation
induced by the order $i_1<
\cdots < i_l$. The right action of $\Sigma_n$ is given by
\[
(\langle\rada{i_1}{i_l}\rangle\times [\sigma])\cdot \rho :=
(\langle\rada{\rho^{-1}(i_1)}{\rho^{-1}(i_l)}\rangle)\times
[\sigma\rho].
\]

\begin{theorem}
\label{vicko}
The collection $\barDelta :=
\{\barDelta_n\}_{\geq 1}$ of cell complexes
has a natural structure of a (right)
module over the operad $\Ass$ for associative algebras. The
traces
over the cellular chain complex $\barD := CC_*(\barDelta)$
are the
objects described by~(\ref{Blatter}).
\end{theorem}

\noindent
{\bf Proof.}
We observed in Remark~\ref{mince}
that the action is determined by a system of `comp'
maps
\[
\circ_i : \barDelta_n
\otimes \Ass(m) \to \barDelta_{m+n-1},\ n,m\geq 1.
\]
Since $\barDelta_n$ is $\Sigma_n$-generated by
$\Delta_n =
\Delta_n \times [\id_n] \subset \barDelta_n$
(the copy corresponding to the
`normal' cyclic order $(\rada 1n)$) and $\Ass(m) =
\bk[\Sigma_n]$ is
$\Sigma_m$-generated by the identity permutation $\id_m
\in \Sigma_m$, it is enough to specify $\circ_i(t,x)$ for
$t\in \Delta_n$ and $x = \id_m$. We define $\circ_i(-,\id_m):
\Delta_n
\to \Delta_{m+n-1}$ to be the unique simplicial map such that
\[
\circ_i(\set j \times [\id_n], [\id_m]) :=
\left\{
\begin{array}{ll}
\set j,& \mbox{ for $1 \leq j \leq i$, and}
\\
\set{j+m-1},& \mbox{ for $i < j \leq n$.}
\end{array}
\right.
\]
In other words, $\circ_i(-,\id_m)$ is the canonical inclusion
$\Delta_n
\hookrightarrow \Delta_{m+n-1}$, identifying $\Delta_n$ to the
$(n-1)$-dimensional face of $\Delta_m$ corresponding to
the subset
$\{\rada 1i,\rada{i+m}{n+m-1}\}$.
The induced map of the cellular chain complex satisfies
\begin{equation}
\label{napoleon}
\circ_i(\langle \rada 1n \rangle \times [\id_n], [\id_m])=
\langle \rada 1i,\rada{i+m}{n+m-1}\rangle \times [\id_{m+n-1}].
\end{equation}

It is a straightforward verification to prove that this
really defines
an $\Ass$-action. The second part will be proved after we
formulate
Theorem~\ref{resiz}.\qed

We are going to give an algebraic characterization of
the $\Ass$-module $\barD$. To do this, we need some more or less
standard
notions, which we will also find useful later. From now on,
if not stated
otherwise, the underlying symmetric monoidal category
will be the category of (differential) graded
vector spaces.

For any collection $E= \coll E$ there exists
the {\em free operad\/} $\fr(E)$ on
$E$~\cite[page~226]{ginzburg-kapranov:DMJ94}. The operad
$\fr(E)$
has the following very explicit description in terms of trees.
Denote by $\tr$ the set of (labeled rooted)
trees and by $\tr_n$ the subset of $\tr$ consisting of
trees having $n$ input edges. Let $E(\tT)$ denote, for $\tT
\in \tr$,
the set of `multilinear' colorings of the vertices of $\tT$
by the elements of
$E$ such that a vertex with $k$ input edges is colored by
an element
of $E(k)$. The free operad $\fr(E)$ on $E$ may be then
described as
\begin{equation}
\label{fax}
\fr(E)(n):=
\bigoplus_{\tT\in\tr_n}E(\tT)
\end{equation}
with the operad structure
on $\fr(E)$ given by the operation of `grafting'
trees. We will in fact always assume that $E(1)= 0$, thus
we consider
in~(\ref{fax}) only trees whose all vertices are at
least binary, i.e.\
they have at least two incoming edges.
The details may be found
in~\cite{ginzburg-kapranov:DMJ94,getzler-jones:preprint}.

\begin{example}
\label{zajicek_usacek}{\rm\
The set $\tr_2$ has only one element (Figure~\ref{t2andt3}) and
$\fr(E)(2)= E(2)$. The set $\tr_3$ has four elements (see again
Figure~\ref{t2andt3})
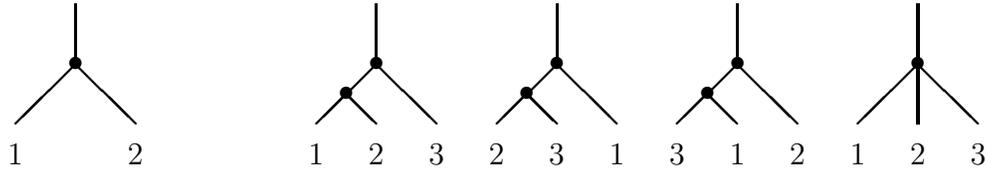
\begin{figure}[hbtp]
\begin{center}
%TexCad Options
%\grade{\off}
%\emlines{\off}
%\beziermacro{\off}
%\reduce{\on}
%\snapping{\off}
%\quality{2.00}
%\graddiff{0.01}
%\snapasp{1}
%\zoom{0.70}
\unitlength 0.80mm
\thicklines
\begin{picture}(157.62,22.95)
\put(7.62,22.95){\line(0,-1){10.00}}
\put(7.62,12.95){\line(-1,-1){10.00}}
\put(7.62,12.95){\line(1,-1){10.00}}
\put(57.62,22.95){\line(0,-1){10.00}}
\put(57.62,12.95){\line(-1,-1){10.00}}
\put(57.62,12.95){\line(1,-1){10.00}}
\put(52.62,7.95){\line(1,-1){5.00}}
\put(7.62,12.95){\makebox(0,0)[cc]{$\bullet$}}
\put(52.62,7.95){\makebox(0,0)[cc]{$\bullet$}}
\put(57.62,12.95){\makebox(0,0)[cc]{$\bullet$}}
\put(-2.38,-2.05){\makebox(0,0)[cc]{$1$}}
\put(17.62,-2.05){\makebox(0,0)[cc]{$2$}}
\put(47.62,-2.05){\makebox(0,0)[cc]{$1$}}
\put(57.62,-2.05){\makebox(0,0)[cc]{$2$}}
\put(67.62,-2.05){\makebox(0,0)[cc]{$3$}}
\put(87.62,22.95){\line(0,-1){10.00}}
\put(87.62,12.95){\line(-1,-1){10.00}}
\put(87.62,12.95){\line(1,-1){10.00}}
\put(82.62,7.95){\line(1,-1){5.00}}
\put(82.62,7.95){\makebox(0,0)[cc]{$\bullet$}}
\put(87.62,12.95){\makebox(0,0)[cc]{$\bullet$}}
\put(77.62,-2.05){\makebox(0,0)[cc]{$2$}}
\put(87.62,-2.05){\makebox(0,0)[cc]{$3$}}
\put(97.62,-2.05){\makebox(0,0)[cc]{$1$}}
\put(117.62,22.95){\line(0,-1){10.00}}
\put(117.62,12.95){\line(-1,-1){10.00}}
\put(117.62,12.95){\line(1,-1){10.00}}
\put(112.62,7.95){\line(1,-1){5.00}}
\put(112.62,7.95){\makebox(0,0)[cc]{$\bullet$}}
\put(117.62,12.95){\makebox(0,0)[cc]{$\bullet$}}
\put(107.62,-2.05){\makebox(0,0)[cc]{$3$}}
\put(117.62,-2.05){\makebox(0,0)[cc]{$1$}}
\put(127.62,-2.05){\makebox(0,0)[cc]{$2$}}
\put(147.62,22.95){\line(0,-1){10.00}}
\put(147.62,12.95){\line(-1,-1){10.00}}
\put(147.62,12.95){\line(1,-1){10.00}}
\put(147.62,12.95){\line(0,-1){10.00}}
\put(147.62,12.95){\makebox(0,0)[cc]{$\bullet$}}
\put(137.62,-2.05){\makebox(0,0)[cc]{$1$}}
\put(147.62,-2.05){\makebox(0,0)[cc]{$2$}}
\put(157.62,-2.05){\makebox(0,0)[cc]{$3$}}
\end{picture}
\end{center}
\caption{The sets $\tr_2$ (left) and $\tr_3$
(right).\label{t2andt3}}
\end{figure}
and $\fr(E)(3)$ consists of three copies of $E(2)\otimes E(2)$
which corresponds to the three binary trees in $\tr_3$ and
one copy of
$E(3)$ corresponding to the corolla (= the tree with one
vertex).
Compare also~\cite[Figure~7]{ginzburg-kapranov:DMJ94}.
}\end{example}

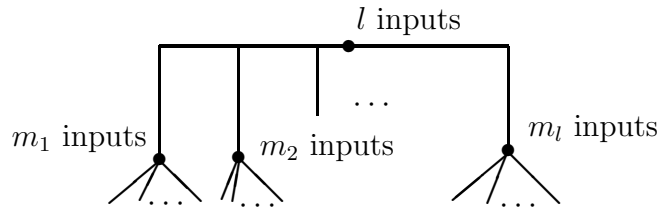
\begin{figure}[hbtp]
\begin{center}
%TexCad Options
%\grade{\off}
%\emlines{\off}
%\beziermacro{\off}
%\reduce{\on}
%\snapping{\off}
%\quality{2.00}
%\graddiff{0.01}
%\snapasp{1}
%\zoom{1.65}
\unitlength 1.00mm
\thicklines
\begin{picture}(59.69,22.89)
\put(6.69,5.39){\line(-6,-5){6.83}}
\put(6.69,5.56){\line(-1,-2){2.75}}
\put(6.69,5.23){\line(1,-1){5.33}}
\put(17.02,5.73){\line(-2,-5){2.20}}
\put(17.02,5.73){\line(1,-1){5.83}}
\put(53.02,6.39){\line(-6,-5){7.50}}
\put(52.85,6.56){\line(-2,-5){2.73}}
\put(52.69,6.73){\line(1,-1){7.00}}
\put(31.69,20.39){\makebox(0,0)[cc]{$\bullet$}}
\put(6.52,5.39){\makebox(0,0)[cc]{$\bullet$}}
\put(17.02,5.73){\makebox(0,0)[cc]{$\bullet$}}
\put(52.85,6.56){\makebox(0,0)[cc]{$\bullet$}}
\put(34.69,12.73){\makebox(0,0)[cc]{$\cdots$}}
\put(7.35,-0.44){\makebox(0,0)[cc]{$\cdots$}}
\put(19.52,-0.61){\makebox(0,0)[cc]{$\cdots$}}
\put(54.19,-0.61){\makebox(0,0)[cc]{$\cdots$}}
\put(32.71,22.89){\makebox(0,0)[lb]{$l$ inputs}}
\put(4.85,7.23){\makebox(0,0)[rb]{$m_1$ inputs}}
\put(19.85,6.06){\makebox(0,0)[lb]{$m_2$ inputs}}
\put(55.52,8.56){\makebox(0,0)[lb]{$m_l$ inputs}}
\put(6.59,5.41){\line(0,1){15.19}}
%\put(17.13,20.59){\line(0,0){0.75}}
\put(52.90,6.41){\line(0,1){14.18}}
\put(52.90,20.59){\line(-1,0){46.32}}
\put(27.50,20.62){\line(0,-1){9.30}}
\put(17.19,5.53){\line(-1,-6){0.92}}
\put(17.02,20.56){\line(0,-1){15.17}}
\end{picture}
\end{center}
\caption{The tree $\tT_{m_1,\ldots,m_l}$. The output vertex is
symbolized as a `rake' with no output edge,
to underline its distinguished character.
\label{thetree}}
\end{figure}

In the same manner, for each operad $\P$ and for each collection
$X =
\coll X$ there exists the {\em free (right) $\P$-module\/}
generated
by the collection $X$, which we denote $X\circ \P$.
An explicit description is~\cite[page~312]{markl:dl}
\begin{equation}
\label{cajicek}
(X\circ \P)(m) =
\def\arraystretch{.7}
\bigoplus
\left(
\mbox{Ind}^{\Sigma_m}_{\Sigma_{m_1}
\times \cdots \times \Sigma_{m_l}}
(X(l)\otimes \P(m_1)\otimes \cdots \otimes \P(m_l))
\right)_{\Sigma_l},
\end{equation}
where the summation is taken over all $m_1+\cdots +m_l =n$,
$l\geq 1$.
On the right-hand side,
$\mbox{Ind}^{\Sigma_m}_{\Sigma_{k_1}
\times \cdots \times \Sigma_{k_l}}(-)$ denotes the induced
representation and $(-)_{\Sigma_l}$ the quotient under the
obvious
action of $\Sigma_l$.
The term $X(l)\otimes \P(m_1)\otimes \cdots \otimes \P(m_l)$
on the right-hand side of~(\ref{cajicek})
can be interpreted as colorings of the tree
$\tT_{m_1,\ldots,m_l}$
from Figure~\ref{thetree} such that the output vertex is colored
by an element of $X(l)$ and
the remaining vertices by elements of $\P$.

\begin{example}
\label{haficek}
{\rm\
We have $(X\circ \P)(1) = X(1)$, corresponding to the tree
$\tT_1$ on
Figure~\ref{odrazka}.
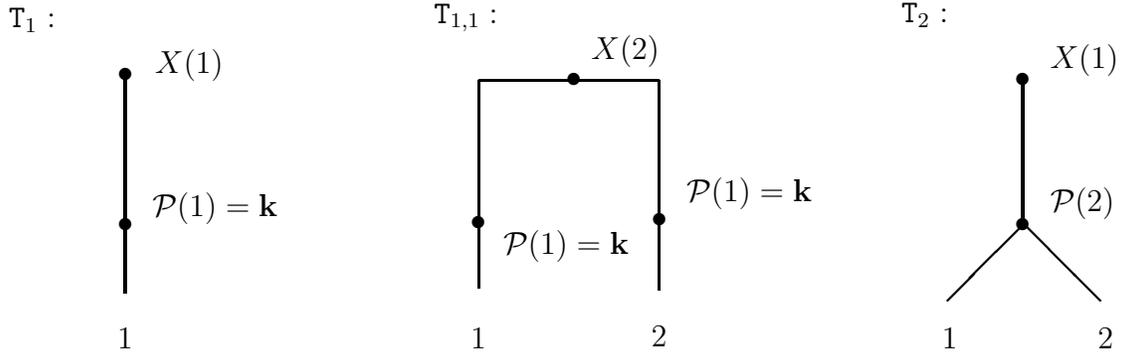
\begin{figure}[hbtp]
\begin{center}
%TexCad Options
%\grade{\off}
%\emlines{\off}
%\beziermacro{\off}
%\reduce{\on}
%\snapping{\off}
%\quality{0.20}
%\graddiff{0.01}
%\snapasp{1}
%\zoom{1.00}
\unitlength 1.00mm
\thicklines
\begin{picture}(142.67,43.00)
\put(131.67,15.00){\line(-1,-1){10.00}}
\put(131.67,15.33){\line(1,-1){10.33}}
\put(12.34,35.00){\makebox(0,0)[cc]{$\bullet$}}
\put(72.01,34.38){\makebox(0,0)[cc]{$\bullet$}}
\put(131.67,34.33){\makebox(0,0)[cc]{$\bullet$}}
\put(12.34,15.00){\makebox(0,0)[cc]{$\bullet$}}
\put(59.29,15.33){\makebox(0,0)[cc]{$\bullet$}}
\put(83.34,15.67){\makebox(0,0)[cc]{$\bullet$}}
\put(131.67,15.00){\makebox(0,0)[cc]{$\bullet$}}
\put(-0.00,42.33){\makebox(0,0)[cc]{$\tT_1:$}}
\put(57.67,42.66){\makebox(0,0)[cc]{$\tT_{1,1}:$}}
\put(118.67,43.00){\makebox(0,0)[cc]{$\tT_2:$}}
\put(16.34,36.33){\makebox(0,0)[lc]{$X(1)$}}
\put(16.00,17.33){\makebox(0,0)[lc]{$\P(1)=\bk$}}
\put(74.36,38.21){\makebox(0,0)[lc]{$X(2)$}}
\put(62.67,12.33){\makebox(0,0)[lc]{$\P(1)=\bk$}}
\put(87.01,19.33){\makebox(0,0)[lc]{$\P(1)=\bk$}}
\put(135.33,37.00){\makebox(0,0)[lc]{$X(1)$}}
\put(135.33,18.00){\makebox(0,0)[lc]{$\P(2)$}}
\put(12.34,0.00){\makebox(0,0)[cc]{$1$}}
\put(59.34,0.00){\makebox(0,0)[cc]{$1$}}
\put(83.34,0.33){\makebox(0,0)[cc]{$2$}}
\put(122.00,0.00){\makebox(0,0)[cc]{$1$}}
\put(142.67,0.00){\makebox(0,0)[cc]{$2$}}
\put(12.34,35.00){\line(0,-1){29.00}}
\put(59.34,34.33){\line(1,0){24.00}}
\put(83.34,34.33){\line(0,-1){28.00}}
\put(59.34,34.33){\line(0,-1){27.67}}
\put(131.67,34.33){\line(0,-1){19.33}}
\end{picture}
\end{center}
\caption{$(X\circ \P)(1)$ (left) and $(X\circ \P)(2)$
(right).\label{odrazka}}
\end{figure}
The vector space $(X\circ \P)(2)$ consists of a
copy of $X(2)$ corresponding to $\tT_{1,1}$ and a copy of
$X(1)\otimes \P(2)$ corresponding to $\tT_2$, see again
Figure~\ref{odrazka}. Note that we still assume that $\P(1)=
\bk$.
}\end{example}

For a graded vector space $V=\bigoplus_p V_p$ let $\susp V$
(resp.
$\desusp V$) be the {\em suspension\/} (resp. the {\em
desuspension\/}) of $V$, i.e.
the graded
vector space defined by $(\susp V)_p := V_{p-1}$
(resp. $(\desusp V)_p
:= V_{p+1}$).
We have the obvious natural maps $\uparrow : V \to \susp V$ and
$\downarrow: V\to \desusp V$.
For a collection $E$, the {\em suspension\/} $\ss E$
is the collection with
\begin{equation}
\label{sign-factor}
(\ss E)(n):=
\sgn\ \otimes \uparrow^{n-1}E(n),
\end{equation}
$n\geq 1$, where
$\uparrow^{n-1}$ is the $(n-1)$-fold
suspension introduced above and $\sgn$ is the signum
representation of
$\Sigma_n$ on $\bk$. The reason why we need the signum
factor is that we intend to apply the suspension to operads and
modules over operads.
Without this factor, the composition induced on
the suspension will
not be equivariant,
compare also~\cite[page~8]{getzler-jones:preprint}. There is an
obvious similar notion of the {\em desuspension\/} $\ss^{-1}
E$ of the
collection $E$.

Let us come back to the promised algebraic characterization
of the
$\Ass$-module $\barD$.
Consider the free $\Ass$-module $\ss\Cycl \circ \Ass$
generated by
the suspension of the collection $\Cycl$ introduced in
Example~\ref{tuzka}.

It will be useful
to have an explicit description of
the elements of $\ss \Cycl \circ
\Ass$. Consider the free graded right $\Sigma_n$-module $H(n)$
generated by the trees $\tT_{m_1,\ldots,m_l}$ introduced in
Figure~\ref{thetree}, with $m_1+\cdots+m_l=n$.
The grading is given by $\deg(\tT_{m_1,\ldots,m_l}):= l-1$
Thus $H(n)$ is, by
definition, the graded vector space with the basis
\[
\left\{
\tT_{m_1,\ldots,m_l}\times \sigma,\ \sigma \in \Sigma_n,\
1\leq l\leq n,\ m_1+\cdots+m_l=n
\right\},\ \deg(\tT_{m_1,\ldots,m_l}\times \sigma)= l-1.
\]
A neat graphical presentation of the symbol
$\tT_{m_1,\ldots,m_l}\times
\sigma$ is the tree $\tT_{m_1,\ldots,m_l}$ with the inputs
labeled by
$\rada{\sigma^{-1}(1)}{\sigma^{-1}(n)}$. For example,
${\tT_{2,1}}\times (123)$ can be depicted as
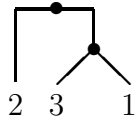
\begin{center}
\unitlength 1.20mm
\thicklines
\begin{picture}(10,13)(53,23)
\put(52.63,33.53){\line(1,0){8.67}}
\put(61.30,33.53){\line(0,-1){4.33}}
\put(61.30,29.19){\line(-1,-1){4.00}}
\put(61.30,29.19){\line(1,-1){4.00}}
\put(52.63,33.53){\line(0,-1){8.00}}
\put(52.66,22.95){\makebox(0,0)[cc]{2}}
\put(57.19,22.95){\makebox(0,0)[cc]{3}}
\put(57.19,33.53){\makebox(0,0)[cc]{$\bullet$}}
\put(61.35,28.95){\makebox(0,0)[cc]{$\bullet$}}
\put(65.21,22.95){\makebox(0,0)[cc]{1}}
\end{picture}
\end{center}

Define  now the left action of the group $\ZZ_l$, considered
as the group
of cyclic permutations of order $l$, on $H(n)$ as
follows. Recall that, for $\zeta \in \ZZ_l \subset \Sigma_l$,
we denoted by $\zeta(\rada{m_1}{m_l}) \in
\Sigma_n$ the permutation which permutes the blocks of
$\rada{m_1}{m_l}$-elements via $\zeta$. Then we put
\begin{equation}
\label{pejska_Mikinka}
\zeta
\left(
\tttr{m_1}{m_l} \times \sigma
\right)
:= \sgn(\zeta) \cdot
\tttr{m_{\zeta^{-1}(1)}}{m_{\zeta^{-1}(l)}} \times
\zeta(\rada{m_1}{m_l})\cdot\sigma.
\end{equation}
The following lemma is an easy consequence of definitions and
formula~(\ref{cajicek}).

\begin{lemma}
\label{MGD}
The graded vector space $(\ss\Cycl
\circ \Ass)(n)$ can be identified with the graded vector space
with the
basis given by equivalence classes of symbols
\begin{equation}
\label{3M}
\tttr{m_1}{m_l} \times \sigma \in H(n),\
\deg(\tT_{m_1,\ldots,m_l})= l-1,
\end{equation}
modulo the left action of the group $\ZZ_l$ defined
in~(\ref{pejska_Mikinka}). Under this identification,
the right action of $\Sigma_n$ on the
equivalence class $[\tttr{m_1}{m_l} \times \sigma]$ is
described as
\[
[\tttr{m_1}{m_l} \times \sigma]\cdot \rho :=
\sgn(\sigma)\cdot [\tttr{m_1}{m_l} \times \sigma\rho].
\]
\end{lemma}

Let $\xi_n\in \Cycl(n)$ be the generator
represented by the identical permutation $\id_n\in \Sigma_n$
and let
$\alpha_2 = \id_2 \in \Ass(2) = \bk[\Sigma_2]$. Define the
differential $\partial$ on $\ss\Cycl \circ \Ass$ by
\begin{eqnarray}
\label{guma}
\partial (\susp^{n-1} \xi_n) &:=& -\sum_{\sigma}
\sgn(\sigma)\cdot \nu(\susp^{n-2}
\xi_{n-1};\alpha_2,\rada 11)\cdot \sigma
\\
\nonumber
&=&-\cyclsum
\nu(\susp^{n-2}\xi_{n-1};\alpha_2,\rada 11).\
\mbox{ /the cyclic sum notation of~(\ref{cyclsum})/}
\end{eqnarray}
Using the identification of
Lemma~\ref{MGD}, this could be also written as
\[
\partial ({\tT}_{\underbrace%
{\mbox{\scriptsize $1,\ldots,1$}}_{n\times}}
\times [\id_n])= -\cyclsum
({\tT}_{2,\underbrace%
{\mbox{\scriptsize $1,\ldots,1$}}_{n-1\times}}\times[\id_n] ),
\]
or, in a diagrammatic shorthand,
\vskip2mm
\begin{center}
%TexCad Options
%\grade{\off}
%\emlines{\off}
%\beziermacro{\off}
%\reduce{\on}
%\snapping{\off}
%\quality{2.00}
%\graddiff{0.01}
%\snapasp{1}
%\zoom{2.00}
\unitlength 1.30mm
\thicklines
\begin{picture}(69.67,5.67)
\put(10.00,5.67){\line(0,0){0.00}}
\put(10.00,5.67){\line(1,0){15.00}}
\put(25.00,5.67){\line(0,-1){4.00}}
\put(10.00,1.67){\line(0,1){4.00}}
\put(13.00,1.67){\line(0,1){4.00}}
\put(8.33,3.67){\makebox(0,0)[rc]{$\partial($}}
\put(26.83,3.67){\makebox(0,0)[lc]{$)$}}
\put(30.67,3.67){\makebox(0,0)[cc]{$=$}}
\put(17.67,5.67){\makebox(0,0)[cc]{$\bullet$}}
\put(10.00,-0.16){\makebox(0,0)[cc]{$1$}}
\put(13.00,-0.16){\makebox(0,0)[cc]{$2$}}
\put(16.00,-0.16){\makebox(0,0)[cc]{$3$}}
\put(25.00,-0.16){\makebox(0,0)[cc]{$n$}}
\put(20.33,2.17){\makebox(0,0)[cc]{$\cdots$}}
\put(16.00,1.67){\line(0,1){4.00}}
\put(16.00,5.67){\line(0,-1){4.00}}

\put(-10,0){
\put(54.67,5.67){\line(1,0){15.00}}
\put(69.67,5.67){\line(0,-1){4.00}}
\put(54.67,5.67){\line(0,-1){1.00}}
\put(59.00,5.67){\line(0,-1){4.00}}
\put(49.67,3.67){\makebox(0,0)[rc]{$-\ \displaystyle\cyclsum$}}
\put(64.50,2.34){\makebox(0,0)[cc]{$\cdots$}}
\put(54.67,4.59){\line(-2,-5){1.20}}
\put(54.67,4.50){\line(1,-3){0.94}}
\put(55.61,1.67){\line(0,0){0.06}}
\put(53.17,0.00){\makebox(0,0)[cc]{$1$}}
\put(56.17,0.00){\makebox(0,0)[cc]{$2$}}
\put(59.17,0.00){\makebox(0,0)[cc]{$3$}}
\put(69.50,0.00){\makebox(0,0)[cc]{$n$}}
\put(62.33,5.67){\makebox(0,0)[cc]{$\bullet$}}
}
\end{picture}
\end{center}
Since the elements $\{ \susp^{n-1} \xi_n\}_{n\geq 1}$ generate
$\ss\Cycl \circ \Ass$, formula~(\ref{guma}) is enough to
determine
the differential $\partial$.
We leave to the reader to verify that the definition is
correct and
that $\partial^2=0$.

\begin{theorem}
\label{resiz}
The $\Ass$-module $\barD = CC_*(\barDelta)$
is isomorphic to the free differential
$\Ass$-module $(\ss\Cycl \circ \Ass,\partial)$ constructed
above.
\end{theorem}

\noindent
{\bf Proof.}
Any differential $\Ass$-module map $\omega:
(\ss\Cycl \circ \Ass,\partial) \to (\barD,\partial_{\barD})$ is
determined by the values $\omega(\susp^{n-1}\xi_n)$, $n\geq 1$.
We define $\omega(\susp^{n-1}\xi_n):= e_n$, where $e_n \in
\barD(n)$
is the top $(n-1)$-dimensional oriented cell
$\langle \rada 1n\rangle$.
We shall verify that the above defined map $\omega$ commutes
with the
differentials,
\[
\omega(\partial \susp^{n-1}\xi_n) = -\omega(\cyclsum
\nu(\susp^{n-2}\xi_{n-1},\alpha_2,\rada11)) =
\partial e_n.
\]
Because $\omega$ is a module homomorphism,
the above equation can be rewritten as
\begin{equation}
\label{hvezdicka}
\cyclsum\nu_{\barD}(e_{n-1};\alpha_2,\rada11) = -\partial e_n.
\end{equation}
The standard formula for the boundary of
$\langle \rada 1n \rangle$~\cite[\S10.1]{switzer:75}
says that
\begin{equation}
\label{ja}
\partial e_n = \partial \langle\rada 1n\rangle
= \sum_{1\leq i \leq n}
\znamenko{i+1}
\langle  \rada 1{i-1},\rada {i+1}n\rangle
\end{equation}
while the defining formula~(\ref{napoleon})
for the $\Ass$-module action on $\barDelta$ gives
\[
\nu_{\barD}(e_{n-1};\alpha_2,\rada11) =
\langle1,\rada3n\rangle.
\]
Now it is enough to observe that
\[
\cyclsum\langle1,\rada3n\rangle= - \sum_{1\leq i
\leq n}\znamenko{i+1}\langle\rada1{i-1},\rada{i+1}n\rangle
\]
which, together with~(\ref{ja}), gives~(\ref{hvezdicka}).

\begin{figure}[hbtp]
\begin{center}
%TexCad Options
%\grade{\off}
%\emlines{\off}
%\beziermacro{\off}
%\reduce{\on}
%\snapping{\off}
%\quality{2.00}
%\graddiff{0.01}
%\snapasp{1}
%\zoom{1.00}
\unitlength 1.20mm
\thicklines
\begin{picture}(91.08,94.33)
\put(-11.20,94.32){\makebox(0,0)[rc]{\fbox{$n=1$:}}}
\put(-1.37,94.15){\makebox(0,0)[cc]{$\bullet$}}
\put(-1.37,86.65){\makebox(0,0)[cc]{$1$}}
\put(-11.20,77.48){\makebox(0,0)[rc]{\fbox{$n=2$:}}}
\put(12.30,77.48){\line(1,0){5.00}}
\put(17.30,77.48){\line(0,-1){5.00}}
\put(12.13,77.48){\line(0,-1){5.00}}
\put(14.72,77.48){\makebox(0,0)[cc]{$\bullet$}}
\put(-3.87,77.48){\line(1,0){5.00}}
\put(1.13,77.48){\line(0,-1){5.00}}
\put(-4.04,77.48){\line(0,-1){5.00}}
\put(-1.45,77.48){\makebox(0,0)[cc]{$\bullet$}}
\put(-4.04,69.65){\makebox(0,0)[cc]{$1$}}
\put(1.13,69.65){\makebox(0,0)[cc]{$2$}}
\put(12.13,69.65){\makebox(0,0)[cc]{$2$}}
\put(17.30,69.65){\makebox(0,0)[cc]{$1$}}
\put(7.97,77.48){\makebox(0,0)[cc]{$=-$}}
\put(20.79,77.48){\makebox(0,0)[lc]{$=\langle12\rangle,$}}
\put(42.46,78.32){\makebox(0,0)[cc]{$\bullet$}}
\put(42.46,74.65){\makebox(0,0)[cc]{$\bullet$}}
\put(45.84,76.48){\makebox(0,0)[lc]{$= \langle1\rangle,$}}
\put(61.76,78.32){\makebox(0,0)[cc]{$\bullet$}}
\put(61.63,74.65){\makebox(0,0)[cc]{$\bullet$}}
\put(65.01,76.48){\makebox(0,0)[lc]{$=-\langle2\rangle$}}
\put(40.13,68.30){\makebox(0,0)[cc]{$1$}}
\put(44.96,68.30){\makebox(0,0)[cc]{$2$}}
\put(59.30,68.30){\makebox(0,0)[cc]{$2$}}
\put(64.30,68.30){\makebox(0,0)[cc]{$1$}}
\put(-11.20,56.52){\makebox(0,0)[rc]{\fbox{$n=3$:}}}
\put(-6.37,56.69){\line(1,0){10.00}}
\put(3.63,56.69){\line(0,-1){5.17}}
\put(-6.37,56.69){\line(0,-1){5.17}}
\put(-6.37,51.52){\line(0,0){0.00}}
\put(-1.37,56.69){\makebox(0,0)[cc]{$\bullet$}}
\put(11.13,56.69){\line(1,0){10.00}}
\put(28.63,56.69){\line(1,0){10.00}}
\put(21.13,56.69){\line(0,-1){5.17}}
\put(38.63,56.69){\line(0,-1){5.17}}
\put(11.13,56.69){\line(0,-1){5.17}}
\put(28.63,56.69){\line(0,-1){5.17}}
\put(11.13,51.52){\line(0,0){0.00}}
\put(28.63,51.52){\line(0,0){0.00}}
\put(16.13,56.69){\makebox(0,0)[cc]{$\bullet$}}
\put(33.63,56.69){\makebox(0,0)[cc]{$\bullet$}}
\put(7.46,57.35){\makebox(0,0)[cc]{$=$}}
\put(25.13,57.52){\makebox(0,0)[cc]{$=$}}
\put(42.80,57.52){\makebox(0,0)[lc]{$=\langle123\rangle$}}
\put(-6.37,49.02){\makebox(0,0)[cc]{$1$}}
\put(-1.37,49.02){\makebox(0,0)[cc]{$2$}}
\put(3.63,49.02){\makebox(0,0)[cc]{$3$}}
\put(11.13,49.02){\makebox(0,0)[cc]{$2$}}
\put(16.13,49.02){\makebox(0,0)[cc]{$3$}}
\put(21.13,49.02){\makebox(0,0)[cc]{$1$}}
\put(28.63,49.02){\makebox(0,0)[cc]{$3$}}
\put(33.63,49.02){\makebox(0,0)[cc]{$1$}}
\put(38.63,49.02){\makebox(0,0)[cc]{$2$}}
\put(-6.37,33.53){\line(1,0){8.67}}
\put(2.30,33.53){\line(0,-1){4.33}}
\put(2.30,29.19){\line(-1,-1){4.00}}
\put(2.30,29.19){\line(1,-1){4.00}}
\put(-6.37,33.53){\line(0,-1){8.00}}
\put(2.30,28.86){\makebox(0,0)[cc]{$\bullet$}}
\put(27.96,33.53){\line(-1,0){8.67}}
\put(19.29,33.53){\line(0,-1){4.33}}
\put(19.29,29.19){\line(1,-1){4.00}}
\put(19.29,29.19){\line(-1,-1){4.00}}
\put(27.96,33.53){\line(0,-1){8.00}}
\put(19.29,28.86){\makebox(0,0)[cc]{$\bullet$}}
\put(11.96,33.53){\makebox(0,0)[cc]{$=-$}}
\put(52.63,33.53){\line(1,0){8.67}}
\put(61.30,33.53){\line(0,-1){4.33}}
\put(61.30,29.19){\line(-1,-1){4.00}}
\put(61.30,29.19){\line(1,-1){4.00}}
\put(52.63,33.53){\line(0,-1){8.00}}
\put(61.30,28.86){\makebox(0,0)[cc]{$\bullet$}}
\put(86.96,33.53){\line(-1,0){8.67}}
\put(78.29,33.53){\line(0,-1){4.33}}
\put(78.29,29.19){\line(1,-1){4.00}}
\put(78.29,29.19){\line(-1,-1){4.00}}
\put(86.96,33.53){\line(0,-1){8.00}}
\put(78.29,28.86){\makebox(0,0)[cc]{$\bullet$}}
\put(70.96,33.53){\makebox(0,0)[cc]{$=-$}}
\put(-6.37,12.19){\line(1,0){8.67}}
\put(2.30,12.19){\line(0,-1){4.33}}
\put(2.30,7.86){\line(-1,-1){4.00}}
\put(2.30,7.86){\line(1,-1){4.00}}
\put(-6.37,12.19){\line(0,-1){8.00}}
\put(2.30,7.52){\makebox(0,0)[cc]{$\bullet$}}
\put(27.96,12.19){\line(-1,0){8.67}}
\put(19.29,12.19){\line(0,-1){4.33}}
\put(19.29,7.86){\line(1,-1){4.00}}
\put(19.29,7.86){\line(-1,-1){4.00}}
\put(27.96,12.19){\line(0,-1){8.00}}
\put(19.29,7.52){\makebox(0,0)[cc]{$\bullet$}}
\put(11.96,12.19){\makebox(0,0)[cc]{$=-$}}
\put(42.40,74.85){\line(-1,-2){2.32}}
\put(42.40,74.85){\line(1,-2){2.32}}
\put(61.67,74.85){\line(-2,-5){1.91}}
\put(61.67,74.70){\line(3,-5){2.78}}
\put(64.45,70.07){\line(0,0){0.06}}
\put(-1.37,94.33){\line(0,-1){5.33}}
\put(42.50,74.73){\line(0,1){3.73}}
\put(61.70,74.73){\line(0,1){3.73}}
\put(-1.37,56.73){\line(0,-1){5.33}}
\put(16.10,56.73){\line(0,-1){5.33}}
\put(33.56,56.73){\line(0,-1){5.33}}
\put(-1.90,33.53){\makebox(0,0)[cc]{$\bullet$}}
\put(23.63,33.53){\makebox(0,0)[cc]{$\bullet$}}
\put(57.03,33.53){\makebox(0,0)[cc]{$\bullet$}}
\put(82.69,33.53){\makebox(0,0)[cc]{$\bullet$}}
\put(-1.90,12.20){\makebox(0,0)[cc]{$\bullet$}}
\put(23.63,12.20){\makebox(0,0)[cc]{$\bullet$}}
\put(52.66,22.95){\makebox(0,0)[cc]{2}}
\put(57.19,22.95){\makebox(0,0)[cc]{3}}
\put(65.21,22.95){\makebox(0,0)[cc]{1}}
\put(74.34,22.95){\makebox(0,0)[cc]{3}}
\put(82.36,22.95){\makebox(0,0)[cc]{1}}
\put(86.89,22.95){\makebox(0,0)[cc]{2}}
\put(91.08,33.58){\makebox(0,0)[lc]{$= \langle23\rangle$}}
\put(-6.44,23.12){\makebox(0,0)[cc]{1}}
\put(-1.73,23.12){\makebox(0,0)[cc]{2}}
\put(6.29,23.12){\makebox(0,0)[cc]{3}}
\put(15.24,23.12){\makebox(0,0)[cc]{2}}
\put(23.26,23.12){\makebox(0,0)[cc]{3}}
\put(27.97,23.12){\makebox(0,0)[cc]{1}}
\put(-6.44,2.03){\makebox(0,0)[cc]{3}}
\put(-1.73,2.03){\makebox(0,0)[cc]{1}}
\put(6.29,2.03){\makebox(0,0)[cc]{2}}
\put(15.24,2.03){\makebox(0,0)[cc]{1}}
\put(23.44,2.03){\makebox(0,0)[cc]{2}}
\put(27.97,2.03){\makebox(0,0)[cc]{3}}
\put(30.93,33.76){\makebox(0,0)[lc]{$=\langle12\rangle$,}}
\put(30.93,12.14){\makebox(0,0)[lc]{$=-\langle13\rangle$}}
\put(6.50,94.33){\makebox(0,0)[cc]{$= \langle1\rangle$}}
\end{picture}
\end{center}
\caption{
A representation of oriented faces of $\Delta_1$, $\Delta_2$
and $\Delta_3$ by
equivalence classes of labeled planar trees, representing
elements of $H(1)$, $H(2)$ and $H(3)$.
\label{peniz}}
\end{figure}
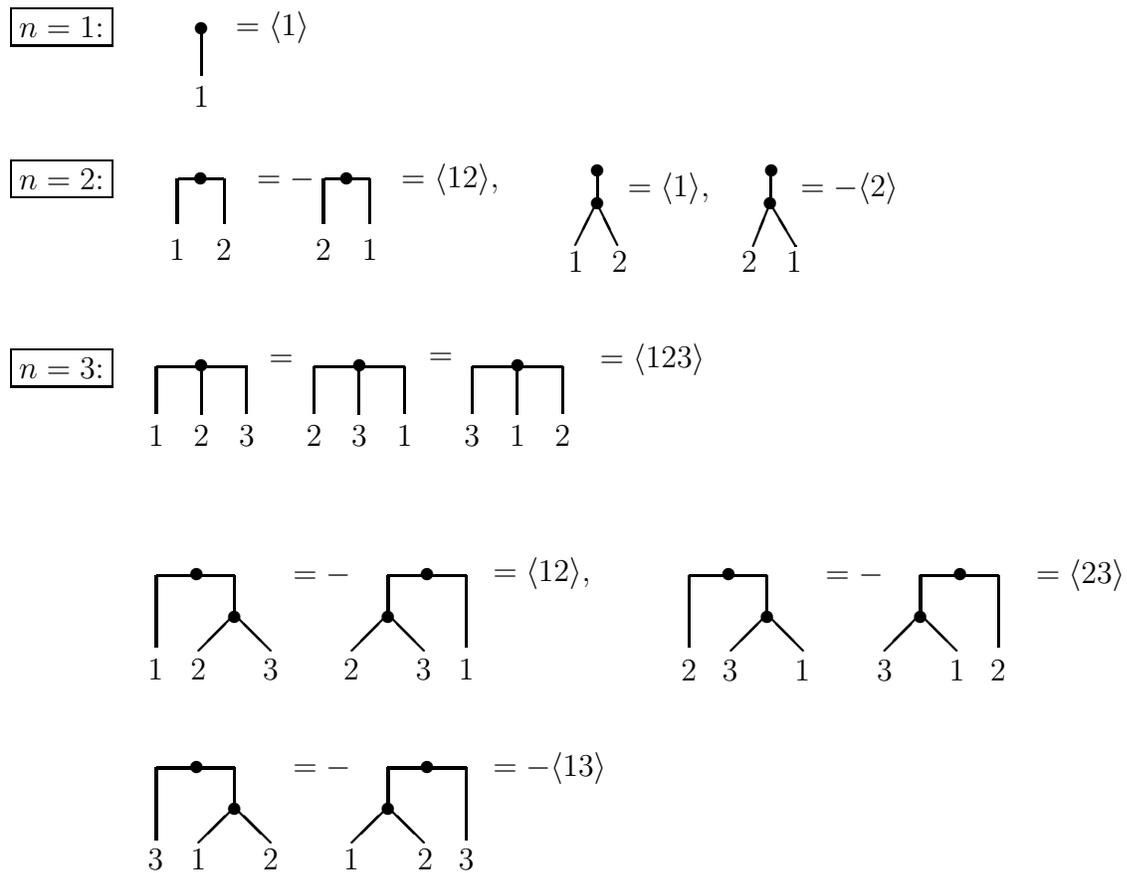

It remains to prove that $\omega$ is an isomorphism.
To this end, we give an explicit formula for the map
$\omega$. Let
$\tttr{m_1}{m_l} \times \sigma \in H(n)$ be as
in Lemma~\ref{MGD}. The
numbers $\rada{m_1}{m_l}$ determine a sequence $\rada{i_1}{i_l}$
by
$i_s:= m_1+\cdots m_{s-1}+1$, $1\leq s\leq
l$. Consider a map $\varphi: H(n)\to \barD(n)$ defined by
\[
\varphi(\tttr{m_1}{m_l} \times \sigma)
:= \langle\rada{\sigma^{-1}(i_1)}{\sigma^{-1}(i_l)}\rangle
\times [\sigma],
\]
where we denoted elements of $\barD(n)$ (= cells of
$\barDelta(n)$) as
in~(\ref{Y}).
It is immediate to see that $\varphi$ is an
$\Sigma_n$-equivariant
epimorphism. For the left action of $\zeta \in \ZZ_l$ we have
\begin{eqnarray*}
\varphi(\zeta(\tttr{m_1}{m_l} \times \sigma)
&=&
\sgn(\zeta)\!\cdot\!
\varphi(
\tttr{m_{\zeta^{-1}(1)}}%
{m_{\zeta^{-1}(l)}} \times \zeta(\rada{m_1}{m_l})\sigma)
\\
&=& \sgn(\zeta)\!\cdot\! \langle
\rada{\sigma^{-1}(i_{\zeta(1)})}{\sigma^{-1}(i_{\zeta(l)})}
\rangle \times
[\sigma]
= \langle\rada{\sigma^{-1}(i_1)}{\sigma^{-1}(i_l)}
\rangle\times
[\sigma],
\end{eqnarray*}
which shows that $\varphi(x) = \varphi(\zeta y)$. On the
other hand,
a moment's
reflection show that $\varphi(x) = \varphi(y)$, for $x,y\in
H(n)$,
implies the existence of some  $\zeta \in \ZZ_l$ such that
$x = \zeta y$.
Thus the map $\varphi$ induces an equivariant
isomorphism $(\ss\Cycl \circ \Ass)(n) = \ZZ_l \backslash H(n)
\cong \barD(n)$, which is exactly our map $\omega$. The nature
of the
map $\omega$ is illustrated on Figure~\ref{peniz}.\qed

A $\barD$-trace $t:\barD \to \E_{A,V}$ is, under the
identification
of Theorem~\ref{resiz}, given by the values $T_n :=
t(\susp^{n-1}\xi_n)$, $n\geq 1$. The axiom~(\ref{Blatter}) then
reflects~(\ref{hvezdicka}).
The relation between the axiom~(\ref{Blatter}) and the
geometry of
the simplex is visualized on Figure~\ref{delta}.
\begin{figure}[hbtp]
\begin{center}
%TexCad Options
%\grade{\off}
%\emlines{\off}
%\beziermacro{\off}
%\reduce{\on}
%\snapping{\off}
%\quality{2.00}
%\graddiff{0.01}
%\snapasp{1}
%\zoom{1.00}
\unitlength 1.00mm
\linethickness{0.4pt}
\begin{picture}(74.16,52.17)
\put(-40.33,49.67){\makebox(0,0)[lc]{$\delta T_1 (a) = 0$:}}
\put(-13.00,49.67){\makebox(0,0)[rc]{$\delta($}}
\put(-10.50,49.67){\makebox(0,0)[cc]{$\bullet$}}
\put(-10.50,52.17){\makebox(0,0)[cc]{$a$}}
\put(-8.00,49.67){\makebox(0,0)[lc]{$)=0$}}
\put(-40.33,39.67){\makebox(0,0)[lc]
{$\delta T_2 (a,b) = T_1(a\cdot b)-\znamenko{|a|\cdot
|b|}T_1(b\cdot a):$}}
\put(37.34,39.67){\makebox(0,0)[lc]{$\delta($}}
\put(43.00,39.67){\makebox(0,0)[cc]{$\bullet$}}
\put(53.00,39.67){\makebox(0,0)[cc]{$\bullet$}}
\put(55.50,39.67){\makebox(0,0)[lc]{$)=$}}
\put(64.66,39.67){\makebox(0,0)[cc]{$\bullet$}}
\put(64.66,43.17){\makebox(0,0)[cc]{$ab$}}
\put(69.16,39.84){\makebox(0,0)[cc]{$-$}}
\put(74.16,39.67){\makebox(0,0)[cc]{$\bullet$}}
\put(74.16,43.17){\makebox(0,0)[cc]{$ba$}}
\put(43.00,39.67){\vector(1,0){10.00}}
\put(-40.50,26.34){\makebox(0,0)[lc]
{$\delta T_3 (a,b,c) = T_2(a\cdot
b,c)+\znamenko{|a|\cdot(|b|+|c|)}
T_2(b\cdot c,a)+\znamenko{|c|\cdot(|a|+|b|)} T_2(c\cdot a,b)$:}}
\put(-40.50,8.67){\makebox(0,0)[lc]{$\delta\left(\rule{0mm}{10mm}\right.$}}
\put(-22.65,14.50){\vector(-3,-4){8.13}}
\put(-14.53,3.67){\vector(-3,4){8.13}}
\put(-30.80,3.67){\vector(1,0){16.25}}
\put(-7.17,8.84){\makebox(0,0)[cc]{$\left.\rule{0mm}{10mm}\right)=$}}
%\put(0.49,8.84){\makebox(0,0)[cc]{$-$}}
\put(4.49,8.84){\vector(1,0){14.00}}
\put(21.99,8.84){\makebox(0,0)[cc]{$+$}}
\put(42.49,8.84){\makebox(0,0)[cc]{$+$}}
\put(24.99,8.84){\vector(1,0){14.00}}
\put(45.49,8.84){\vector(1,0){14.00}}
\put(-26.45,10.67){\makebox(0,0)[rb]{$(ca)b$}}
\put(-18.61,10.34){\makebox(0,0)[lb]{$(bc)a$}}
\put(-22.78,2.84){\makebox(0,0)[ct]{$(ab)c$}}
\put(11.50,11.50){\makebox(0,0)[cb]{$(ab)c$}}
\put(32.00,11.67){\makebox(0,0)[cb]{$(bc)a$}}
\put(52.00,11.67){\makebox(0,0)[cb]{$(ca)b$}}
\put(42.99,43.17){\makebox(0,0)[cc]{$ab$}}
\put(53.16,43.00){\makebox(0,0)[cc]{$ba$}}
%\put(3.17,8.84){\makebox(0,0)[cc]{$($}}
%\put(61.67,8.84){\makebox(0,0)[cc]{$)$}}
\end{picture}
\end{center}
\caption{
Relation of the axiom~(9)
to the geometry of the
simplex.
\label{delta}
}
\end{figure}

\section{Quadratic operads and modules; modules associated
to cyclic
operads}
\label{cervena-tuzka}

Each operad $\P$ can be presented as a quotient $\fr(E)/I$,
for a
collection $E$ and an `operadic' ideal $I$. The operad $\P$
is {\em quadratic\/}~\cite[page~228]{ginzburg-kapranov:DMJ94}
if it has a presentation
such that the
collection $E$ is concentrated in degree $2$, $E = E(2)$,
and the
ideal $I$ is generated by a subspace $R \subset \fr(E)(3)$. In
this
case we write $\P = \prez ER$.

\begin{example}{\rm\
\label{genius}
Quadratic operads are omnipresent. Just
recall that $\Ass = \prez ER$ for $E = E(2) = \bk[\Sigma_2]$,
the
regular representation of $\Sigma_2$. Choosing a generator
$\mu \in E$, we can write $E = \Span(\mu,\mu S_{21})$, where
$S_{21}\in \Sigma_2$
is the transposition. Then $R\subset \fr(E)(3)$ is the
$\Sigma_3$-subspace generated by $\mu(1,\mu)- \mu(\mu,1)$. If we
think of $\mu$ as corresponding to a multiplication, then the
generator of $R$ expresses the associativity. Sometimes
we simplify
the notation and write $\Ass= \prez \mu{\mu(1,\mu)-
\mu(\mu,1)}$.
}\end{example}

Similarly, for each $\P$-module $M$ there exists a collection
$X$ such
that $M = (X\circ \P)/J$ for some right submodule
$J\subset X\circ \P$.
The following definition was introduced independently also
in~\cite{ginzburg-voronov}.

\begin{definition}
\label{celenka}
The module $M$ is called quadratic, if, in the above
presentation,
the collection $X$ is
concentrated in degree $1$, $X= X(1)$, and the right submodule
$J$
is generated by a subspace $G \subset (X \circ \P)(2)$.
In this case we write $M = \prezmod X{\P}G$.
\end{definition}

\begin{example}{\rm\
\label{nabijecka}
Let $X=X(1)$ be generated by one element $g$; $X(1)=
\Span(g)$. Let
$G \subset (X\circ \Ass)(2)$ be defined as $G = \Span(g(\mu)(1-
S_{21}))$, where $\mu$ and $S_{21}\in \Sigma_2$
has the same meaning as in Example~\ref{genius}.
Then it is not
difficult to see that
$\Cycl = \prezmod X{\Ass}G$ or, in a more explicit notation,
\begin{equation}
\label{Spitfire}
\Cycl = \prezmod g{\Ass}{g(\mu)- g(\mu)S_{21}}.
\end{equation}
Now we can give the characterization of $\Cycl$ traces
promised in Example~\ref{why-traces}.
Since the $\Ass$-module $\Cycl$ is generated by
$g\in \Cycl(1)$, any $\Cycl$-trace $t: \Cycl \to \E_{A,W}$ is
determined by the image $T := t(g)\in \E_{A,V}(1)= \Hom AW$. The
symmetry~(\ref{nabla}) then follows from the condition
$t(g(\mu)(1-
S_{21})) = 0$.
}\end{example}

In fact, we show that $\Cycl$ is a very special case of a
module over
the {\em cyclic\/} operad $\Ass$. Cyclic operads were
introduced by
E.~Getzler and
M.M.~Kapranov~\cite[Section~2]{getzler-kapranov:cyclic}.
The definition we are going to recall is based on the following
convention.

We interpret the symmetric group $\Sigma_{n+1}$ as the group of
permutations of the set $\{\rada 0n\}$ and $\Sigma_n$
as the subgroup of $\Sigma_{n+1}$ consisting of permutations
$\sigma \in \Sigma_{n+1}$ with $\sigma(0)= 0$. If $\tau_n
\in \Sigma_{n+1}$ is the cycle $(\rada 0n)$, then $\tau_n$
and $\Sigma_n$ generate $\Sigma_{n+1}$.

\begin{figure}[hbtp]
\begin{center}
%TexCad Options
%\grade{\off}
%\emlines{\off}
%\beziermacro{\off}
%\reduce{\on}
%\snapping{\off}
%\quality{2.00}
%\graddiff{0.01}
%\snapasp{1}
%\zoom{2.00}
\unitlength 1.00mm
\thinlines
\begin{picture}(60.84,24.84)
\put(49.67,16.00){\line(-1,-1){10.17}}
\put(39.50,5.83){\line(1,0){20.00}}
\put(59.50,5.83){\line(-1,1){10.00}}
\put(49.50,15.83){\line(0,1){7.33}}
\put(41.33,5.83){\line(0,-1){4.83}}
\put(44.83,5.83){\line(0,-1){4.83}}
\put(48.00,5.83){\line(0,-1){4.83}}
\put(57.33,5.83){\line(0,-1){4.83}}
\put(53.00,3.33){\makebox(0,0)[cc]{$\cdots$}}
\put(9.83,16.00){\line(-1,-1){10.17}}
\put(-0.33,5.83){\line(1,0){20.00}}
\put(19.67,5.83){\line(-1,1){10.00}}
\put(9.67,15.83){\line(0,1){7.33}}
\put(1.50,5.83){\line(0,-1){4.83}}
\put(5.00,5.83){\line(0,-1){4.83}}
\put(8.17,5.83){\line(0,-1){4.83}}
\put(17.50,5.83){\line(0,-1){4.83}}
\put(13.17,3.33){\makebox(0,0)[cc]{$\cdots$}}
\put(39.08,0.91){\oval(4.50,3.17)[b]}
\put(36.83,1.00){\line(0,1){22.33}}
\put(55.17,23.09){\oval(11.33,3.50)[t]}
\put(60.83,23.67){\line(0,-1){23.00}}
\put(9.67,10.33){\makebox(0,0)[cc]{$p$}}
\put(49.50,10.33){\makebox(0,0)[cc]{$p\tau_n$}}
\put(27.00,10.33){\makebox(0,0)[cc]{$\longmapsto$}}
\put(1.50,-2.67){\makebox(0,0)[cc]{$1$}}
\put(5.00,-2.67){\makebox(0,0)[cc]{$2$}}
\put(8.17,-2.67){\makebox(0,0)[cc]{$3$}}
\put(17.50,-2.67){\makebox(0,0)[cc]{$n$}}
\put(44.83,-2.67){\makebox(0,0)[cc]{$1$}}
\put(48.00,-2.67){\makebox(0,0)[cc]{$2$}}
\put(60.83,-2.67){\makebox(0,0)[cc]{$n$}}
\end{picture}
\end{center}
\caption{A `visualization' of the action of $\tau_n$. The
element
$\tau_n$ turns
$p\in \P(n)$, represented as a `thing' with $n$ inputs and one
output, a bit so that the first
input becomes the output and the output
becomes the last input of $p\cdot \tau_n$.\label{feel}}
\end{figure}
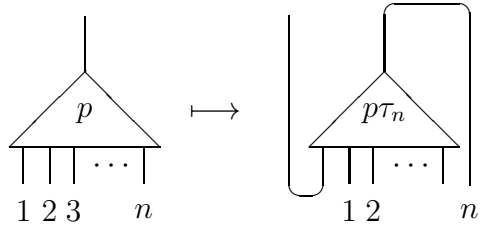

\begin{definition}
\label{tabal}
Cyclic operad is an ordinary operad $\P$ such that the right
$\Sigma_n$-action on $\P(n)$ extends, for $n\geq 1$, to
an action
of $\Sigma_{n+1}$ such that
\begin{itemize}
\item[(i)]
$\tau_1(1)= 1$, where $1\in \P(1)$ is the unit and
\item[(ii)]
for $p\in \P(m)$ and $q\in \P(n)$,
\[
\gamma(p;q,\rada 11)\cdot \tau_{m+n-1}= \znamenko{|p|\cdot|q|}
\gamma(q\cdot\tau_n;\rada11,p\cdot\tau_m).
\]
\end{itemize}
\end{definition}
An intuitive feeling for the action of $\tau_n$ is suggested by
Figure~\ref{feel}.

Let $\P = \prez ER$ be a quadratic operad. Thus $E(2)$ is a
$\Sigma_2= \Z_2$ space and the homomorphism $\sgn:
\Sigma_3\to \Z_2$ equips $E(2)$ with a $\Sigma_3$-action which
induces on $\fr(E)$ a cyclic operad structure. We say, according
to~\cite[\S3.2]{getzler-kapranov:cyclic},
that $\P$ is {\em cyclic quadratic\/},
if the subspace $R
\subset \fr(E)(3)$ is $\Sigma_4$-invariant.
In this case the operad $\P$ carries a natural cyclic structure
induced from the cyclic structure of $\fr(E)$.

An example of a cyclic quadratic operad is the operad $\Ass$,
see~\cite[Proposition~2.4]{getzler-kapranov:cyclic}
for a very explicit description of the cyclic
structure, and also for other examples of cyclic operads.

Let $\P = (\P,\gamma)$ be a cyclic operad in
the sense of Definition~\ref{tabal}. Define the $\P$-module
$M_{\P}$ as follows. As a collection, $M_{\P}(n+1)= \P(n)$,
for $n\geq 0$. The structure maps are given by
\begin{eqnarray}
\label{eqv}
\nu(x;p_0,\rada 11)&:=&\znamenko{|p_0|\cdot |x|}
\gamma(p_0\cdot\tau_{m_0};\rada 11,x), \mbox{ and}
\\
\nonumber
\nu(x;1,\rada{p_1}{p_n}) &:=& \gamma(x;\rada{p_1}{p_n}),
\end{eqnarray}
for $x\in M_{\P}(n+1)=\P(n)$, $p_i\in \P(m_i)$, $0\leq i\leq n$.
These conditions, along with the module axioms
(Definition~\ref{Amphora}), imply that
\[
\nu(x;\rada{p_0}{p_n})= \znamenko{|p_0|\cdot |x|}
\gamma(p_0\cdot\tau_{m_0};\rada11, \nu(x;1,\rada{p_1}{p_n})).
\]
The verification of module axioms of $M_{\P}$
is routine.

\begin{definition}
\label{Turmo}
We call
the module $M_\P$ the module associated to the cyclic
operad $\P$.
\end{definition}

We will also consider {\em unital\/} operads
which describe algebras {\em with unit\/}.
Unital operad is an
operad such that $\P(0)$ is nonempty, generated by an element
$\vartheta$, encoding the unit $\bk \to A$ in the corresponding
algebra
$A$. The element $\vartheta$ determines,
for $n\geq 1$ and $1\leq i\leq n$, the
`degeneracy' maps $s_i : \P(n)\to \P(n-1)$ by $s_i(p):=
\gamma(p;\rada 11,\vartheta,\rada 11)$ ($\vartheta$ at the
$i$-th
position). These maps satisfy certain commutation
relations~\cite[page~278, Proposition~3]{stasheff:TAMS63}
which follow from the axioms of an operad. For us, the most
important
is the relation $s_1(p)= s_2(p S_{21})$ for $p\in \P(2)$, which
follows from the equivariance of structure maps. This,
together with
a natural requirement that $s_1 = s_2$ on $\P(2)$, gives
that the
maps $s_1=s_2$ are $\Sigma_2$-equivariant on $\P(2)$. We
denote both
maps $s_1$ and $s_2$ (which are the same) by $s$.

\begin{example}
\label{vypocetni}
{\rm\
There is the operad $\UAss$ for associative algebras with unit,
which is the same as the operad $\Ass$ except that $\UAss(0)=
\Span(\vartheta)$. The map $s :\Ass(2) = \bk[\Sigma_2]\to \bk$
is the
standard augmentation of the group ring $\bk[\Sigma_2]$.

Another example is the operad $\UComm$ for unital commutative
algebras. It has, for $n\geq 1$, $\UComm(n)= \Comm(n)= \id$ (the
trivial
one-dimensional representation), and
$\UComm(0)= \Span(\vartheta)$. The map $s : \UComm(2)= \bk
\to \bk$
is the identity.

A more complicated example is the operad $\UPoiss$ for
Poisson algebras with unit. Here by a Poisson algebra with
unit we
mean an ordinary Poisson
algebra~\cite[Example~3.3]{markl:dl} $P = (P, \cdot, [-,-])$
with a distinguished element $1\in P$ which is a two-sided
unit for
the commutative multiplication $\cdot$, while $[x,1]= [1,x]=
0$, for
all $x\in P$. The operad $\UPoiss$ coincides with the operad
$\Poiss$
for Poisson algebras (which is very explicitly described
in~\cite{fox-markl:ContM97}),
except that $\UPoiss(0) =\Span(\vartheta)$. The
component $\UPoiss(2)$ is the direct sum $\id \oplus \sgn$,
with the trivial one-dimensional representation $\id$
concentrated in degree zero, and the signum representation
$\sgn$ in degree 1.
The map $s : \UPoiss(2)\to \bk = \UPoiss(1)$ is the
projection on the zero-dimensional component.
}\end{example}

We saw that a natural way to construct cyclic operads was to
take a
quadratic operad $\P = \prez ER$ for which the relations
$R$ were
invariant under the natural $\Sigma_4$-action; the operad
$\P$ had
then a natural cyclic structure. There is a similar approach to
unital operads.

So, let $\P = \prez ER$ be a quadratic operad and suppose we are
given an epimorphism $s :E \to \bk$. This will be a model for
the two
degeneracy maps $s_1 = s_2 : \P(2)= E \to \P(1)= \bk$.
These two maps induce degeneracy maps on the free operad
$\fr(E)$ satisfying the correct commutation relations. The
following
definition expresses the conditions assuring that this structure
preserves the relations $R$.

\begin{definition}
\label{virgo}
Let $\P = \prez ER$ be a quadratic (resp.~cyclic quadratic)
operad.
Suppose that we are given an epimorphism $s : E \to \bk$
such that
\begin{itemize}
\item[(i)]
$s$ is invariant under the $\Sigma_2$ (resp.~$\Sigma_3$, in the
cyclic case) action, and
\item[(ii)]
the induced maps $s_1,s_2,s_3: \fr(E)(3)\to \fr(E)(2)$ send the
subspace $R\subset \fr(E)(3)$ to zero.
\end{itemize}
Then the collection $\UP$, defined by $\UP(n):= \P(n)$ for
$n\geq 1$
and $\UP(0)= \bk$, has a natural structure of a unital
operad. We
call operads of this form quadratic unital (resp.~cyclic
quadratic
unital)
operads.
\end{definition}

All the three operads $\UAss$, $\UComm$ and $\UPoiss$ from
Example~\ref{vypocetni} are cyclic quadratic unital operads
in the sense of
Definition~\ref{virgo}.

\begin{proposition}
\label{zabacek}
Let $\P = \prez ER$ be a cyclic unital quadratic operad in
the sense
of Definition~\ref{virgo}. Then the associated module
$M_{\UP}$ is
quadratic,
\[
M_{\UP} =
\prezmod{\Span(\vartheta)}{\P}{{\rm Ker}(s):E \to {\bf k}}.
\]
\end{proposition}

\noindent
{\bf Proof.}
Let $X = X(1):= \Span(g)$ and consider the map $\psi :
X\circ \P \to M_{\UP}$ defined by $p(g) :=
\vartheta \in M_{\UP}(1)= \UP(0)$.
Because clearly $X\circ \P = \P$, $\psi$ is, by~(\ref{eqv}),
given as
\[
(X\circ \P)(n)= \P(n)\ni p \longmapsto \nu(\vartheta;p) =
\gamma(p\cdot\tau_n; \rada 11,\vartheta)\in M_{\UP}(n)
= \UP(n-1).
\]
Thus $\psi$ is a sequence
$\psi(n):\P(n)\mapsto \UP(n-1)$ given by
$\psi(n)(p):= \gamma(p\cdot\tau_n;
\rada 11,\vartheta)$. These maps are epimorphisms, because $E$
generates $\P$ and $s:E \to \bk$ (= the composition with
$\vartheta$)
is epi, by assumption.

On the other hand, $\psi(2):\P(2)=E \to \UP(1)=\bk$ is exactly
the map
$s$, thus the submodule of $X\circ \P$ generated by $\Ker(s)$ is
certainly contained in the kernel of $\psi$. A moment's
reflections
shows that the whole kernel of $\psi$ is generated by
$\Ker(s)$.\qed

\begin{example}{\rm\
\label{kacirek}
If $\P = \UAss$ is the operad for unital associative algebras
from
Example~\ref{vypocetni},
then ${\rm Ker}(s:E= \bk[\Sigma_2] \to {\bf k})=
\sgn$. So, by Proposition~\ref{zabacek}, $M_{\UAss}=
\prezmod{\vartheta}{\Ass}{\sgn}$, thus $M_{\UAss} = \Cycl$,
by~(\ref{Spitfire}).
For the operad $\UComm$ the kernel of $s$ is trivial,
hence
\[
M_{\UComm} = \Span(\vartheta)\circ \Comm \cong \Comm.
\]
In other words, $M_{\UComm}$ is the operad
$\Comm$ considered as a right
module over itself. We recommend to the reader to make the
similar
discussion for the operad $\UPoiss$.
}\end{example}

\section{Cyclohedron as the cobar construction}
\label{penezenka}

For a graded vector space $V$, let
$\dual V$ denote the graded dual of $V$, i.e.\ the graded
vector space $\bigoplus_p(\dual V)_p$ with $(\dual V)_p
={\rm Hom}^{-p}(V,{\bk})={\rm Hom}(V_{p},{\bk})$, the space
of linear maps from $V_{p}$ to $\bk$. If $V$ is a right
$\Sigma_n$-module, then we equip $\dual V$ with the transposed
$\Sigma_n$-action.
We believe that there is no real risk of confusion of the ${}^*$
indicating the dual with the star indicating the degree.

Recall
also~\cite[\S3.5.1]{ginzburg-kapranov:DMJ94}
that a {\em cooperad\/} is a collection $\Q=\coll {\Q}$
together with a system of maps
\[
\omega = \omega_{\rada{m_1}{m_l}}: \Q(m_1+\cdots+m_l)
\to \Q(l)\otimes \Q(m_1)\otimes \cdots \otimes \Q(m_l),
\]
which satisfy the axioms which are exactly the duals of
the axioms
for an operad. A typical example of a cooperad is the dual
$\dual\P$ of
an operad $\P$, i.e.~the collection $\dual\P =
\{\dualI{\P(n)}\}_{n\geq 1}$
with the cooperad structure defined by the dualization of the
structure maps of $\P$; here an obvious finite type assumption
is
necessary, but it will be always satisfied in the paper and
we will
make no explicit comments about it.

Observe that if $\P$ is an operad, then both the suspension
$\ss \P$ and the desuspension $\ss^{-1} \P$ introduced
in Section~\ref{hrnicek1}, with the sign convention
of~(\ref{sign-factor}),
have a natural operad
structure induced from the operad structure on $\P$.

The structure maps of a cooperad $\Q$ determine (and are
determined
by) a map ${\overline \nu}:
\Q\to \fr(\Q)$ of collections. Composing this
map with the (de)suspensions gives a degree -1
map $\desusp \Q \to \fr(\desusp \Q)$
which uniquely (because of the freeness of $\fr(\desusp \Q)$)
extends
to
a degree -1 derivation $\dcob$ of the operad $\fr(\desusp
\Q)$ which
satisfies, as a consequence of the axioms of a cooperad,
$\dcob \circ\dcob =0$.
The differential graded operad
\[
\Cob(\Q) := (\fr(\desusp \Q),\dcob)
\]
is called the {\em cobar construction\/} on the cooperad
$\Q$~\cite[\S3.2.12]{ginzburg-kapranov:DMJ94}.

Let $\P = \prez ER$ be a quadratic operad. Take the dual $\dual
E$ of
$E$ and let $R^{\perp}$ be the annihilator of the space
$R\subset \fr(E)(3)$ in
$\fr(\dual E \otimes \sgn)(3)$. The quadratic operad $\P^! :=
\prez {\dual E \otimes \sgn}{R^{\perp}}$ is,
according to~\cite[\S2.1.9]{ginzburg-kapranov:DMJ94}, called the
{\em Koszul\/} (or {\em quadratic\/}) {\em dual\/}
of the operad $\P$. We always have a map
$\desusp \dualI{\ss \P} \to \P^!$ of collections defined as
the composition
\[
\desusp \dualI{\ss \P}\stackrel{\rm proj.}{\longrightarrow}
(\desusp \dualI{\ss \P)(2)} = \dual E\otimes \sgn =
\P^!(2) \hookrightarrow \P^!
\]
which extends, by the freeness of
$\fr(\desusp \dualI{\ss \P})$, to a differential graded
operad map
\begin{equation}
\label{budicek}
\pi : \Cob(\dualI{\ss \P}) = (\fr(\desusp \dualI{\ss \P}),\dcob)
\to (\P^!,\partial=0).
\end{equation}
The quadratic operad $\P$ is called
{\em Koszul\/}~\cite[Definition~4.1.3]{ginzburg-kapranov:DMJ94}
if the map in~(\ref{budicek}) is a homology isomorphism.
In the rest of the paper, the space of generators $E=E(2)$ of a
quadratic
operad will always be ungraded, concentrated in degree
zero. Then the
components of
both the operad $\P$ and its dual $\P^!$ are concentrated
in degree
zero as well, and the Koszulness
implies that the complex $(\Cob(\dualI{\ss \P})(n),
\dcob)$ is acyclic in positive dimensions, for all $n$. On
the other
hand, as shown
in~\cite[Theorem~4.1.13]{ginzburg-kapranov:DMJ94},
this acyclicity condition implies the
Koszulness of $\P$.

\begin{example}{\rm\
\label{yhr}
The operad $\Ass$ is well-known to be Koszul self-dual, $\Ass =
\Ass^!$~\cite[Theorem~2.1.11]{ginzburg-kapranov:DMJ94}, and
Koszul~\cite[Corrolary~4.2.7]{ginzburg-kapranov:DMJ94}.
We have already remarked that the
operad $\barA$ of cellular chains on the (symmetrized)
associahedron
$\barK$ coincides, as a differential graded operad, with
the cobar
construction $\Cob(\dualI{\ss\Ass})$. Let us give an explicit
illustration of this statement.

The $n$-th piece $(\desusp \dualI{\ss \Ass})(n)$ of the
collection
$(\desusp \dualI{\ss \Ass})$ is isomorphic to
one copy of the regular
representation $\bk[\Sigma_n]$ concentrated in degree $(n-2)$.
The isomorphism is not unique, but we may choose, for example,
$\lambda : \bk[\Sigma_n] \to \desusp (\ss \Ass)^*(n)$ given by
$\lambda(\sigma)(\susp^{n-2}\rho) := \sgn(\sigma)\cdot
\delta_{\sigma,\rho}$, where $\sigma, \rho \in \Sigma_n$ and the
meaning of the `Kronecker delta' $\delta_{\sigma,\rho}$
is clear.

Formula~(\ref{fax})
then describes $\fr(\desusp \dualI{\ss \Ass})(n)$ as the
vector space
spanned by the set of all {\em planar\/} (rooted, labeled)
$n$-trees,
having at least binary vertices. The identification of these
trees
with the bracketings of $\rada{\sigma(1)}{\sigma(n)}$,
$\sigma \in \Sigma_n$, i.e.~with the elements of the set
$\barB(n)$,
is classical -- see~\cite[\S1.4]{boardman-vogt:73};
two examples are
shown on Figure~\ref{hreben}.
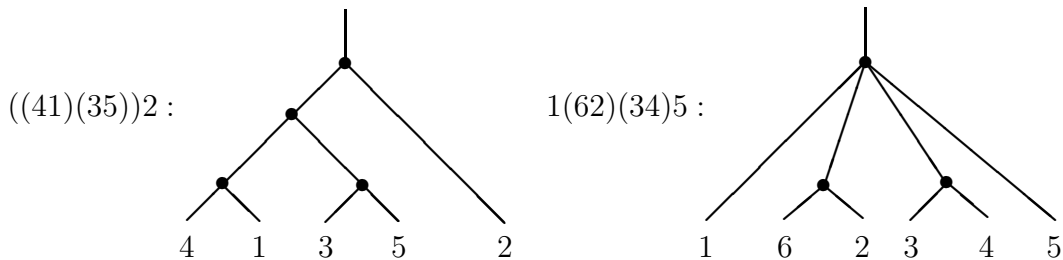
\begin{figure}[hbtp]
\begin{center}
%TexCad Options
%\grade{\off}
%\emlines{\off}
%\beziermacro{\off}
%\reduce{\off}
%\snapping{\off}
%\quality{2.00}
%\graddiff{1.00}
%\snapasp{1}
%\zoom{1.91}
\unitlength 0.70mm
\thicklines
\begin{picture}(174.84,50.34)
\put(40.00,50.00){\line(0,-1){10.00}}
\put(40.00,40.00){\line(-1,-1){30.00}}
\put(40.00,39.83){\line(1,-1){30.33}}
\put(30.00,30.00){\line(1,-1){20.17}}
\put(29.83,30.00){\makebox(0,0)[cc]{$\bullet$}}
\put(40.00,39.67){\makebox(0,0)[cc]{$\bullet$}}
\put(16.74,16.91){\makebox(0,0)[cc]{$\bullet$}}
\put(43.26,16.57){\makebox(0,0)[cc]{$\bullet$}}
\put(10.00,5.00){\makebox(0,0)[cc]{$4$}}
\put(23.67,5.00){\makebox(0,0)[cc]{$1$}}
\put(36.33,4.83){\makebox(0,0)[cc]{$3$}}
\put(50.17,5.00){\makebox(0,0)[cc]{$5$}}
\put(70.33,5.00){\makebox(0,0)[cc]{$2$}}
\put(7.67,31.33){\makebox(0,0)[rc]{$((41)(35))2:$}}
\put(108.67,10.00){\line(1,1){30.17}}
\put(138.84,40.17){\line(0,1){10.17}}
\put(138.84,40.17){\line(-1,-3){7.89}}
\put(131.00,16.67){\line(-6,-5){7.67}}
\put(130.84,16.67){\line(6,-5){7.50}}
\put(138.84,40.00){\line(2,-3){15.56}}
\put(154.00,16.67){\line(-1,-1){6.67}}
\put(154.17,17.17){\line(6,-5){7.83}}
\put(138.84,40.00){\line(6,-5){36.00}}
\put(138.84,39.83){\makebox(0,0)[cc]{$\bullet$}}
\put(130.84,16.50){\makebox(0,0)[cc]{$\bullet$}}
\put(154.17,17.00){\makebox(0,0)[cc]{$\bullet$}}
\put(108.50,5.00){\makebox(0,0)[cc]{$1$}}
\put(123.50,5.00){\makebox(0,0)[cc]{$6$}}
\put(138.17,5.00){\makebox(0,0)[cc]{$2$}}
\put(147.50,4.83){\makebox(0,0)[cc]{$3$}}
\put(161.84,4.83){\makebox(0,0)[cc]{$4$}}
\put(174.67,4.83){\makebox(0,0)[cc]{$5$}}
\put(108.33,31.34){\makebox(0,0)[rc]{$1(62)(34)5:$}}
\put(43.24,16.56){\line(-1,-1){6.97}}
\put(16.74,16.74){\line(1,-1){6.97}}
\end{picture}
\end{center}
\caption{Two examples of an identification of planar trees with
elements of $\B$.\label{hreben}}
\end{figure}
The fact that the cellular differential coincides with
$\dcob$ is a routine combinatorics.
See~\cite[Example~4.1]{markl:zebrulka} for details.

The identification above shows that the operad $\Ass$
is Koszul. More precisely, we know that $H_*(\barK)= H_0(\barK)=
\Ass$, because $\barK(n)$ is the union of convex polyhedra
indexed by
the elements of the symmetric group. On the other hand, due
to the
identification above, $H_*(\Cob(\dualI{\ss \P}))= H_*(\barK)$,
thus the map
$\pi$ of~(\ref{budicek}) is a homology isomorphism.
}\end{example}

A {\em comodule\/} over a cooperad $\Q$
is a collection $N = \coll N$ together with structure maps
\[
\kappa = \kappa_{\rada{m_1}{m_l}}: N(m_1+\cdots+m_l)
\to N(l)\otimes \Q(m_1)\otimes \cdots \otimes \Q(m_l),
\]
satisfying axioms dual to the axioms of a module over an operad.
An example is the dual $\dual \M$ of a $\P$-module $\M$,
which is a
natural comodule over the cooperad $\dual\P$.

Let $N$ be a $\Q$-comodule. As in the case of
cooperads, the structure maps induce a degree -1 differential
$\dcob$
on the free module $N\circ \fr(\desusp \Q)$. The right
differential graded $\Cob(\Q)$-module
\[
\Cob(N;\Q) :=(N\circ \fr(\desusp \Q),\dcob)
\]
is called the (right) {\em cobar construction\/} on the
$\Q$-comodule
$N$.

Suppose that $M = \prezmod X{\P}G$ is a
quadratic $\P$-module, in the
sense of Definition~\ref{celenka}, over a quadratic operad $\P =
\prez ER$. Take
the dual $\dual X$ and let $G^{\perp}$ be the annihilator of
$G$ in $(\dual
X\circ (\dual E \otimes \sgn))(2)$. Then the quadratic
$\P^!$-module
$M^! := \prezmod {\dual X}{\P^!}{G^{\perp}}$
is called the {\em Koszul\/} (or {\em quadratic\/})
{\em dual\/} of $M$. The above definitions were
independently made in~\cite{ginzburg-voronov}.

\begin{proposition}
The $\Ass$-module $\Cycl$ is Koszul self-dual, $\Cycl^! =
\Cycl$.
\end{proposition}

\noindent
{\bf Proof.}
Under the notation of Example~\ref{nabijecka}, the
$\Sigma_2$-space
$(X\circ \Ass)(2)$ is the direct sum $\id \oplus \sgn$ of
the trivial
and signum representations, while clearly $R= \sgn$,
generated by
$g(\mu) - g(\mu)S_{21}$. Then $R^\perp$ is easily seen to
be $\sgn$
and the proposition follows.\qed

Observe that, for a $\P$-module $M = \coll M$, the collection
$\ss M$ has an induced
$\ss \P$-module structure.
For a quadratic $\P$-module $\prezmod X{\P}G$ over a quadratic
operad
$\P =
\prez ER$ we have, as in~(\ref{budicek}), the map
\begin{equation}
\label{Opusem}
\pi :\Cob(\dualI {\ss M},\dualI {\ss \P}) \to (M^!,\partial
= 0),
\end{equation}
induced by the composition
\[
\dualI {\ss M} \stackrel{\rm proj.}{\longrightarrow}
\dualI{\ss M}(1) = \dual X = M^!(1) \hookrightarrow M^!.
\]
The following definition was independently made
in~\cite{ginzburg-voronov}.

\begin{definition}
\label{telefon}
A quadratic module $M$ over a quadratic operad $\P$
is called Koszul if the
map $\pi$ in~(\ref{Opusem}) is a homology isomorphism.
\end{definition}

\begin{theorem}
\label{Amphora1}
The module $\Cycl$ is Koszul.
\end{theorem}

The theorem will follow from Theorem~\ref{ucpavka} which
explicitly identifies the cobar
construction $\Cob(\dualI{\ss\Cycl}, \dualI{\ss\Ass})$
to the cellular chain
complex
$\M$ of the symmetrized cyclohedron $\barW$. Thus the
Koszulness of
$\Cycl$ is, as in the case of the operad $\Ass$,
a consequence of the fact that the cyclohedron is a convex
polyhedron. At the end of the
section we formulate a more general statement which also implies
Theorem~\ref{Amphora1}.

\begin{theorem}
\label{ucpavka}
The cobar construction $\Cob(\dualI{\ss\Cycl}, \dualI{\ss
\Ass})$ is
isomorphic
to the cellular chain complex $\M = CC_*(\barW)$
of the cyclohedron $\barW$.
\end{theorem}

\noindent
{\bf Proof.}
As we observed in Example~\ref{yhr}, the collection
$\desusp\dualI{\ss \Ass}(n)$ consists of one copy of the regular
representation $\bk[\Sigma_n]$ concentrated in degree
$(n-2)$. Similarly,
the collection $\dualI{\ss \Cycl}$ is
isomorphic to the collection $\ss \Cycl$.
We are going to give an explicit description of the space
$(\dualI{\ss\Cycl}\
\circ \desusp\dualI{\ss \Ass})(n)$, similar to that of
Lemma~\ref{MGD}.
Suppose that $T_i$ is,
for each $1\leq i \leq l$, a planar $m_i$-tree,
whose all vertices are at least binary. Let $R(\rada{T_1}{T_l})$
be
the tree obtained by grafting the trees $\rada{T_1}{T_l}$ at the
inputs of the `$l$-rake,' or, pictorially:
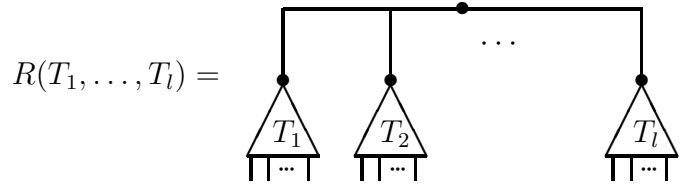
\begin{center}
\def\SetFigFont#1#2#3{{}}
\def\smash#1{{#1}}
\setlength{\unitlength}{0.0062500in}%
\begin{picture}(370,213)(20,440)
\thicklines
\put(150,540){\line( 0, 1){63}}
\put(150,603){\line( 1, 0){300}}
\put(450,603){\line( 0,-1){63}}
\put(150,540){\line(-1,-2){ 30}}
\put(120,480){\line( 1, 0){ 60}}
\put(180,480){\line(-1, 2){ 30}}
\put(123,480){\line( 0,-1){ 21}}
\put(138,480){\line( 0,-1){ 21}}
\put(171,480){\line( 0,-1){ 21}}
\put(240,603){\line( 0,-1){63}}
\put(240,540){\line(-1,-2){ 30}}
\put(210,480){\line( 1, 0){ 60}}
\put(270,480){\line(-1, 2){ 30}}
\put(450,540){\line(-1,-2){ 30}}
\put(420,480){\line( 1, 0){ 60}}
\put(480,480){\line(-1, 2){ 30}}
\put(261,480){\line( 0,-1){ 21}}
\put(471,480){\line( 0,-1){ 21}}
\put(426,480){\line( 0,-1){ 21}}
\put(441,480){\line( 0,-1){ 21}}
\put(216,480){\line( 0,-1){ 21}}
\put(231,480){\line( 0,-1){ 21}}
\put(153,462){\makebox(0,0)[b]{
             \smash{
                    \SetFigFont{12}{14.4}{rm}$\cdot\!\!\cdot\!\!\cdot$
                               }
                               }}
\put(330,568){\makebox(0,0)[b]{\smash{\SetFigFont{12}{14.4}{rm}$\cdots$}}}
\put(246,462){\makebox(0,0)[b]%
{\smash{\SetFigFont{12}{14.4}{rm}$\cdot\!\!\cdot\!\!\cdot$}}}
\put(456,462){\makebox(0,0)[b]%
{\smash{\SetFigFont{12}{14.4}{rm}$\cdot\!\!\cdot\!\!\cdot$}}}

\put(300,597){\makebox(0,0)[b]{\smash{\SetFigFont{12}{14.4}{rm}$\bullet$}}}

\put(150,537){\makebox(0,0)[b]{\smash{\SetFigFont{12}{14.4}{rm}$\bullet$}}}
\put(154,492){\makebox(0,0)[b]{\smash{\SetFigFont{12}{14.4}{rm}$T_1$}}}

\put(240,537){\makebox(0,0)[b]{\smash{\SetFigFont{12}{14.4}{rm}$\bullet$}}}
\put(244,492){\makebox(0,0)[b]{\smash{\SetFigFont{12}{14.4}{rm}$T_2$}}}

\put(450,537){\makebox(0,0)[b]{\smash{\SetFigFont{12}{14.4}{rm}$\bullet$}}}
\put(454,492){\makebox(0,0)[b]{\smash{\SetFigFont{12}{14.4}{rm}$T_l$}}}

\put(10,540  ){\makebox(0,0)[b]%
{\smash{\SetFigFont{12}{14.4}{rm}$R(\rada{T_1}{T_l}) =$}}}

\end{picture}
\end{center}
Let $J(n)$ be the free graded right $\Sigma_n$-module generated
by the
symbols $R(\rada{T_1}{T_l})$, $m_1+\cdots+ m_l = n$, with
the degree
defined by
\[
\deg(R(\rada{T_1}{T_l}))= n-1-\sum_{i=1}^l \#{\rm vert}(T_i),
\]
where $ \#{\rm vert}(T_i)$ is the number of vertices of
the tree $T_i$. For $\zeta \in
\ZZ_l$ we put (see~(\ref{pejska_Mikinka}) for the notation)
\begin{equation}
\label{monce1}
\zeta(R(\rada{T_1}{T_l}) \times \sigma) :=\sgn(\zeta) \cdot
R(\rada{T_{\eta^{-1}(1)}}{T_{\eta^{-1}(l)}}) \times
\zeta(\rada{m_1}{m_l}) \sigma.
\end{equation}

Then the graded vector space $(\dualI{\ss\Cycl}\
\circ \desusp\dualI{\ss \Ass})(n)$ is spanned by equivalence
classes of
elements
\[
R(\rada{T_1}{T_l}) \times \sigma \in J(n),
\]
modulo the left action of the group $\ZZ_l$
introduced in~(\ref{monce1}). The right
action of the group $\Sigma_n$ is given by
\[
[R(\rada{T_1}{T_l}) \times \sigma]\cdot \rho :=
[R(\rada{T_1}{T_l}) \times \sigma\rho].
\]

\begin{figure}[hbtp]
\begin{center}
%TexCad Options
%\grade{\off}
%\emlines{\off}
%\beziermacro{\off}
%\reduce{\on}
%\snapping{\off}
%\quality{2.00}
%\graddiff{0.01}
%\snapasp{1}
%\zoom{1.00}
\unitlength 1.20mm
\thicklines
\begin{picture}(91.08,167.33)
\put(-11.20,167.32){\makebox(0,0)[rc]{\fbox{$n=1$:}}}
\put(-1.37,167.15){\makebox(0,0)[cc]{$\bullet$}}
\put(-1.37,159.65){\makebox(0,0)[cc]{$1$}}
\put(-11.20,150.48){\makebox(0,0)[rc]{\fbox{$n=2$:}}}
\put(12.30,150.48){\line(1,0){5.00}}
\put(17.30,150.48){\line(0,-1){5.00}}
\put(12.13,150.48){\line(0,-1){5.00}}
\put(14.72,150.48){\makebox(0,0)[cc]{$\bullet$}}
\put(-3.87,150.48){\line(1,0){5.00}}
\put(1.13,150.48){\line(0,-1){5.00}}
\put(-4.04,150.48){\line(0,-1){5.00}}
\put(-1.45,150.48){\makebox(0,0)[cc]{$\bullet$}}
\put(-4.04,142.65){\makebox(0,0)[cc]{$1$}}
\put(1.13,142.65){\makebox(0,0)[cc]{$2$}}
\put(12.13,142.65){\makebox(0,0)[cc]{$2$}}
\put(17.30,142.65){\makebox(0,0)[cc]{$1$}}
\put(7.97,150.48){\makebox(0,0)[cc]{$=-$}}
\put(20.79,150.48){\makebox(0,0)[lc]{$=12,$}}
\put(42.46,151.32){\makebox(0,0)[cc]{$\bullet$}}
\put(42.46,147.65){\makebox(0,0)[cc]{$\bullet$}}
\put(45.84,149.48){\makebox(0,0)[lc]{$= (12),$}}
\put(61.76,151.32){\makebox(0,0)[cc]{$\bullet$}}
\put(61.63,147.65){\makebox(0,0)[cc]{$\bullet$}}
\put(65.01,149.48){\makebox(0,0)[lc]{$= (21)$}}
\put(40.13,141.30){\makebox(0,0)[cc]{$1$}}
\put(44.96,141.30){\makebox(0,0)[cc]{$2$}}
\put(59.30,141.30){\makebox(0,0)[cc]{$2$}}
\put(64.30,141.30){\makebox(0,0)[cc]{$1$}}
\put(-11.20,129.52){\makebox(0,0)[rc]{\fbox{$n=3$:}}}
\put(-6.37,129.69){\line(1,0){10.00}}
\put(3.63,129.69){\line(0,-1){5.17}}
\put(-6.37,129.69){\line(0,-1){5.17}}
\put(-6.37,124.52){\line(0,0){0.00}}
\put(-1.37,129.69){\makebox(0,0)[cc]{$\bullet$}}
\put(11.13,129.69){\line(1,0){10.00}}
\put(28.63,129.69){\line(1,0){10.00}}
\put(21.13,129.69){\line(0,-1){5.17}}
\put(38.63,129.69){\line(0,-1){5.17}}
\put(11.13,129.69){\line(0,-1){5.17}}
\put(28.63,129.69){\line(0,-1){5.17}}
\put(11.13,124.52){\line(0,0){0.00}}
\put(28.63,124.52){\line(0,0){0.00}}
\put(16.13,129.69){\makebox(0,0)[cc]{$\bullet$}}
\put(33.63,129.69){\makebox(0,0)[cc]{$\bullet$}}
\put(7.46,130.35){\makebox(0,0)[cc]{$=$}}
\put(25.13,130.52){\makebox(0,0)[cc]{$=$}}
\put(42.80,130.52){\makebox(0,0)[cc]{$=$}}
\put(48.30,130.52){\makebox(0,0)[cc]{$123$}}
\put(-6.37,122.02){\makebox(0,0)[cc]{$1$}}
\put(-1.37,122.02){\makebox(0,0)[cc]{$2$}}
\put(3.63,122.02){\makebox(0,0)[cc]{$3$}}
\put(11.13,122.02){\makebox(0,0)[cc]{$2$}}
\put(16.13,122.02){\makebox(0,0)[cc]{$3$}}
\put(21.13,122.02){\makebox(0,0)[cc]{$1$}}
\put(28.63,122.02){\makebox(0,0)[cc]{$3$}}
\put(33.63,122.02){\makebox(0,0)[cc]{$1$}}
\put(38.63,122.02){\makebox(0,0)[cc]{$2$}}
\put(-1.37,61.85){\line(-1,-1){6.17}}
\put(-1.37,61.85){\line(1,-1){6.00}}
\put(-1.37,61.69){\makebox(0,0)[cc]{$\bullet$}}
\put(-1.37,65.19){\makebox(0,0)[cc]{$\bullet$}}
\put(25.96,61.85){\line(-1,-1){6.17}}
\put(53.30,61.85){\line(-1,-1){6.17}}
\put(25.96,61.85){\line(1,-1){6.00}}
\put(53.30,61.85){\line(1,-1){6.00}}
\put(25.96,61.69){\makebox(0,0)[cc]{$\bullet$}}
\put(53.30,61.69){\makebox(0,0)[cc]{$\bullet$}}
\put(25.96,65.19){\makebox(0,0)[cc]{$\bullet$}}
\put(53.30,65.19){\makebox(0,0)[cc]{$\bullet$}}
\put(-7.54,52.02){\makebox(0,0)[cc]{$1$}}
\put(-1.37,52.02){\makebox(0,0)[cc]{$2$}}
\put(4.63,52.02){\makebox(0,0)[cc]{$3$}}
\put(8.96,62.19){\makebox(0,0)[cc]{$=(123),$}}
\put(19.80,52.19){\makebox(0,0)[cc]{$3$}}
\put(25.96,52.19){\makebox(0,0)[cc]{$1$}}
\put(31.96,52.19){\makebox(0,0)[cc]{$2$}}
\put(36.96,62.35){\makebox(0,0)[cc]{$=(312)$,}}
\put(47.13,52.19){\makebox(0,0)[cc]{$1$}}
\put(53.30,52.19){\makebox(0,0)[cc]{$2$}}
\put(59.30,52.19){\makebox(0,0)[cc]{$3$}}
\put(64.13,62.19){\makebox(0,0)[cc]{$=(123)$}}
\put(-1.61,36.99){\line(-1,-1){6.31}}
\put(-1.61,41.25){\makebox(0,0)[cc]{$\bullet$}}
\put(-1.61,36.61){\makebox(0,0)[cc]{$\bullet$}}
\put(26.22,36.99){\line(-1,-1){6.31}}
\put(54.04,36.99){\line(-1,-1){6.31}}
\put(26.22,41.25){\makebox(0,0)[cc]{$\bullet$}}
\put(54.04,41.25){\makebox(0,0)[cc]{$\bullet$}}
\put(26.22,36.61){\makebox(0,0)[cc]{$\bullet$}}
\put(54.04,36.61){\makebox(0,0)[cc]{$\bullet$}}
\put(-8.10,27.90){\makebox(0,0)[cc]{$1$}}
\put(5.81,27.90){\makebox(0,0)[cc]{$3$}}
\put(-6.37,106.53){\line(1,0){8.67}}
\put(2.30,106.53){\line(0,-1){4.33}}
\put(2.30,102.19){\line(-1,-1){4.00}}
\put(2.30,102.19){\line(1,-1){4.00}}
\put(-6.37,106.53){\line(0,-1){8.00}}
\put(2.30,101.86){\makebox(0,0)[cc]{$\bullet$}}
\put(27.96,106.53){\line(-1,0){8.67}}
\put(19.29,106.53){\line(0,-1){4.33}}
\put(19.29,102.19){\line(1,-1){4.00}}
\put(19.29,102.19){\line(-1,-1){4.00}}
\put(27.96,106.53){\line(0,-1){8.00}}
\put(19.29,101.86){\makebox(0,0)[cc]{$\bullet$}}
\put(11.96,106.53){\makebox(0,0)[cc]{$=-$}}
\put(52.63,106.53){\line(1,0){8.67}}
\put(61.30,106.53){\line(0,-1){4.33}}
\put(61.30,102.19){\line(-1,-1){4.00}}
\put(61.30,102.19){\line(1,-1){4.00}}
\put(52.63,106.53){\line(0,-1){8.00}}
\put(61.30,101.86){\makebox(0,0)[cc]{$\bullet$}}
\put(86.96,106.53){\line(-1,0){8.67}}
\put(78.29,106.53){\line(0,-1){4.33}}
\put(78.29,102.19){\line(1,-1){4.00}}
\put(78.29,102.19){\line(-1,-1){4.00}}
\put(86.96,106.53){\line(0,-1){8.00}}
\put(78.29,101.86){\makebox(0,0)[cc]{$\bullet$}}
\put(70.96,106.53){\makebox(0,0)[cc]{$=-$}}
\put(-6.37,85.19){\line(1,0){8.67}}
\put(2.30,85.19){\line(0,-1){4.33}}
\put(2.30,80.86){\line(-1,-1){4.00}}
\put(2.30,80.86){\line(1,-1){4.00}}
\put(-6.37,85.19){\line(0,-1){8.00}}
\put(2.30,80.52){\makebox(0,0)[cc]{$\bullet$}}
\put(27.96,85.19){\line(-1,0){8.67}}
\put(19.29,85.19){\line(0,-1){4.33}}
\put(19.29,80.86){\line(1,-1){4.00}}
\put(19.29,80.86){\line(-1,-1){4.00}}
\put(27.96,85.19){\line(0,-1){8.00}}
\put(19.29,80.52){\makebox(0,0)[cc]{$\bullet$}}
\put(11.96,85.19){\makebox(0,0)[cc]{$=-$}}
\put(42.40,147.85){\line(-1,-2){2.32}}
\put(42.40,147.85){\line(1,-2){2.32}}
\put(61.67,147.85){\line(-2,-5){1.91}}
\put(61.67,147.70){\line(3,-5){2.78}}
\put(64.45,143.07){\line(0,0){0.06}}
\put(-1.59,36.88){\line(1,-1){6.50}}
\put(26.24,36.88){\line(1,-1){6.50}}
\put(54.08,36.88){\line(1,-1){6.50}}
\put(0.58,34.55){\line(-1,-1){3.67}}
\put(0.58,34.55){\makebox(0,0)[cc]{$\bullet$}}
\put(28.41,34.55){\line(-1,-1){3.67}}
\put(56.24,34.55){\line(-1,-1){3.67}}
\put(28.41,34.55){\makebox(0,0)[cc]{$\bullet$}}
\put(56.24,34.55){\makebox(0,0)[cc]{$\bullet$}}
\put(-1.46,17.82){\line(1,-1){6.31}}
\put(-1.46,22.08){\makebox(0,0)[cc]{$\bullet$}}
\put(-1.46,17.44){\makebox(0,0)[cc]{$\bullet$}}
\put(-1.48,17.71){\line(-1,-1){6.50}}
\put(-3.65,15.38){\line(1,-1){3.67}}
\put(-3.65,15.38){\makebox(0,0)[cc]{$\bullet$}}
\put(26.37,17.82){\line(1,-1){6.31}}
\put(54.20,17.82){\line(1,-1){6.31}}
\put(26.37,22.08){\makebox(0,0)[cc]{$\bullet$}}
\put(54.20,22.08){\makebox(0,0)[cc]{$\bullet$}}
\put(26.37,17.44){\makebox(0,0)[cc]{$\bullet$}}
\put(54.20,17.44){\makebox(0,0)[cc]{$\bullet$}}
\put(26.35,17.71){\line(-1,-1){6.50}}
\put(54.18,17.71){\line(-1,-1){6.50}}
\put(24.18,15.38){\line(1,-1){3.67}}
\put(52.01,15.38){\line(1,-1){3.67}}
\put(24.18,15.38){\makebox(0,0)[cc]{$\bullet$}}
\put(52.01,15.38){\makebox(0,0)[cc]{$\bullet$}}
\put(-1.37,167.33){\line(0,-1){5.33}}
\put(42.50,147.73){\line(0,1){3.73}}
\put(61.70,147.73){\line(0,1){3.73}}
\put(-1.37,129.73){\line(0,-1){5.33}}
\put(16.10,129.73){\line(0,-1){5.33}}
\put(33.56,129.73){\line(0,-1){5.33}}
\put(-1.90,106.53){\makebox(0,0)[cc]{$\bullet$}}
\put(23.63,106.53){\makebox(0,0)[cc]{$\bullet$}}
\put(57.03,106.53){\makebox(0,0)[cc]{$\bullet$}}
\put(82.69,106.53){\makebox(0,0)[cc]{$\bullet$}}
\put(-1.90,85.20){\makebox(0,0)[cc]{$\bullet$}}
\put(23.63,85.20){\makebox(0,0)[cc]{$\bullet$}}
\put(-1.37,65.33){\line(0,-1){3.60}}
\put(25.83,62.00){\line(0,1){3.20}}
\put(53.30,61.87){\line(0,1){3.33}}
\put(-1.37,61.73){\line(0,-1){6.13}}
\put(25.83,61.73){\line(0,-1){6.13}}
\put(53.30,61.73){\line(0,-1){6.27}}
\put(-1.64,41.33){\line(0,-1){4.40}}
\put(-1.50,22.13){\line(0,-1){4.53}}
\put(26.23,41.33){\line(0,-1){4.40}}
\put(54.10,41.33){\line(0,-1){4.40}}
\put(26.36,22.13){\line(0,-1){4.53}}
\put(54.23,22.13){\line(0,-1){4.53}}
\put(52.66,95.95){\makebox(0,0)[cc]{2}}
\put(57.19,95.95){\makebox(0,0)[cc]{3}}
\put(65.21,95.95){\makebox(0,0)[cc]{1}}
\put(74.34,95.95){\makebox(0,0)[cc]{3}}
\put(82.36,95.95){\makebox(0,0)[cc]{1}}
\put(86.89,95.95){\makebox(0,0)[cc]{2}}
\put(91.08,106.58){\makebox(0,0)[lc]{= 2(31)}}
\put(-6.44,96.12){\makebox(0,0)[cc]{1}}
\put(-1.73,96.12){\makebox(0,0)[cc]{2}}
\put(6.29,96.12){\makebox(0,0)[cc]{3}}
\put(15.24,96.12){\makebox(0,0)[cc]{2}}
\put(23.26,96.12){\makebox(0,0)[cc]{3}}
\put(27.97,96.12){\makebox(0,0)[cc]{1}}
\put(-6.44,75.03){\makebox(0,0)[cc]{3}}
\put(-1.73,75.03){\makebox(0,0)[cc]{1}}
\put(6.29,75.03){\makebox(0,0)[cc]{2}}
\put(15.24,75.03){\makebox(0,0)[cc]{1}}
\put(23.44,75.03){\makebox(0,0)[cc]{2}}
\put(27.97,75.03){\makebox(0,0)[cc]{3}}
\put(30.93,106.76){\makebox(0,0)[lc]{$=1(23)$,}}
\put(30.93,85.14){\makebox(0,0)[lc]{$=3(12)$}}
\put(-3.13,28.13){\makebox(0,0)[cc]{$2$}}
\put(19.89,28.13){\makebox(0,0)[cc]{$2$}}
\put(24.77,28.13){\makebox(0,0)[cc]{$3$}}
\put(32.79,28.13){\makebox(0,0)[cc]{$1$}}
\put(47.78,28.31){\makebox(0,0)[cc]{$3$}}
\put(52.49,28.31){\makebox(0,0)[cc]{$1$}}
\put(60.51,28.31){\makebox(0,0)[cc]{$2$}}
\put(3.67,38.24){\makebox(0,0)[lc]{$=(1(23))$,}}
\put(30.87,38.42){\makebox(0,0)[lc]{$=(2(31))$,}}
\put(58.94,38.59){\makebox(0,0)[lc]{$=(3(12))$}}
\put(3.72,19.67){\makebox(0,0)[lc]{$=(1(23))$,}}
\put(30.88,19.83){\makebox(0,0)[lc]{$=(2(31))$,}}
\put(58.88,20.00){\makebox(0,0)[lc]{$=((31)2)$}}
\put(-8.10,8.40){\makebox(0,0)[cc]{$1$}}
\put(5.81,8.40){\makebox(0,0)[cc]{$3$}}
\put(0.04,8.63){\makebox(0,0)[cc]{$2$}}
\put(19.89,8.63){\makebox(0,0)[cc]{$2$}}
\put(27.94,8.63){\makebox(0,0)[cc]{$3$}}
\put(32.79,8.63){\makebox(0,0)[cc]{$1$}}
\put(47.78,8.81){\makebox(0,0)[cc]{$3$}}
\put(55.66,8.81){\makebox(0,0)[cc]{$1$}}
\put(60.51,8.81){\makebox(0,0)[cc]{$2$}}
\put(6.50,167.33){\makebox(0,0)[cc]{$= (1)$}}
\end{picture}
\end{center}
\caption{A representation of elements of $\BC(1)$, $\BC(2)$ and
$\BC(3)$ by equivalence classes of planar trees.\label{master}}
\end{figure}
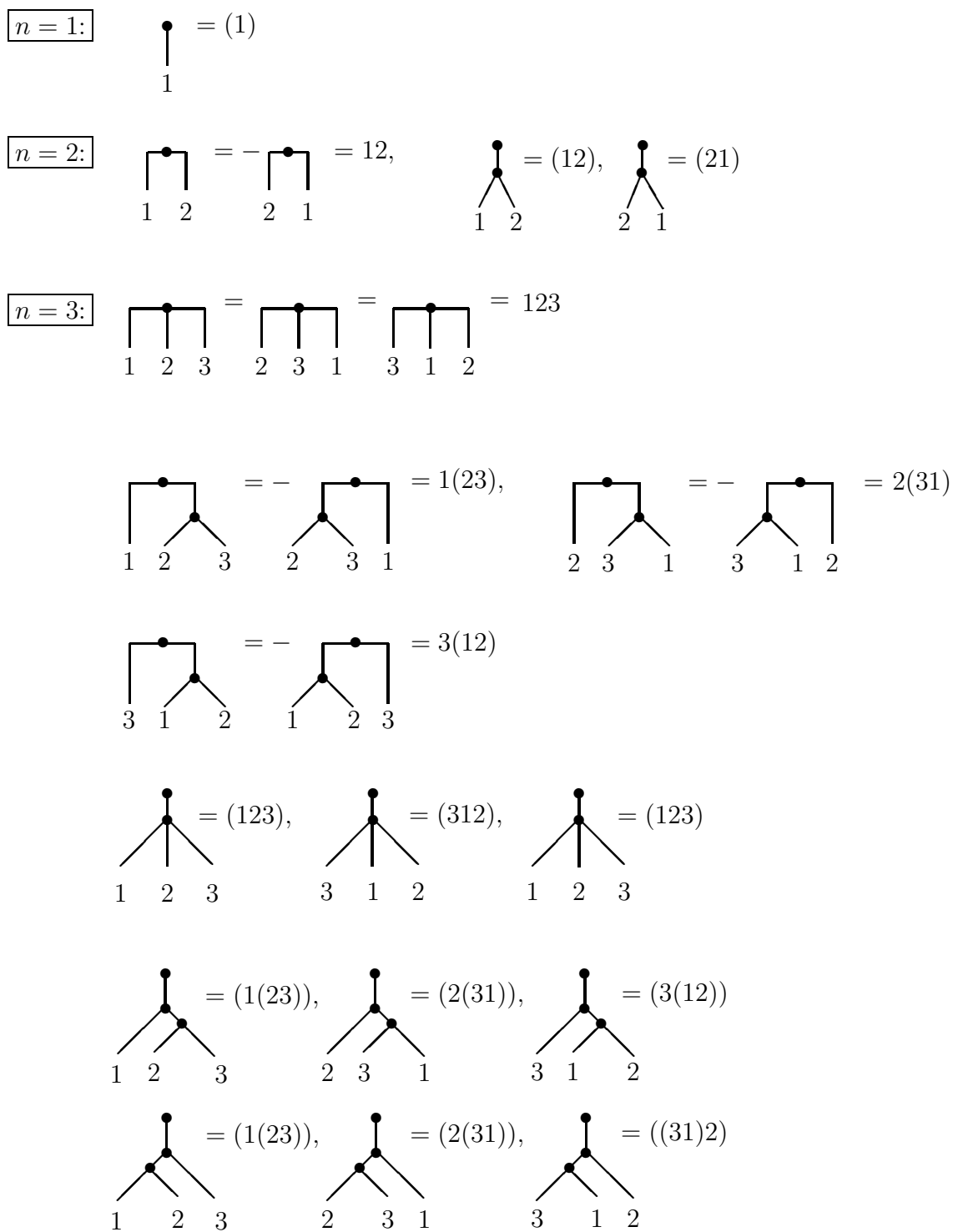

We may symbolize the element $R(\rada{T_1}{T_l}) \times
\sigma \in
J(n)$ as the tree $R(\rada{T_1}{T_l})$ with the inputs
labeled by $\rada{\sigma^{-1}(1)}{\sigma^{-1}(l)}$.
There is an almost obvious one-to-one correspondence between
these
labeled planar $n$-trees and cyclic bracketings of $n$
indeterminates from
$\barBC(n)$. This
becomes absolutely clear after
looking at Figure~\ref{master}.

To be more formal, the isomorphism $\varomega :
(\dualI{\ss\Cycl}\
\circ \desusp\dualI{\ss \Ass}) \to \M$ is defined by
$\varomega(\susp^{n-1}
\xi_n) := f_n$, where $f_n = (\rada 1n)$ is the top
$n$-dimensional
cell of the cyclohedron $W_n$ and $\xi_n$ the generator of
$\Cycl(n)$
represented by $\id_n \in \Sigma_n$. We must specify also the
orientation of $f_n$. In Observation~\ref{sdff} we
constructed points $\rada {P_1}{P_n}$ spanning a simplex
$\Delta_{f_n}\in f_n$. We orient $f_n$ coherently with the
orientation
of $\Delta_{f_n}\in f_n$ induced by the order $P_1 < \cdots <
P_n$ of
its vertices.

The cobar
differential $\partial_{\Cob}$ is given by
\begin{equation}
\label{tuti}
\partial_{\Cob}(\susp^{n-1} \xi_n) =
\cyclsum \sum_{1\leq k\leq n}\znamenko{n+k}
\nu(\susp^{n-k}\xi_{n-k+1};\alpha_k,\rada11),
\end{equation}
where $\alpha_k= \id_k\in \bk[\Sigma_k]= \Ass$ is the generator,
$k\geq 1$. We shall compare now~(\ref{tuti}) to the geometric
boundary of the top-dimensional cell $f_n$ of the cyclohedron
$W_n$.
This can be done exactly as in the proof of Theorem~\ref{resiz}
and we
leave it to the reader.\qed

The {\em proof of Theorem~\ref{nuzky}\/} is now immediate. An
$\M$-trace $T: \M\to \E_{A,W}$ is determined by a system
$\{T_n :=
t(\susp^{n-1}\xi_n): A^{\otimes n}\to W\}_{n \geq 1}$. The
axiom~(\ref{Ax}) then reflects~(\ref{tuti}).

We finish this section by the following theorem whose proof,
based on
a straightforward but involved spectral sequence argument,
we omit.

\begin{theorem}
\label{Katalogizacni}
Let $M_{\UP}$ be a module associated to a cyclic unital
quadratic
operad $\P$. Then $M_{U\P}$ is Koszul if and only if $\P$ is.
\end{theorem}

Because $\Cycl = M_{\UAss}$
(Example~\ref{kacirek}) and the operad $\Ass$ is
Koszul~\cite[Corollary~4.2.7]{ginzburg-kapranov:DMJ94},
Theorem~\ref{Katalogizacni} gives an alternative proof of
Theorem~\ref{Amphora1}.

\section{Cyclohedron as a compactification of the simplex}
\label{22}

For a compact Riemannian manifold $V$, let $C^0_n(V)
=\{(\rada{v_1}{v_n});\ v_i\not= v_j\}$ be the configuration
space
of $n$ distinct points in $V$. Axelrod and Singer constructed
in~\cite{axelrod-singer:preprint}
a compactification $C_n(V)$ of this space,
by adding to $C^0_n(V)$
the blow-ups along the diagonals. The space $C_n(V)$ is
a manifold
with corners, whose open part (= top-dimensional stratum) is
$C^0_n(V)$.

There exists a similar compactification of the {\em moduli
space\/}
$\osfF_m(n)$ of configurations of $n$
distinct points in the $m$-dimensional
Euclidean
plane ${\bf R}^m$ modulo the action of the affine group,
described by Getzler and
Jones in~\cite[\S3.2]{getzler-jones:preprint}
and denoted by $\sfF_m(n)$. The authors
of~\cite{getzler-jones:preprint} also
observed that the collection
$\sfF_m := \{\sfF_m(n)\}_{n\geq 1}$ has a
natural structure of a topological operad. In~\cite{markl:cf} we
proved
that

\begin{theorem}
\label{1725}
If $V$ is an $m$-dimensional parallelizable Riemannian
manifold, then
the collection $C(V) := \{C_n(V)\}_{n\geq 1}$ forms a right
module
over the operad $\sfF_m$ in the category of manifolds with
corners.
\end{theorem}

There is also a `framed' version of the above
theorem for manifolds which are not parallelizable, but we
will not
need it.

Take $V=S^1$. Then it is immediately seen that the space
$C_n^0(V)$
has $(n-1)!$ components indexed by cyclic orders of $n$
points on the circle. Each of these components is isomorphic to
$\oDelta^n \times S^1$,
where $\oDelta^n$
is the open $n$-dimensional simplex. It is `well-known' (see
Remark~\ref{beuo}
below) that the compactification
$C_n(S^1)$ is isomorphic to $\barW_n \times S^1$, the product
of the
symmetrized cyclohedron with the
circle~\cite[page~5249]{bott-taubes:JMP94}.

Similarly, $\osfF_1(n)$ is easily seen to have $n!$ components
indexed
by orders of the set of $n$ points on the line, each component
being
isomorphic to $\oDelta^{n-2}$. Again, it is `well-known'
that the compactification $\sfF_1(n)$ is the (symmetrized)
associahedron $\barK_n$~\cite[3.2(1)]{getzler-jones:preprint}.
This assumed,
our statement (Theorem~\ref{sirky}) about the existence of a
$\barK$-module structure on the cyclohedron follows from
Theorem~\ref{1725} applied on $C(S^1)= \barW \times S^1$
(the extra factor $S^1$ plays no r\^ole).

\begin{remark}{\rm
\label{beuo}
The Axelrod-Singer compactification $C_n(S^1)$ is a manifold
with
corners, constructed by a very explicit sequence of blow-ups. We
do
not know any `universal' characterization of this space. Thus
to prove that $C_n(S^1) \cong \barW_n \times S^1$ would require
an explicit construction of an isomorphism of two manifolds with
corners,
which is certainly not a tempting challenge. But a reflection
on the
structure of these two object `proves' the isomorphism
`beyond any
doubts', which is the opinion shared by many authors. The
same remark applies also to the isomorphism $\sfF_1(n)\cong
\barK_n$.
}\end{remark}

\begin{remark}{\rm
It follows from general properties of manifold-with-corners
that both
$K_n$ and $W_n$ are truncations of a
simplex~\cite[Proposition~6.1]{markl:cf}, but this existence
statement
says
nothing about an explicit linear convex realization of
Section~\ref{1968}.
}\end{remark}

As we observed above,
the cyclohedron $W_n$ can be viewed as the simplex $\Delta^n$,
some faces of whose were blown-up. In the rest of this
section we
show that a very natural spectral sequence related to the cobar
construction can be interpreted as an inverse process --
`deblowing-up' of the cyclohedron back to the closed simplex.

For a collection $X$ and $p\geq 1$, let $\skel Xp\subset X$
be the subcollection
defined by $\skel Xp(n) = X(n)$ for $n\leq p$ and $\skel Xp(n)
= 0$
otherwise.

Let us consider, for a comodule $N$ over an cooperad $\Q$
and for a natural $n$, the subspace
\[
F_p(n) := (\skel N{p+1} \circ \fr(\desusp \Q))(n) \subset
(N \circ \fr(\desusp \Q))(n).
\]
It is easily seen that $F_p(n)$ is $\dcob$-invariant, thus
$\{F_p(n)\}_{p\geq 0}$ is an increasing filtration
of the $n$-th piece $\Cob(N,\Q)(n)$ of the
cobar construction $\Cob(N,\Q) = (N \circ \fr(\desusp \Q),
\dcob)$.
Let $\SS(n) = (E^r_{pq}(n), d^r)$ be the corresponding spectral
sequence.
The following lemma is an easy exercise.

\begin{lemma}
\label{whoop}
The spectral sequence $\SS(n) = (E^r_{pq}(n), d^r)$
constructed above converges to $H_*(\Cob(N,\Q)(n))$.
The first term $E^1$ is described as
\[
E^1_{pq}(n)= (N(p+1)\circ H_*(\Cob(\Q)))(n)_{p+q},
\]
the space of elements of degree $p+q$ in the $n$-th piece of
the free
$H_*(\Cob(\Q))$-module on the $(p+1)$-th piece of the
collection $N$.
\end{lemma}

If the cooperad $\Q$ and the comodule $N$ are Koszul, the
spectral
sequence above collapses at the 1st term, which has a very
explicit
description. Since we did not formulate the Koszulness for
cooperads
and comodules (though the definition is an exact dual),
we suppose
from now on that
$N = \dualI{\ss M}$ and $\Q = \dualI{\ss \P}$, for a module $M$
over an
operad $\P$. We also suppose that $M$ and $\P$ are not graded,
i.e.~
that both $M(n)$ and $\P(n)$ are concentrated in degree $0$,
$n\geq 1$.

\begin{proposition}
\label{myska}
If the operad $\P$ is Koszul, then $E^1_{pq}(n)= 0$ for
$q\geq 1$,
while
\[
E^1_{p0}(n)= (\dualI{\ss M}(p+1) \circ \P^!)(n),
\]
and the spectral sequence collapses at $E^1$.
The module $M$ is Koszul if
and only if the complex
\[
0 \stackrel{d^1}{\longleftarrow} E^1_{00}(n)
\stackrel{d^1}{\longleftarrow} E^1_{10}(n)
\stackrel{d^1}{\longleftarrow} E^1_{20}(n) \longleftarrow \cdots
\]
is acyclic in positive dimensions, for all $n\geq 1$.
\end{proposition}

\noindent
{\bf Proof.}
If $\P$ is Koszul, then $H_*(\Cob(\P))=
H_0(\Cob(\P))=\P^!$, by definition. Thus, by Lemma~\ref{whoop},
$E^1_{pq}(n)= (\dualI{ \ss M}(p+1)\circ \P^!)(n)_{p+q}$ which
may be
nonzero only for $q=0$, because $\dualI{ \ss M}(p+1)$ is
concentrated
in degree $p$.
The collapsing is obvious from degree
reasons.
The second part of the statement follows immediately
from the definition of the
Koszulness of a module (Definition~\ref{telefon}).\qed

\begin{figure}[hbtp]
\begin{center}
%TexCad Options
%\grade{\off}
%\emlines{\off}
%\beziermacro{\off}
%\reduce{\on}
%\snapping{\off}
%\quality{2.00}
%\graddiff{0.01}
%\snapasp{1}
%\zoom{1.00}
\unitlength 1.50mm
\thicklines
\begin{picture}(91.66,52.08)
\put(0.50,37.00){\line(1,1){10.00}}
\put(10.72,47.41){\line(1,0){20.00}}
\put(30.72,47.08){\line(1,-1){10.00}}
\put(40.39,37.41){\line(-1,-1){10.00}}
\put(30.72,27.08){\line(-1,0){20.00}}
\put(10.86,27.21){\line(-1,1){10.00}}
\put(0.72,37.08){\makebox(0,0)[cc]{$\bullet$}}
\put(10.72,47.08){\makebox(0,0)[cc]{$\bullet$}}
\put(30.72,47.08){\makebox(0,0)[cc]{$\bullet$}}
\put(40.72,37.08){\makebox(0,0)[cc]{$\bullet$}}
\put(30.72,27.08){\makebox(0,0)[cc]{$\bullet$}}
\put(10.72,27.08){\makebox(0,0)[cc]{$\bullet$}}
\put(8.72,52.08){\makebox(0,0)[cc]{$((12)3)$}}
\put(32.72,52.08){\makebox(0,0)[cc]{$(1(23))$}}
\put(42.72,37.08){\makebox(0,0)[lc]{$1)((23)$}}
\put(32.72,22.08){\makebox(0,0)[cc]{$1))(2(3$}}
\put(8.72,22.08){\makebox(0,0)[cc]{$1)2)((3$}}
\put(-1.28,37.08){\makebox(0,0)[rc]{$(12))(3$}}
\put(20.72,52.08){\makebox(0,0)[cc]{$(123)$}}
\put(20.72,22.08){\makebox(0,0)[cc]{$1)2(3$}}
\put(7.72,42.08){\makebox(0,0)[lc]{$(12)3$}}
\put(7.72,32.08){\makebox(0,0)[lc]{$12)(3$}}
\put(33.72,42.08){\makebox(0,0)[rc]{$1(23)$}}
\put(33.72,32.08){\makebox(0,0)[rc]{$1)(23$}}
\put(20.72,37.08){\makebox(0,0)[cc]{$123$}}
\put(55.83,4.33){\line(5,6){16.67}}
\put(72.50,24.33){\line(5,-6){16.81}}
\put(89.30,4.16){\line(-1,0){33.64}}
\put(10.72,46.75){\line(1,0){20.00}}
\put(40.86,36.94){\line(-1,-1){10.00}}
\put(10.39,26.75){\line(-1,1){10.00}}
\put(72.33,24.11){\makebox(0,0)[cc]{$\bullet$}}
\put(55.66,4.11){\makebox(0,0)[cc]{$\bullet$}}
\put(89.44,4.11){\makebox(0,0)[cc]{$\bullet$}}
\put(53.00,1.89){\makebox(0,0)[cc]{$\{3\}$}}
\put(91.66,2.11){\makebox(0,0)[cc]{$\{2\}$}}
\put(72.55,27.44){\makebox(0,0)[cc]{$\{1\}$}}
\put(63.25,16.77){\makebox(0,0)[rb]{$\{13\}$}}
\put(81.44,17.00){\makebox(0,0)[lb]{$\{12\}$}}
\put(72.33,1.44){\makebox(0,0)[ct]{$\{23\}$}}
\put(72.33,12.33){\makebox(0,0)[cc]{$\{123\}$}}
\put(44.84,29.43){\line(6,-5){9.02}}
\put(44.20,28.74){\line(6,-5){9.01}}
\put(53.93,21.00){\line(-1,0){2.61}}
\put(53.93,21.00){\line(0,1){2.16}}
\end{picture}
\end{center}
\caption{$\Delta_3$ as deblowing-up of $W_3$. The faces of
$W_3$ which
are contracted by $d^0$ to a vertex
are indicated by double lines, $(123)$ is contracted to
$\set 1$,
$1)(23$ to $\set 2$ and $12)(3$ to $\set 3$.\label{deblow}}
\end{figure}
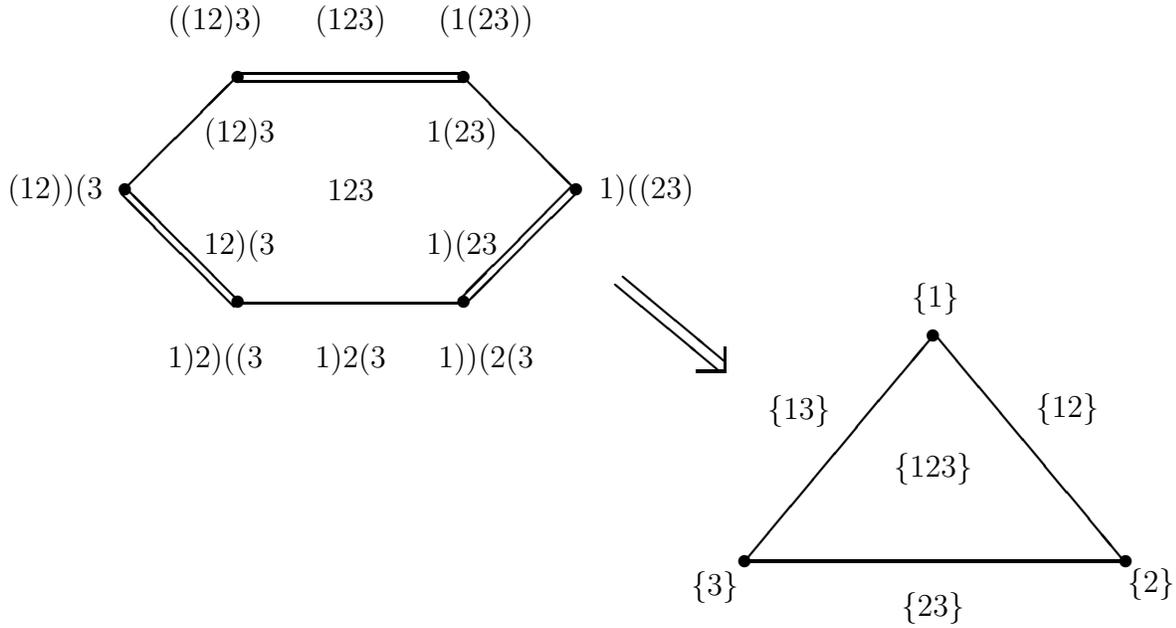

Our spectral sequence has,
for $\P = \Ass$ and $M =\Cycl$, a beautiful
geometric meaning.
The initial term $\SS^0 = (E^0_{pq},d^0)$ is the cobar
construction $\Cob(\dualI{\ss\Cycl}, \dualI{\ss\Ass})$
which is isomorphic, by
Theorem~\ref{ucpavka}, to the
cellular chain complex of the cyclohedron $\barW_n$, while
$\SS^1 =
(E^1_{pq},d^1) = (\ss\Cycl \circ \Ass,\partial)$ is isomorphic,
by Theorem~\ref{resiz}, to the cellular chain
complex of the simplex $\barDelta_n$. The passage from
$\SS^0$ to
$\SS^1$ can be interpreted as the `deblowing-up' of
the cyclohedron back
to the simplex. This process is visualized on
Figure~\ref{deblow}.

\section*{Appendix: Traces versus invariant bilinear forms.}

Let us recall the following notion
of~\cite[Definition~4.1]{getzler-kapranov:cyclic}.
If $\P$ is a cyclic operad and $A$ a $\P$-algebra, then
a bilinear
form $B: A\otimes A\to W$ with values in a vector space $W$
is called
{\em invariant\/} if, for all $n\geq 0$, the map
$B_n:\P(n)\otimes A^{\otimes (n+1)} \to W$ defined by the
formula
\begin{equation}
\label{plus}
B_n(p\ot x_0 \ot x_1 \ot \cdots \ot x_n):=
\znamenko {|x_0|\cdot |p|}B(x_0,p(\rada {x_1}{x_n})),
\end{equation}
is invariant under the action of the symmetric group
$\Sigma_{n+1}$
on $\P(n)\ot A^{\otimes (n+1)}$.

\vskip2mm
\noindent
{\bf Proposition A.1.}
{\it Let $\P$ be a cyclic operad and let $M_{\P}$ be the
associated
module introduced in Definition~\ref{Turmo}.
Let $A$ be a $\P$-algebra. Then
there exists a 1-1 correspondence between $M_{\P}$-traces on the
$\P$-algebra $A$ in the sense of Definition~\ref{el}, and
invariant bilinear forms on $A$.}

\vskip2mm
\noindent
{\bf Proof.}
Let $t : M_\P \to \E_{A,W}$ be an $M_\P$-trace. Because
$M_\P(n+1)= \P(n)$, the trace is represented by a system $\{t_n
:\P(n-1)\to \Hom{A^{\otimes n}}W\}_{n\geq 2}$ of linear maps. We
claim that $B := t_2(1): A\otimes A \to W$, where $1\in \P(1)$
is the
unit, is an invariant bilinear form.
To see it, observe that~(\ref{plus}) can be rewritten as
\[
B_n(p\otimes x_0 \otimes \cdots \otimes x_n) =
\nu_{\E_{A,W}}(t_2(1);1,p)(\rada{x_0}{x_n}),
\]
while
\[
\nu_{\E_{A,W}}(t_2(1);1,p) =
\nu_{\E_{A,W}}(\nu_{M_\P}(1;1,p)) = t_{n+1}(p \cdot \tau_n),
\]
thus
\begin{equation}
\label{xplus}
B_n(p\otimes x_0 \otimes \cdots \otimes x_n) =
t_{n+1}(p \cdot \tau_n)(\rada{x_0}{x_n})
\end{equation}
and the equivariance of $B_n$ follows from the equivariance of
$t_{n+1}$. On the other hand, if $B$ is an invariant bilinear
form,
then~(\ref{xplus}) defines a trace.\qed

%\bibliography{b}

\catcode`\@=11

\noindent
Mathematical institute of the Academy, %\hfill\break
\v Zitn\'a 25, %\hfill\break
115 67 Prague 1, %\hfill\break
The Czech Republic\hfill\break
e-mail: {\tt markl@math.cas.cz}

\end{document}